%% file: MarkovichNatalia1.tex
\RequirePackage{fix-cm}
\documentclass[smallextended]{svjour3}       
\smartqed  
\usepackage{graphicx}
\usepackage{subfigure}
\usepackage{algorithm}
\usepackage{algpseudocode}

\newcommand{\Ii}{\mathbf{1}}
\begin{document}
\title{Extremal properties of evolving networks: local dependence and heavy tails
}


\author{Natalia Markovich
}


\institute{N. Markovich \at
              V.A.Trapeznikov Institute of Control Sciences  Russian Academy of Sciences \\
              Tel.: +7(495)3348820\\
              \email{nat.markovich@gmail.com}           
}

\date{Received: date / Accepted: date}

\maketitle

\begin{abstract}
A network evolution with predicted tail and extremal indices of PageRank and the Max-Linear Model used as node influence indices in random graphs is considered.  The tail index shows a heaviness of the distribution tail. The extremal index is a measure of clustering (or local dependence) of the stochastic process. The cluster implies a set of consecutive exceedances of the process over a sufficiently high threshold. Our recent results  concerning sums and maxima of non-stationary random length sequences of regularly varying random variables are extended to random graphs. Starting with a set of connected stationary seed communities  as a hot spot and ranking them with regard to  their tail indices,  the tail and extremal indices of  new nodes that are appended  to the network  may be determined. This procedure allows us to predict a temporal network evolution in terms of tail and extremal indices. The extremal index determines  limiting distributions of a maximum  of the PageRank and the Max-Linear Model of newly attached nodes. The exposition is provided by algorithms and examples. To validate our theoretical results, our simulation and real data study concerning a linear preferential attachment as a tool for network growth are provided.
\end{abstract}
\input{introduction1.tex}
\input{preliminaries1.tex}
\section{Main results}\label{Sec3}
\input{iterations1.tex}
\input{dependence1.tex}
\input{modification1.tex}
\input{simulation1.tex}
\input{realdata1.tex}
\input{conclusions1.tex}
\input{appendix.tex}

\input{reference1.tex}
\begin{acknowledgements}
The author was supported by the Russian Science Foundation
(grant \mbox{No.\,22-21-00177}). The author would like to thank anonymous reviewers for  useful comments.
\end{acknowledgements}
\end{document}

%% file: introduction1.tex
\section{Introduction}
Extreme value theory concerning  the sums and maxima of  random  sequences  attracts a lot of interest due to numerous applications (see, 
Asmussen \& Foss, 2018; Jessen \& Mikosch, 2006; Lebedev, 2015;  Markovich \& Rodionov, 2020a; Olvera-Cravioto, 2012; Robert \& Segers, 2008; Tillier \& Wintenberger, 2018). It has progressed in recent years from finite to random lengths sequences, particularly with regard to an application in random graphs and networks (Jelenkovic \& Olvera-Cravioto, 2010; Jelenkovic \& Olvera-Cravioto, 2015; Garavaglia et al., 2020; Volkovich \& Litvak, 2010).
Tail and extremal indices of the sums and maxima of  non-stationary  sequences of regularly varying random variables (r.v.s) with finite and random lengths were obtained in Goldaeva (2013),  Markovich and Rodionov (2020a), Markovich (2021, 2022). While the tail index shows the heaviness of the distribution tail, the extremal index  is a local dependence measure that shows the cluster structure of a stationary distributed sequence.
\par
The distribution tail of a non-negative r.v. $X$ is called regularly varying $RV_{-\alpha}$ with the tail index $\alpha$ if it holds
\begin{eqnarray*}
  \overline{F}(x)&=& P\{X>x\}=x^{-\alpha}\ell(x),
\end{eqnarray*}
where the function $\ell(x)$ is  slowly varying, i.e. $\lim_{x\to\infty}\ell(tx)/\ell(x)=1$  holds for any  $t>0$.
\\
Let $X^n= \{X_i\}_{i=1}^n$ 
be a sample of 
r.v.s with cumulative distribution function (cdf) $F(x)$.
By Leadbetter et al. (1983, p. 67) 
the stationary sequence  $\{X_n\}_{n\ge 1}$ is said to have extremal index
 $\theta\in [0,1]$
 if
for each $0<\tau <\infty$ there is a sequence of real numbers $u_n=u_n(\tau)$ such that it holds
\begin{eqnarray}\label{11}&&\lim_{n\to\infty}n(1-F(u_n))=\tau, \qquad
\lim_{n\to\infty}P\{M_n\le u_n\}=e^{-\tau\theta},\end{eqnarray}
where $M_n=\max\{X_1,...,X_n\}$. Particularly, the non-stationarity of $\{X_n\}$ causes the non-existence of the extremal index.
In case the extremal index exists, it allows us to obtain a limiting distribution of $M_n$, namely, 
$P\{M_n\le u_n\} = F^{n\theta}(u_n) + o(1)$ holds. 
For independent r.v.s $\{X_1,...,X_n\}$ it holds $\theta=1$. The reciprocal of $\theta$ approximates the mean cluster size. In this sense, it measures a local dependence. Throughout the article the cluster of exceedances defines a set of consecutive observations exceeding a threshold between two consecutive non-exceedances. Since such clusters 
may cause destructive events, the extremal index plays an important theoretical and practical role.
\par
 A network evolution arises in many applications like the World Wide Web, urban transport networks, citations between scientific articles, percolation theory to site and bond percolation (Bollob\'{a}s \& Riordan, 2006; Newman, 2018), 
information spreading (Censor-Hiller \& Shachnai, 2010; Mosk-Aoyama \& Shah, 2006), economic networks of trades (da Cruz \& Lind, 2013) and many others.
 The 
 evolution is of main interest, particularly with regard to the brain neurological  networks (Bagrow \& Brockmann, 2013; McCormick \& Contreras, 2001), the infection spreading (Holme \& Litvak, 2017), and the popularity of Web pages.
The appending of new nodes and edges may be modelled by a preferential attachment (PA) 
(Ghoshal et al., 2013; Krapivsky \& Redner, 2001; Newman, 2018; Norros \& Reittu, 2006; Samorodnitsky et al., 2016; Wan et al., 2020; da Cruz \& Lind, 2013) or the attachment depending on the clustering coefficient 
(Bagrow \& Brockmann, 2013; Schroeder et al., 2022). The choice of a seed network as a  hot spot is also important for future attachments since the network may be non-homogeneous.
\par
Our objective is to obtain an evolved network with predicted tail and extremal indices of the PageRank (PR) and the Max-Linear Model (MLM) that are used as influence indices of the nodes.
We apply the results obtained in Markovich and Rodionov (2020a), Markovich (2021, 2022) to the graph enlargement and assume that 
the node PRs 
of a seed network are regularly varying distributed r.v.s (see Appendices \ref{Sec5.1} and \ref{Sec5.2} for details).
\par
We begin the attachment from a seed network consisting of stationary communities of nodes that may be connected by a few edges. The community consists of sets of nodes that are strongly connected with each other and weakly connected with nodes from other communities (Fortunato, 2010;  Leskovec et al., 2009; Mester et al., 2021). A Directed Louvain's Algorithm is a powerful tool to divide a graph into non-overlapping and weakly connected communities, Dugu\'{e}  \& Perez (2015) (see Markovich et al. (2022) for an implementation). The definition and testing of the stationarity in the graphs remain an open problem. One can determine that a graph is stationary if for all finite subsets of vertices with the same adjacency matrices the joint distributions of their in- and out-degrees are the same (Markovich et al., 2022).
\par
The communities can be ranked regarding their tail indices.
The community with the minimum tail index determines the tail index of PRs and the MLMs of newly appended nodes that have each at least one edge with this 
"dominating"  community.
Tail indices of node's in- and out-degrees for 
their  power-law
distributions were obtained in Samorodnitsky et al. (2016), Wan et al. (2020) depending on parameters of linear PA tools used for growing networks.
In Banerjee and Olvera-Cravioto (2021) the tail behavior of the power law distribution of the PR of a uniformly chosen vertex in a directed preferential attachment (DPA) graph is obtained. It is shown that this power law is heavier than the tail of the limiting in-degree distribution. The DPA has a  limited application since it assigns to each vertex
a deterministic out-degree and it produces graphs without directed cycles. The PA schemes by Samorodnitsky et al. (2016), Wan et al. (2020) used here are free from these restrictions.
To our best knowledge, the evolution of the tail and extremal indices of the PR and the MLM is considered here at the first time.
\par
A bridge between the sums and maxima  of  random  sequences of random lengths on one side and the PR and the MLM
 on the other side was given in  Markovich (2022) by finding the conditions when the tail and extremal indices of the sum and maximum in the right-hand sides of equations (\ref{6}) and (\ref{6a}) given in Appendix \ref{Sec5.1} are the same. This property was proved under 
 practically plausible assumptions (see Appendix \ref{Sec5.2}).
 In this paper this approach is extended further to random graphs.
\par
\begin{figure}[h!]
\begin{minipage}[t]{\textwidth}
\centering
 \subfigure[]{\includegraphics[width=0.45\textwidth]{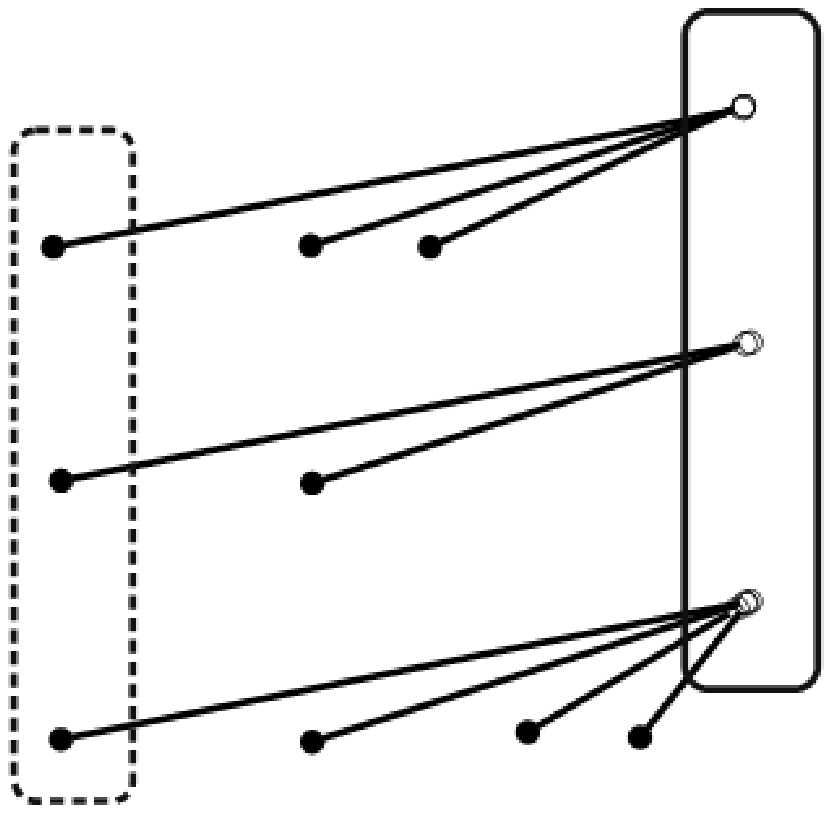} 
\label{fig:1a}}
\subfigure[]{\includegraphics[width=0.45\textwidth]{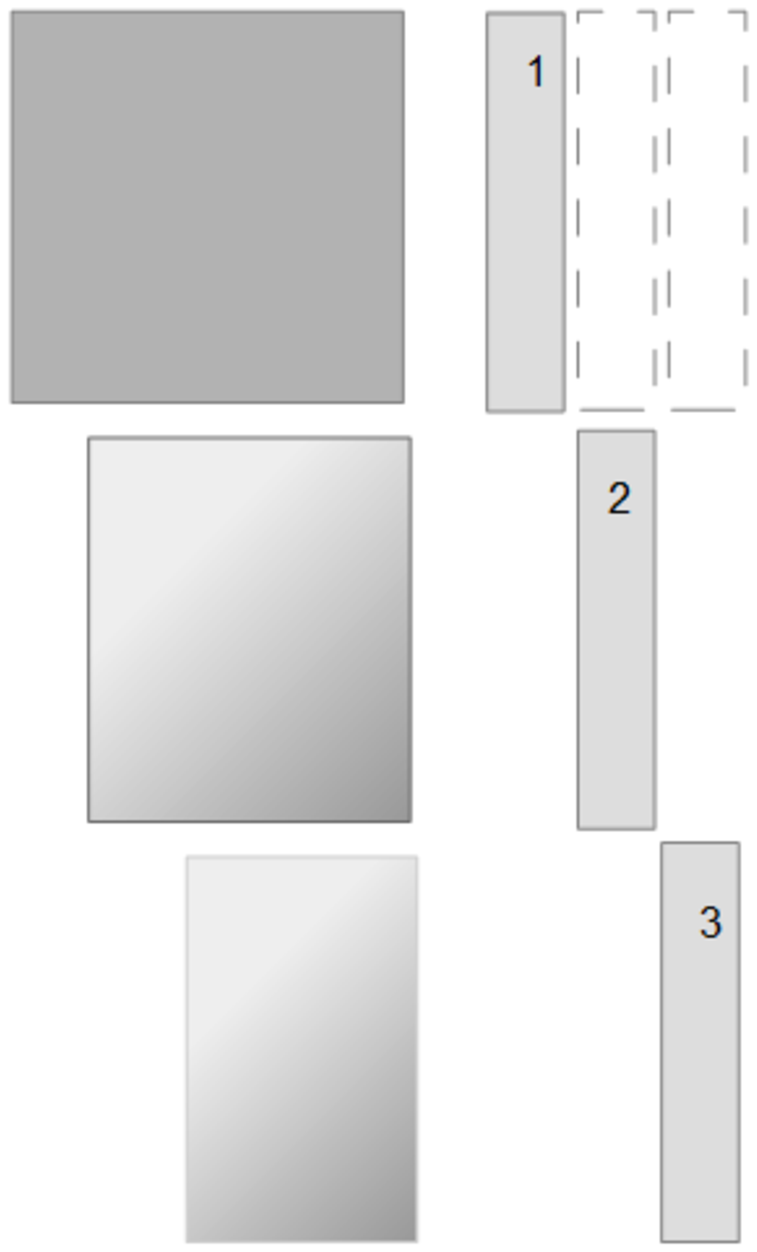}
\label{fig:1b}}
    \caption{The set of root nodes (open circles) and their nearest neighbors (filled circles) where the most heavy-tailed community is marked by a rectangle with a dotted line (Fig. \ref{fig:1a}); the creation of the enumerated columns of the next iteration matrix by the initial matrix and its submatrices by the "domino principle" (Fig. \ref{fig:1b}).
    }\label{fig:1}
    \end{minipage}
\end{figure}

We determine the mean size of the cluster of exceedances in the graph from  perspectives of extreme value theory. 
It is not  evident how to identify clusters of high level exceedances in graphs and to calculate the mean size of the clusters  as it is done for sequences of r.v.s due to an arbitrary enumeration of nodes in the graphs. Considering PRs and the MLMs of a set of root nodes as sequences of sums and maxima of PRs of their nearest neigbors with in-coming links to the roots, one may determine the extremal index of PRs and the MLMs of the roots. 
We call this value the extremal index of the (sub)graph.
It follows from Markovich (2021, 2022) that the extremal index of the subgraph (if it does exist) is determined by  the extremal index of the  most heavy-tailed 
("dominating")
sets of nodes  within the subgraph, see Fig.\ref{fig:1a}.
These sets can be obtained in the same way as communities. We assume that each 
"dominating"
community contains a stationary regularly varying distributed set of node PRs  with   a  minimum tail index.

The extremal index  depends strongly on the mutual dependence between "dominating" communities. The latest results in Markovich (2022) contains conditions to obtain the extremal index of the PRs and MLMs of the root nodes, see Appendix \ref{Sec5.2}.
\par
The paper is organized as follows. A problem description is given in Sect. \ref{Sec2_1}.
Theoretical constraints which are
important for graphs are mentioned in Sect. \ref{Sec2_3}. In Sect. \ref{Sec3} the results obtained in Markovich and Rodionov (2020a), Markovich (2021, 2022) are further developed to obtain the tail and extremal indices of PR and the MLM in 
an enlarged network started from a seed subgraph. 
Sect. \ref{Sec3} includes estimation methods and algorithms to implement the ideas to simulated and  real networks. Sect. \ref{Sec4} contains some conclusions.
The 
paper is finalized by proofs in the Appendix, where necessary theoretical results concerning  sums and maxima of random length sequences of regularly varying distributed r.v.s as well as methods relating to graphs are 
recalled.

%% file: preliminaries1.tex
\section{Preliminaries}\label{Sec2}
\subsection{Problem Description}\label{Sec2_1}
As in Markovich and Rodionov (2020a), Markovich (2021, 2022) we focus on  a doubly-indexed array $\{Y_{n,i}: n,i\ge 1\}$  of nonnegative r.v.s in which the "row index" $n$ corresponds to time, and the "column index" $i$  enumerates the series. The length $N_n$ of "row" sequences $\{Y_{n,i}: i\ge 1\}$ for each $n$ is generally random. Namely, $\{N_n: n\ge 1\}$ is a sequence of non-negative integer-valued r.v.s. For each $i$ the "column" sequence $\{Y_{n,i}: n\ge 1\}$ is assumed to be  strict-sense stationary with extremal index $\theta_i$ having a regularly varying distribution tail 
\begin{eqnarray}\label{11a} P\{Y_{n,i}>x\}&=&\ell_i(x)x^{-k_i}\end{eqnarray}
with tail index $k_i>0$ and  a slowly varying function $\ell_i(x)$. There are no assumptions on the dependence structure in $i$. Following  Markovich and Rodionov (2020a), Markovich (2022) we consider the weighted sums and maxima
\begin{eqnarray}\label{3}
&&Y_{n}^*(z, N_n) =\max(z_1Y_{n,1},...,z_{N_n}Y_{n,N_n}),\nonumber
\\
&&
Y_{n}(z, N_n)=z_1Y_{n,1}+...+z_{N_n}Y_{n,N_n}
\end{eqnarray}
for positive constants 
$\{z_i\}$, $z_i>0$, $i=1,2,...$.
\par
Let $G_n=(V_n, E_n)$ be a directed graph with a set of vertices $V_n=\{1,...,n\}$, and a set of directed edges $E_n$.
The $Y_{n}^*(z, N_n)$ and $Y_{n}(z, N_n)$ may be interpreted  as sums and maxima at the right-hand side of (\ref{6}) and (\ref{6a}) (see Appendix \ref{Sec5.1}). Each sequence $z_1Y_{n,1},...,z_{N_n}Y_{n,N_n}$, $n\ge 1$ represents the weighted influence indices of nodes in the one-link neighborhood from the root node $n$. These neighbor nodes are marked by filled circles in Fig. \ref{fig:1a}. $N_n$ denotes an in-degree of the root node $n$ that is the number of its nearest neighbors with in-coming links to the root.
\par
 In Markovich (2021) $A_jR_j$,
$j\in\{1,..., N_i\}$ in (\ref{6}) and (\ref{6a}) is denoted as $z_jY_{i,j}$ with $z_j=c$.  By the definition of Google's PR it follows that $A_j=c/D_j$ holds, where $D_j$ is the out-degree of the node $j$ and $c>0$ is  a damping factor, the only parameter of the Personalized PR (Volkovich \& Litvak, 2010). By Lemma A.1 (iii) in Volkovich and Litvak (2010) $Y_{i,j}=R_j/D_j$ has the same tail index as $R_j$ since $R_j$ and $D_j$ are assumed to be mutually independent and $E(1/D_j)<1$ holds. 
Hence, the tail of $Y_{i,j}$ is dominated by the tail of $R_j$.
One can rewrite the right-hand sides of (\ref{6}) and (\ref{6a}) as
\begin{eqnarray}\label{9a}Y_i(c, N_i)&=&c\sum_{j=1}^{N_i}Y_{i,j}+Q_i,
~~~~Y^*_i(c, N_i)=c\bigvee_{j=1}^{N_i}Y_{i,j}\vee Q_i,~~ i\in\{1,...,n\}.
\end{eqnarray}
In Markovich (2022) the conditions were found when $Y_i(c, N_i)$ and $Y^*_i(c, N_i)$ have the same tail and extremal indices (see Appendix \ref{Sec5.2}). 
Each node in a random network is considered as a root of some directed graph of its followers which may contain cycles. As in Markovich (2021) we consider graph communities as the "column" series.
\par
We aim to get the tail and extremal indices of an evolved graph starting from  a seed set of nodes with 
known tail and extremal indices.
\subsection{Important constraints}\label{Sec2_3}
Let us mention the constraints of Theorem \ref{T3} 
(Appendix \ref{Sec5.2}) that are important for graphs.
\begin{enumerate}
\item The stationarity of the node's in-degrees $\{N_n\}$ is not assumed, but $N_n'$s have the same  tail index and  their distribution tail has to be lighter than the tail of the node PRs.
\item The in-degree $\{N_i\}$ of the $i$th node and the PRs $\{Y_{i,j}\}$ of the $j$th nodes that link to node $i$ in (\ref{9a}) are independent. 
\item
The mutual pair-wise dependence between elements of the stationary $d$ "column" sequences  with minimum tail index has to be the same. Otherwise, the sequences of sums and maxima over $d$ "row" elements corresponding to these columns are non-stationary and thus, the extremal indices of such sequences do not exist.
\item For each row at least one element corresponding to the "column" sequences with the minimum tail index has to be non-zero. Since the tail index has to be estimated one may deal with a unique community  with a minimum tail index in a graph.
\item
Elements of the "column" sequences with non-minimum tail index may be arbitrarily dependent. They may have different tail indices larger than the minimum tail index and hence, these  "column" sequences may be non-stationary and have no extremal indices.
\item
In terms of graphs, the communities may be considered as "column" sequences 
and only the $d$ 
"dominating"
communities are required to be stationary distributed. Note that d is generally a random variable.
\end{enumerate}

%% file: iterations1.tex
\subsection{Iterations}\label{Sec3_1}
We focus on a directed graph $G_n$ with $n$ nodes.
Theorem \ref{T3} (Appendix \ref{Sec5.2}) relates 
to a single iteration by ranks of 
the
one-link neighbors of root nodes. It states that 
 the tail index of PRs and the MLMs of the roots is determined by 
the 
tail index of their 
most heavy-tailed nearest neighbors.
PR and the MLM of the root nodes  have the same tail index and in some cases the same extremal index.
\\
The connection between the PR of a node and the solution of (\ref{6}) is proved by convergence in distribution of the $m$th iteration  
\begin{eqnarray*}\label{1a}R^{(m)}&=&\sum_{j=1}^{N}A_{j}R_j^{(m-1)}+Q, ~~m\ge 1,\end{eqnarray*}
corresponding to the Galton-Watson tree 
to $R$ as $m\to\infty$ starting from an initial distribution $R^{(0)}$. The r.v.s $\{R_j^{(m-1)}\}$ are assumed to be independent identically distributed (iid) copies of $R^{(m-1)}$ (Jelenkovic \& Olvera-Cravioto, 2010; Volkovich \& Litvak, 2010).
Theorem 3.2 in Volkovich \& Litvak (2010)
 states that the tail behavior of $R^{(m)}$ is determined by the asymptotic of the r.v. with the heaviest tail among $N$ and $Q$.
 The initial distribution of $R^{(0)}$ is assumed to have a lighter tail than $N$ or $Q$ which both are regularly varying r.v.s, i.e. the iterations may start with $R^{(0)}\equiv 1$.
 \\
  The tail behavior of $R^{(m)}$ is proved in Theorem \ref{T3} (Appendix \ref{Sec5.2}) omitting the independence Assumptions A (Appendix \ref{Sec5.1}) and assuming that the r.v.s $\{A_{j}R_j^{(m-1)}\}$ are non-stationary regularly varying distributed and a random number of the most heavy-tailed r.v.s $\{A_{j}R_j^{(m-1)}\}$ are independent or weakly dependent (conditions (A1) or (A2) in Appendix \ref{Sec5.2}). The statement of Theorem \ref{T3} (Appendix \ref{Sec5.2}) is similar to Proposition 3.1 in Volkovich and Litvak (2010),
  where 
 $R^{(0)}$ is assumed to be a regularly varying r.v. with tail index $\alpha_R>0$. If $P\{N>x\}=o(P\{R^{(0)}>x\})$ and $P\{Q>x\}=o(P\{R^{(0)}>x\})$, then 
 $P\{R^{(m)}>x\}\sim C_R^{(m)}P\{R^{(0)}>x\}$ for all $m\ge 1$ as $x\to\infty$  is stated. $C_R^{(m)}$ is a constant depending on $m$, $\alpha_R$ and $E(N)$.
 The distribution of $R^{(\infty)}$, the unique nontrivial solution of (\ref{6}), does not depend on the distribution of $R^{(0)}$ assuming that $E(R^{(0)})= 1$, $E(A)=(1-E(Q))/E(N)<1$ and Assumptions A hold (Volkovich \& Litvak, 2010, Theorem 3.1).
\\
The convergence of the maximum recursion
 \begin{eqnarray*}
 R^{(m)}&=&\left(\bigvee_{j=1}^{N}A_{j}R_j^{(m-1)}\right)\vee Q,
 \end{eqnarray*}
 to the solution $R$ of (\ref{6a}) under the Assumptions A as $m\to\infty$
  and provided that the initial values corresponding to leafs 
 possess a moment condition  is derived in Jelenkovic and Olvera-Cravioto (2015).
 $\{R_j^{(m-1)}\}$ are independent copies of $R^{(m-1)}$ corresponding to the tree starting with an individual node $j$ in the first generation and ending at the $m$th generation.
\unskip
\subsection{Evolution from a 
seed subgraph (the "domino principle")}\label{Sec3_2}
Let us describe the evolution process of the network by means of changing of matrices corresponding to the doubly-indexed array $\{Y_{n,i}\}$ determined in Sect. \ref{Sec2_1}. Then $m$ below is connected with the time.
Using the notations of Sect. \ref{Sec2_1} we  deal with recursions
\begin{eqnarray}\label{1b} Y_{i,j}^{(m)}&=&c\sum_{s=j}^{N_i}Y_{i,s}^{(m-1)}+Q_i,
\end{eqnarray}
\begin{eqnarray}\label{1e} X_{i,j}^{(m)}&=&\left(c\bigvee_{s=j}^{N_i}X_{i,s}^{(m-1)}\right)\vee Q_i,~~~\{X_{i,j}^{(0)}\}\equiv \{Y_{i,j}^{(0)}\},
\end{eqnarray} $m, i, j\ge 1$.
Let us consider 
matrices
related to the scheme of series $\{Y^{(0)}_{n,i}: n,i\ge 1\}$ and corresponding tail and extremal indices $( k^{(0)}_i, \theta^{(0)}_i)$:
 \begin{equation}\label{4c}A^{(0)}=\left(
    \begin{array}{cccccc}
      Y^{(0)}_{1,1} & Y^{(0)}_{1,2} & Y^{(0)}_{1,3} &...  & 0 & Q_1 \\
      Y^{(0)}_{2,1} & 0 & Y^{(0)}_{2,3} & ...  & Y^{(0)}_{2,N_2} & Q_2 \\
      ...& ... & ... & ... & ... & ... \\
      Y^{(0)}_{n,1} & Y^{(0)}_{n,2} & Y^{(0)}_{n,3} & ...  &  Y^{(0)}_{n,N_n} & Q_n \\
    \end{array}
  \right),
\end{equation}
\begin{eqnarray*}\label{4d}
&&\qquad ~~~\left(
    \begin{array}{cccccc}
      k^{(0)}_1 & k^{(0)}_2 & k^{(0)}_{3} & ... &  k^{(0)}_{N}  & ~k^{(0)}_{N+1} \\
      \theta^{(0)}_1 & \theta^{(0)}_2 & \theta^{(0)}_{3} & ... &  ~\theta^{(0)}_{N} & 1 \\
          \end{array}
  \right).
\end{eqnarray*}
 $\{Q_i\}$ is a sequence of iid r.v.s and thus, its extremal index is equal to $1$. 
 We start the evolution with a seed network that can be represented by a matrix $A^{(0)}$.
 Network communities may be interpreted as columns of $A^{(0)}$.
  A zero $s$th element in the $i$th row $Y^{(0)}_{i,s}$, $s\ge 1$ 
  of $A^{(0)}$ means that the $i$th root node has no followers in the $s$th community or there is no link between them. For instance, if  a row corresponds to a set of papers citing a book, then zero implies that the book is not cited by a paper from the corresponding community.
 \par
 Without loss of generality we assume that $k_1^{(0)}\le k_2^{(0)}\le k_3^{(0)}\le ...$ holds. We may assume that the first
  $d_1^{(0)}$ columns of $A^{(0)}$ are the most heavy-tailed distributed with the minimum tail index $k_1^{(0)}$, the next $d_2^{(0)}$ columns have the second minimum tail index $k_2^{(0)}$ etc.
 Generally, $\{d_i^{(0)}\}$ are r.v.s.
 \par We further consider the evolution iterations as follows: starting with this 
 seed
network, each time appending a set of nodes  which transform the matrix $A^{(m)}$. The nodes are appended by some attachment tool.
The $j$th column $\{Y_{i,j}^{(m)}\}_{i\ge 1}$ (or $\{X_{i,j}^{(m)}\}_{i\ge 1}$) of the matrix $A^{(m)}$  is defined by  (\ref{1b}) (or (\ref{1e})) using the submatrix $\{Y^{(m-1)}_{n,i}: n\ge 1,i\ge j\}$ (or $\{X_{n,i}^{(m-1)}\}: n\ge 1,i\ge j\}$) of the matrix $A^{(m-1)}$.
It is assumed that the newly attached nodes that built the matrix $A^{(m)}$ have edges only with nodes of the subgraph corresponding to the matrix $A^{(m-1)}$. The evolution looks like the "domino principle" as $j$ in (\ref{1b}) or (\ref{1e}) increases (Fig. \ref{fig:1b}).
What would be the extremal and tail indices of the "column" sequences for the next iteration matrices $\{A^{(m)}\}$, $m\ge 1$?
\par Returning to the citation example, the "domino principle" means that some books may be cited by representative papers from all considered communities including the most influential ones or by papers from parts of the communities which are not so distinguish as the first ones. For instance,
$\{Y_{i,1}^{(1)}\}$ is calculated by all elements of the
$i$th row of $A^{(0)}$, as far as $\{Y_{i,2}^{(1)}\}$ by the same row elements starting from the second one etc.
\par
 In Item (ii) of the next theorem, conditions $(A1)-(A4)$ (Appendix \ref{Sec5.2}) with regard to elements of matrix $A^{(0)}$ are assumed where $d$ is replaced by $\lfloor d_n-1\rfloor$, where $d_n = \min(C, l_n)$, 
$C>1$ and $l_n$ satisfies (\ref{27}), (\ref{chi}). 
 \begin{theorem}\label{Prop1}
 Let the conditions of Theorem  \ref{T3} (Appendix \ref{Sec5.2}) with regard to $\{Y^{(0)}_{n,i}: n,
  i\ge 1\}$
 be fulfilled and  at least one element  in each row of $\{Y^{(0)}_{n,i}\}_{i\ge 1}$
 among the columns with tail index $k_j^{(0)}$, i.e.
 $\{Y^{(0)}_{n,i}: n\ge 1,d_{j-1}^{(0)}+1\le i\le d_{j-1}^{(0)}+d_j^{(0)}\}$, $d_{0}^{(0)}=0$, 
 $j\in\{1,2,...\}$ for each $n$
 be non-zero. Assume that $d_j^{(0)}$ and $\{Y_{n,i}^{(0)}\}$  are independent for each $j\ge 1$.
 \begin{enumerate}
 \item[(i)]  If $d_j^{(0)}=1$, $j\in\{1,2,...\}$ almost surely (a.s.),
 then $\{Y_{i,j}^{(m)}\}_{i\ge 1}$ and $\{X_{i,j}^{(m)}\}_{i\ge 1}$ calculated by (\ref{1b}) and (\ref{1e}) have
 the same tail index $k_j^{(0)}$
and  the same extremal index 
  $\theta_j^{(0)}$  for any $m\ge 1$;
 \item[(ii)] Let $\{d_j^{(0)}\}$, $j\in\{1,2,...\}$ be  bounded discrete r.v.s  such that
$1<d_j^{(0)}<d_n = \min(C, l_n)$, 
$C>1$ holds.
\begin{enumerate}
\item
    If (A1) or (A2) for any $d_j^{(0)}\in\{2,3,...,\lfloor d_n-1\rfloor\}$, $j\in\{1,2,...\}$
holds and
$N_n$ and $\{Y_{n,i}^{(0)}\}$ are  independent, then $\{Y_{i,j}^{(m)}\}_{i\ge 1}$ and $\{X_{i,j}^{(m)}\}_{i\ge 1}$ have the  tail index $k_j^{(0)}$ for any $j\ge 1$ and $m\ge 1$.
\item If (A4) where in (\ref{11b}) $d_j^{(0)}$, $j\in\{1,2,...\}$  is replaced by $\lfloor d_n-1\rfloor$
holds, then $\{X_{i,s}^{(m)}\}_{i\ge 1}$  has the  extremal index $\theta_s^{(0)}$ for any $s\ge 1$ and $m\ge 1$. If, in addition, (A1) (or (A2)) for any 
$d_j^{(0)}\in\{2,3,...,\lfloor d_n-1\rfloor\}$ holds, then $\{Y_{i,s}^{(m)}\}_{i\ge 1}$  has the same extremal index  for any $s\ge 1$ and $m\ge 1$.
\end{enumerate}
\end{enumerate}
 \end{theorem}
 Theorem  \ref{Prop1} is valid assuming that each row of  $A^{(0)}$  contains at least one non-zero element in the most heavy-tailed columns, i.e. the columns with the minimum tail index. 
 Otherwise, the sums and maxima over rows may be non-stationary distributed with different tail indices.
 The elements of the matrices $\{A^{(m)}\}$, $m\ge 1$ of the next iterations are non-zero by their construction as row sums or maxima. Theorem \ref{Prop1} states that the limit distributions of 
 recursions (\ref{1b}) and (\ref{1e}) 
 depend on 
 distributions of columns $\{Y^{(0)}_{n,i}: n\ge 1,i\ge 1\}$.
  \par
  Due to the complex nature of real-world networks the "non-zero assumption" may be rather restrictive when the number $d$ of the most heavy-tailed columns is small.
 To overcome the problem, one can permute the rows of $A^{(0)}$ to have blocks of rows with non-zero elements at least in one of the most heavy-tailed  columns in the block.
\begin{example}\label{Exam1}
An example of such permutation is given by matrix $A_{*}^{(0)}$:
\begin{eqnarray*}
&&A_{*}^{(0)}=\left(
                        \begin{array}{cccccc}
                          Y^{(0)}_{1,1} & Y^{(0)}_{1,2} & Y^{(0)}_{1,3}&... & 0 & Q_1 \\
                          Y^{(0)}_{2,1} & 0 & Y^{(0)}_{2,3} &...& Y^{(0)}_{2,N_2} & Q_2 \\
                          ... & ... & ... & ... & ... \\
                          Y^{(0)}_{n_1,1} &  Y^{(0)}_{n_1,2} &  0 & ... & 0 & Q_{n_1} \\
                         \hline
                          0 & Y^{(0)}_{n_1+1,2} & Y^{(0)}_{n_1+1,3} & ... & 0 & Q_{n_1+1}\\
                          0 & ... & ... & ... & ... & ...\\
                          0 &  Y^{(0)}_{n_2,2} & Y^{(0)}_{n_2,3} & ... 
                          & Y^{(0)}_{n_2,N_{n_2}} & Q_{n_2}\\
                          \cline{2-6}
                          0 & 0 & Y^{(0)}_{n_2+1,3} & ... & 0 & Q_{n_2+1}\\
                          ... & ... & ... & ... & ... & ...\\
                          0 &  0 & Y^{(0)}_{n_3,3} & ... 
                          & Y^{(0)}_{n_3,N_{n_3}} & Q_{n_3}\\
                          ... & ... & ... & ... & ... & ...\\
                        \end{array}
                      \right),
                      \end{eqnarray*}
                      where
             \begin{eqnarray*}      &&  A_{*}^{(1)}=\left(
                        \begin{array}{ccccc}
                          Y^{(1)}_{1,1} & Y^{(1)}_{1,2} & Y^{(1)}_{1,3}&... &Q_1 \\
                          Y^{(1)}_{2,1} & Y^{(1)}_{2,2} & Y^{(1)}_{2,3} &...& Q_2\\
                          ... & ... & ... & ...  \\
                          Y^{(1)}_{n_1,1} &  Y^{(1)}_{n_1,2} &   Y^{(1)}_{n_1,3} & ... &Q_{n_1} \\
                          \hline
                          0 & Y^{(1)}_{n_1+1,2} & Y^{(1)}_{n_1+1,3} & ... & Q_{n_1+1}\\
                          0 & ... & ... & ... &... \\
                          0 & Y^{(1)}_{n_2,2} & Y^{(1)}_{n_2,3} & ... & Q_{n_2}\\
                          \cline{2-5}
                          0 & ... & ... & ... \\
                          0 & 0 & Y^{(1)}_{n_3,3} & ... & Q_{n_3}\\
                          ... & ... & ... & ... & \\
                        \end{array}
                      \right),
\end{eqnarray*}
is the matrix of the next iteration,
$n_i=o(n)$, $i\in\{1,2,...\}$ as $n\to\infty$. The "column" series are assumed to be stationary regularly varying distributed and their tail indices  $k^{(0)}_1<k^{(0)}_2<k^{(0)}_3<...$  are increasing. We assume for simplicity that there is a unique column with a minimum tail index in each block. Zeroes in  $A_{*}^{(0)}$ imply that the corresponding node has no out-going links to other nodes, particularly, to newly appearing ones. By Theorem 4 in Markovich \& Rodionov (2020)
the sums and maxima over each of the first $n_1$ rows in the first block have the tail index $k^{(0)}_1$, and over rows with numbers $\{n_1+1,...,n_2\}$ in the second block  the tail index $k^{(0)}_2$ etc. as $n=\sum_i n_i\to\infty$. Elements of the $i$th block of  $A_{*}^{(0)}$ may contain zeroes apart of the $i$th 
"dominating"
column, $i\in\{1,2,...\}$, but they do not impact the tail index of the $i$th column of  $A_{*}^{(1)}$ since they have non-minimum tail indices.
In terms of citations it means that a set of books cited by the most influential community inherits its tail and extremal indices.
The part of the $i$th column within the $i$th block of $A_{*}^{(1)}$ corresponding to rows $\{n_{i-1}+1,...,n_{i}\}$, $i\ge 1$ with $n_0=0$   
has the tail index $k^{(0)}_i$ 
as $n\to\infty$. The non-zero rest of the $i$th column in $A_{*}^{(1)}$ corresponding to rows $1,2,...,\sum_{j=1}^{i-1}n_j$ may be non-stationary distributed by Remark \ref{Rem1} (Appendix \ref{Sec5.2}).
\\
 In case 
 the first $d>1$ columns of $A_{*}^{(0)}$ have the minimum tail index, e.g., $k^{(0)}_1=...=k^{(0)}_d<k^{(0)}_{d+1}<...$, the first rows in $A_{*}^{(0)}$ providing non-zero elements in the first $d$ column are selected into one block.
 $A_{*}^{(1)}$ can be partitioned  into blocks with non-zero elements and used further for the next iterations.
\end{example}
\begin{example}\label{Exam2} Let us consider a citation network (Krapivsky \& Redner, 2001; Newman, 2018). Each newly appearing paper cites a list of selected papers. In terms of networks the number of cited papers implies  the out-degree of the new node (the paper), 
Fig. \ref{fig:2}. Since papers are published continuously, the citation to some paper can appear in  papers published later.  PRs of the newly appearing
papers in a unit time build a row in matrix $A^{(0)}_*$.
The sum and maxima of these PRs  over the row provide the PR and the MLM of the  
paper published earlier.
At time $n$ one can build $n$ rows of the matrix $A^{(0)}$. A column of the matrix $A^{(1)}_*$ of influence measures of previously published cited papers is built by $A^{(0)}_*$. One has to detect the columns of $A^{(0)}_*$ with the minimum tail index. It is assumed that each column of 
newly appearing
papers in $A^{(0)}_*$ with a minimum tail index is stationary distributed. According to Corollary \ref{Cor2} (Appendix \ref{Sec5.2}) other columns may be non-stationary distributed  and include r.v.s with different tail indices. 
This allows us to obtain the tail and extremal indices of the sequence of previously published cited
papers, i.e. 
the first "column" of $A^{(1)}_*$. To predict the tail and extremal indices of other columns, the rest of the columns of $A^{(0)}_*$ has to be stationary distributed.
\par
The number of "dominating" communities (the  columns of $A^{(0)}_*$ with the minimum tail index) is plausible to be a bounded r.v.. There exists a "top" community among the latter ones such that its maximum PR is the largest. Then the MLMs   of  cited papers published earlier have the extremal and tail indices of
 the "top" community. If the communities with  the minimum tail index are independent, then the PRs   of  cited papers have the same extremal and tail indices. This follows by Markovich (2022), see Appendix \ref{Sec5.2}. It implies that the citation by the "top dominating" community of newly appearing papers impacts the influence of the cited papers.
\par
The network can be re-directed such that the rows of the matrix $A^{(0)}_*$ contain previously published cited papers and the columns of $A^{(1)}_*$ include the PR and the MLM of the  newly appearing papers, see Fig. \ref{fig:2}.
\par
The evolution of the citation network is determined here by the dynamics of the tail and extremal indices which may predict the long-term citation impact of a set of publications.
The mutual pair-wise dependence between r.v.s in the most heavy-tailed columns has to be the same and a part of each row  related to the latter columns has to contain at least one non-zero element. Otherwise, sums $Y_{n}(z,N_n)$ and maxima $Y_{n}^*(z,N_n)$ are not stationary distributed, and the extremal index does not exist, see Example 1 in Markovich (2021).
\begin{figure}[h!]
   \begin{center}
    \includegraphics[width=0.44\textwidth]{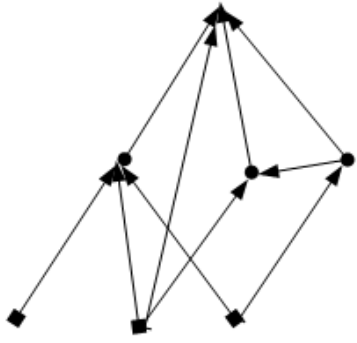}
    \end{center}
    \caption{The scheme of a citation network with newly appearing papers marked by squares and cited papers marked by black circles; cited papers  refer to papers published earlier. The scheme can be re-directed to the "ancestor" nodes taking a root of the tree as a newly appearing paper.
    }\label{fig:2}
\end{figure}
\end{example}
 \begin{remark}
Theorem  \ref{Prop1} is based on the 
assumption that  
 node PRs of the initial seed graph
 corresponding to  (\ref{4c}) are regularly varying distributed. 
 This 
 is  supported
 by empirical studies where
 the in-degree and PR of the Web are shown to follow a power law with the same exponent (Litvak et al., 2007; Pandurangan et al., 2002; Volkovich \& Litvak, 2010).
 Marginal degree power laws for growing random networks were established in Krapivsky and Redner  (2001). 
 In Samorodnitsky et al. (2016), 
 the joint distribution of in-
and out- degrees in networks growing by a linear PA model is proved to have 
regularly varying tails.
 To obtain regularly varying PRs, 
 one can consider nodes with PRs in  (\ref{4c}) as roots of classical branching trees whose leafs have PRs equal to $1$ as in Chen (2014), Volkovich and Litvak (2010).
\end{remark}

%% file: dependence1.tex
\subsection{Selection of the seed matrix $A^{(0)}$}
  \par Next, we  study how to define elements of the matrix $A^{(0)}$ in (\ref{4c}) and to find the tail and extremal indices of its columns. 
  One can collect nodes by an attachment tool and calculate their PRs. Then, one can partition the obtained seed graph into 
  communities 
  (Clauset et al., 2004; Coscia, 2011; Dugu\'{e} \& Perez, 2015). The node's PRs of the communities can be considered as the columns of $A^{(0)}$.
  \par
  For a given  personalization vector $Q_i = 1/n, 1 \leq i \leq n = \mid V_n\mid,$
the scale-free PR
$R_i^{(n)} = n R_i$
of a node $i$  can be computed iteratively (Chen et al., 2014)  by
\begin{eqnarray}
\label{8}
&&\widehat{R}_i^{(n,0)} = 1,
~~\widehat{R}_i^{(n,k)} = \sum_{j\to i}\frac{c}{D_j}\widehat{R}_j^{(n,k-1)} + (1-c),~~k>0,
\end{eqnarray}
until the difference between two consecutive iterations will be small enough.
Here, $j\to i$ implies that node $j$
links to node $i$, i.e. $(j, i) \in E_n$.
\\
The mean field analysis is based on an aggregation of Web pages in classes according to pairs $k=(k_{in}, k_{out})$ of their in- and out-degrees and using  averages of PRs within each $k$-degree class to calculate the PR (Fortunato et al., 2011).
\\
The estimation of the tail index does not require an enumeration of elements of the sample. Thus, it can be estimated by one of the nonparametric estimators. The extremal index determines roughly the inverse of the expected cluster length. Its estimation may depend on the node enumeration.
We propose a specific intervals estimator in Section \ref{ModInterEst} that allows us to avoid the node enumeration and select only a threshold $u$ as one parameter.

\subsection{Empirical estimation}\label{Sec3_2}
\par The following problems have to be solved regarding the estimation of the model parameters: (i) stationarity testing  of communities; (ii) detection of pair-wise 
dependence between the elements of the most heavy-tailed communities that has to be the same. 
Stationarity of a community is required for the tail and extremal indices to exist according to their definitions. Condition (ii) has to be checked since the sequences of row sums and maxima  have to be stationary distributed, otherwise their extremal index does not exist. Testing (ii) is limited due to the complexity of real-world networks. If the most heavy-tailed community is unique, then (ii) is naturally omitted. 
\par
The well-known nonparametric estimators of the extremal index for random sequences like blocks and runs estimators require usually the choice of a threshold $u$ and/or
a declustering parameter, e.g. the block size (Beirlant et al., 2004).
The intervals estimator  (Ferro \& Segers, 2003), estimators introduced in Robert (2009) and the $K-$gaps estimator (S\"{u}veges \& Davison, 2010)) are threshold-based, i.e. they require a choice of $u$ as a single parameter. Since nodes belonging to a graph community cannot be enumerated, the estimators used for random sequences require a modification.
\par The tail index of a community of nodes may be estimated by one of the nonparametric  methods based on the upper order statistics of the sample which is more appropriate for dependent data, e.g.
 the ratio estimator in Novak (2002), Resnick and St\v{a}ric\v{a} (1999) or the SRCEN estimator in McElroy and  Politis (2007). The Hill's estimator
 \begin{eqnarray}\label{7}\widehat{\alpha}(n,k)&=&\left(\frac{1}{k}\sum_{j=1}^k\log(\frac{X_{(n-i+1)}}{X_{(n-k)}})\right)^{-1}\end{eqnarray}
 based on the $k$ upper order statistics $X_{(1)}\le X_{(2)}\le...\le X_{(n)}$
  still works in practice, despite the network data may not be iid (Wang and Resnick, 2019, 2020; Wan et al., 2020). The value $k$
 can be selected by minimizing the bootstrap mean squared error (MSE)  (Markovich, 2007; Markovich et al., 2017) or by minimizing the Kolmogorov-Smirnov distance (Clauset et al., 2009; Drees et al., 2020; Wan et al., 2020).
\par
We aim to predict 
the tail and extremal indices  of the PR and the MLM of  newly attached nodes of an evolving graph.
The main restriction is that each 
new node is considered as a root of one-link graph of its neighbors  
and must have at least one neighbor with the minimum tail index of its PR. We propose first to partition a seed network into stationary distributed communities, e.g.
by maximizing the modularity (Clauset et al., 2004;  Newman, 2018). This can be done efficiently by a directed Louvain Method (Dugu\'{e} \& Perez, 2015).
 The stationarity of communities can be roughly checked by the mean excess function. 
 Our idea is to find a set of "dominating" communities in the graph with a minimum tail index
 and start an attachment of new nodes to an existing graph 
 to be sure that each newly appearing node has at least one link to the 
 "dominating"
  communities. 
We allow a newly attached node to link to nodes of any community of the seed network 
and to assign the new node to class $i$ if it has at least one directed link to nodes of the $i$th community. 
The classes may correspond to in- or out-degree. The $i$th community corresponds to the $i$th column of matrix $A^{(0)}_*$ in Example \ref{Exam1}. For instance, if $i=1$ holds, then Class 1 includes a set of nodes with PRs  calculated by the upper block  of $A^{(0)}_*$. 
\par
 Beforehand, we assign the code $00...0$ of the length $N_C$ to each new node.
  $N_C$ is a  number of classes equal to the number of communities. The communities are numbered in ascending order of their tail indices. Once a new node has an edge  to the $i$th community,  its code is replaced by $00...iJ...J$,  $J\in \{0,1,2,...,N_C\}$, where $i$ stands at the $i$th position.  After $n$ nodes are appended, we obtain $n$ codes and  classify the nodes. For example, for $N_C=3$ we assign the nodes with codes $123$, $103$, $120$, $100$ to Class $1$,  with codes $023$, $020$ to Class $2$ and  with code $003$ to Class $3$. One can distinguish classes regarding  in-coming or out-going edges to the $i$th community.
  Since an attachment tool establishes an appearance  time  of new nodes, one can estimate the extremal index by an estimator that is determined for sequences, e.g., by the intervals estimator.
   \begin{algorithm}
  \caption{Classification of newly appended nodes}\label{Algor2}
  \begin{algorithmic}[1]
\State\label{S1} Select an 
initial 
directed graph with $n>1$ nodes 
as a seed network 
and calculate the PRs of its nodes.
\State Partition the seed network into $N_C$ communities, where  $N_C$ is predefined beforehand.
\State Check the stationarity of the communities and in case of non-stationarity select another seed network. 
\State\label{S4} Estimate the tail index of each community, e.g. by the Hill estimator (\ref{7}), and rank communities in ascending order of their tail indices. 
\State\label{Al7} Attach $N_0$ new nodes and corresponding new edges by the PA schemes (see Appendix \ref{PA}) to the 
nodes 
of the communities. 
\State Encode each newly appended node 
according to its edges to the $i$th community. 
\State Assign each new node to one of $N_C$ classes according to its code of the length $N_C$: codes $1J...J$ imply Class $1$ and codes $00..0iJ...J$, $J\in \{0,1,2,...,N_C\}$ with $i-1$ zeroes to Class $i$, $2\le i\le N_C$.
\State Estimate the extremal index of each community, 
e.g., by a modified intervals estimator presented in Algorithm \ref{Algor3}.
\State \label{S7} Assign the minimum tail index to PRs and the MLMs of the new nodes from Class 1 and calculate their extremal index: 
(a) if the community with minimum tail index $k_1$ is unique, then the extremal index 
is equal to  the extremal index $\theta_1$ of the 
"dominating"
community;
(b) if there is a random number of "dominating" communities, then the MLMs of the new nodes from Class 1 have the extremal index of the  "dominating" community with the maximum PR, and if, in addition, arbitrary enumerated sequences of node PRs of the "dominating" communities satisfy conditions (A1) or (A2) (see, Appendix \ref{Sec5.2}), then the PRs  of the new nodes from Class 1 have the same extremal index as MLMs; the  conditions (A1) or (A2) for the "dominating" communities provide the same minimum tail index $k_1$ for the PRs and the MLMs of the new nodes from Class 1.
\State 
Class 2 obtains the second minimum tail index corresponding to the next set of communities in the range and the respective extremal index as in item \ref{S7} and in the same way classes with numbers $i>2$ obtain their tail and extremal indices.

\end{algorithmic}
\end{algorithm}
The classification that is described in Algorithm \ref{Algor2} is supported by Theorem \ref{Prop1}. It is explored in Sections \ref{Sim}, \ref{RealData}. Item  \ref{S7},(a) is justified by Theorem 4 in Markovich \& Rodionov (2020), and  Item  \ref{S7},(b) by Theorem \ref{T3}.

%% file: modification1.tex
\subsubsection{Modification of the intervals estimator for graphs}\label{ModInterEst}
\par 
If the  extremal index of the community with enumerated nodes exists, then it can be estimated by the intervals estimator (Ferro \& Segers, 2003)
 \begin{eqnarray}\label{15}
\hat{\theta}_n(u)&=&\Big\{
\begin{array}{ll}
\min(1,\hat{\theta}_n^1(u)), \mbox{ if } \max\{(T(u))_i : 1 \leq i \leq L -1\} \leq 2,\\
\min(1,\hat{\theta}_n^2(u)), \mbox{ if } \max\{(T(u))_i : 1 \leq i \leq L -1\} > 2,
\end{array}
\end{eqnarray}
where
\begin{eqnarray}\label{16}&&\hat{\theta}_n^1(u)=\frac{2(\sum_{i=1}^{L-1}(T(u))_i)^2}{(L-1)\sum_{i=1}^{L-1}(T(u))_i^2},\end{eqnarray}
\begin{eqnarray}\label{17}&&\hat{\theta}_n^2(u)=\frac{2(\sum_{i=1}^{L-1}((T(u))_i-1))^2}{(L-1)\sum_{i=1}^{L-1}((T(u))_i-1)((T(u))_i-2)},\end{eqnarray}
$L-1$ is the random number of the inter-exceedance times $\{(T(u))_i\}$.
\begin{figure}[tbp]
 \begin{minipage}[t]{\textwidth}
\centering
 \subfigure[]{\includegraphics[width=0.45\textwidth]{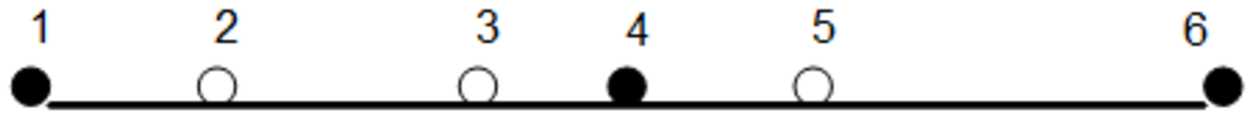} 
\label{fig:3a} 
}
%
%
\subfigure[]{\includegraphics[width=0.45\textwidth]{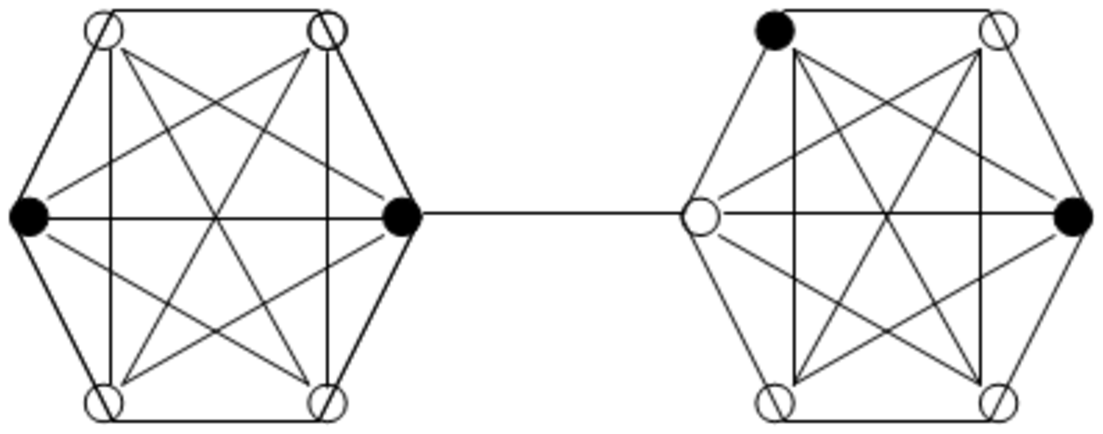} 
\label{fig:3b}
}
  \caption{The simplest chain graph (\ref{fig:3a}) and the $n$-barbell graph (\ref{fig:3b}) with significant and non-significant nodes
  marked by filled and open circles, respectively.}
\label{fig:3}
\end{minipage}
\end{figure}
  $T(u)$ denotes 
  a r.v. equal in
distribution to r.v.
\begin{eqnarray*}&&\min\{j\ge 1: X_{j+1}>u\}~~\mbox{given}~~ X_{1}>u,\end{eqnarray*}
or it holds
\begin{eqnarray*}P\{T(u)=n\}&=&P\{M_{1,n}\le u, X_{n+1}>u|X_{1}>u\}~~\mbox{for}~~ n\ge 1,\end{eqnarray*}
where $M_{i,j}=\max\{X_{i+1},...,X_{j}\}$, $M_{1,1}=-\infty$,
for the underlying sequence $\{X_n\}$ with the cdf $F(x)$.
For exceedance
times $1\le S_1<...<S_L\le n$ it follows
\begin{equation}\label{ti1}T(u)_i=S_{i+1}-S_{i}, \qquad i\in\{1,...,L-1\}.\end{equation}
$T(u)$ implies the number of observations running under
$u$ between two consecutive exceedances.
 $T(u_n)$ normalized by the tail function $\{Y=\overline{F}(u_n)T(u_n)\}$
 is derived to be asymptotically exponentially distributed with a weight $\theta$ and with an atom at zero with a weight $1-\theta$ (Ferro \& Segers, 2003). Here, $\theta$ is the extremal index.
The advantage of the intervals estimator is that it 
requires only the threshold value  $u$ as a parameter.  $u$ can be found as a high quantile of  $\{X_n\}$ since then the inter-exceedance times are approximately independent. The event $\{T(u)=1\}$ 
corresponds to neighbor exceedances. Such $\{T(u)\}$ have to be excluded from consideration (S\"{u}veges \& Davison, 2010).
\par
The intervals estimator was proposed for random sequences. Let us give an intuition, why the application of the intervals estimator to random graphs is plausible.
Clusters of exceedances form asymptotically a Poisson process with rate $\tau\theta$, where $\tau$ is taken from (\ref{11}) and $\theta$ is the extremal index of the underlying process $X_i$ (Beirlant et al., 2004; Rootz\'{e}n, 1988).
\par
A naive approach is to determine locations of nodes with PRs consistently exceeding a sufficiently high $u$ by a Poisson  process $P_n$ of rate $\tau\theta$. Roughly speaking, we superimpose the $P_n$ on the graph. The weight of an edge of the graph defined on a Poisson number of points may be taken equal to the inter-exceedance times. The latter may be measured as the number of nodes on the path 
(not equal to one) between a pair of nodes with exceedances of PRs. A rigorous proof is out of scope of the paper.
\par
Nodes in random graphs can be arbitrarily enumerated. Thus, the intervals estimator has to be  modified to avoid an enumeration of nodes (Algorithm \ref{Algor3}).
Regarding a graph community, $T(u)$ can be taken equal to the length of the 
path expressed
in edges between two nodes whose influence
indices exceed the threshold $u$.  All internal nodes along the path should have the influence indices less than $u$.
\\
Regarding a simple example, the
chain graph $G=(V,E)$ with edges $E=\{\{1,2\}, \{2,3\},..., \{m-1,m\}\}$, $m=6$ is shown in Fig. \ref{fig:3a}. 
We have $T(u)_1=3$ and
$T(u)_2=2$. For the $2-$barbell graph in Fig. \ref{fig:3b} we get 
$\{T(u)\}=\{2,2,2,2,3,3,3,3,3,3,4,4,4,4\}$ 
excluding  single edges. 
\begin{algorithm}
\caption{Modified intervals estimator}\label{Algor3}
  \begin{algorithmic}[1]
\State Let $\{X_i\}$ be 
 influence node indices, 
 e.g., PRs.
\State Take a high quantile of $\{X_i\}$ as threshold $u^*$.
\State Find nodes with exceedances, i.e. such that $X_i>u^*$ holds.
\State Set sequences $X_{xy}=\{X_{x}, X_{i_1}, X_{i_2}, \dots, X_{i_m}, X_{y}\}$, $m\ge 1$ corresponding to 
paths $x\rightarrow \dots\rightarrow y$ from a 
node $x$ to each node $y$. 
\State Define $\{T(u^*)_i\}$, $i\in\{1,2,...L\}$, by (\ref{ti1}) for all sequences $X_{xy}$ such that influence indices $X_x$ and $X_y$ 
exceed $u^*$ but $X_{i_1}, X_{i_2}, \dots, X_{i_m}$ do not, where $L$ is a total number of inter-exceedance times. 
\State Estimate the extremal index $\widehat{\theta}=\widehat{\theta}(u^*)$ of  the graph by (\ref{15})-(\ref{17}). 
 \end{algorithmic}
\end{algorithm}
The number of inter-exceedance times $\{(T(u))_i\}$ used in the modified intervals estimator may be larger than in the initial intervals estimator 
and its variance 
would be smaller.

%% file: simulation1.tex
\subsection{Inference and simulation}\label{Sim}
\unskip
We provide simulation results to support our theoretical conclusions. Our simulation concerns PRs of newly appended nodes only, but not their MLMs since the latter have the same tail and extremal indices by 
Markovich and Rodionov (2020a), 
Markovich (2021) and Theorem \ref{Prop1}. 
\par
We generate three Thorny Branching Tree (TBT) graphs. The TBT is a variation
of a branching tree where each node has an edge pointing to its parent. But it also has a certain number of unpaired outbound links  that are pointing outside of the tree (Chen et al., 2014a). 
Algorithm 1 proposed in Chen et al. (2014b) generates a bi-degree sequence in such a way that the in- and out-degrees of nodes follow closely the desired regularly varying distributions and  that
the sums of in- and out-degrees 
are the same. The bi-degree sequences are later used to
construct random graphs using the configuration model and TBT
simultaneously.
Our $TBT_1$, $TBT_2$, $TBT_3$ have power-law distributed  in- and out-degrees  with tail indices $(\iota_{in},\iota_{out})$ equal to $(3.8,2.0)$, $(2.5,2.5)$, $(3.0,4.5)$, respectively, Tab. \ref{Table1}. Due to the condition $\sum_{s=1}^n\iota_{in_s}=\sum_{s=1}^n\iota_{out_s}$ the tail indices of  simulated in- and out-degrees may slightly differ from the initial values  $(\iota_{in},\iota_{out})$.
\par
\begin{figure}[tbp]
 \begin{minipage}[t]{\textwidth}
\centering
 \subfigure[]{\includegraphics[width=0.3\textwidth]{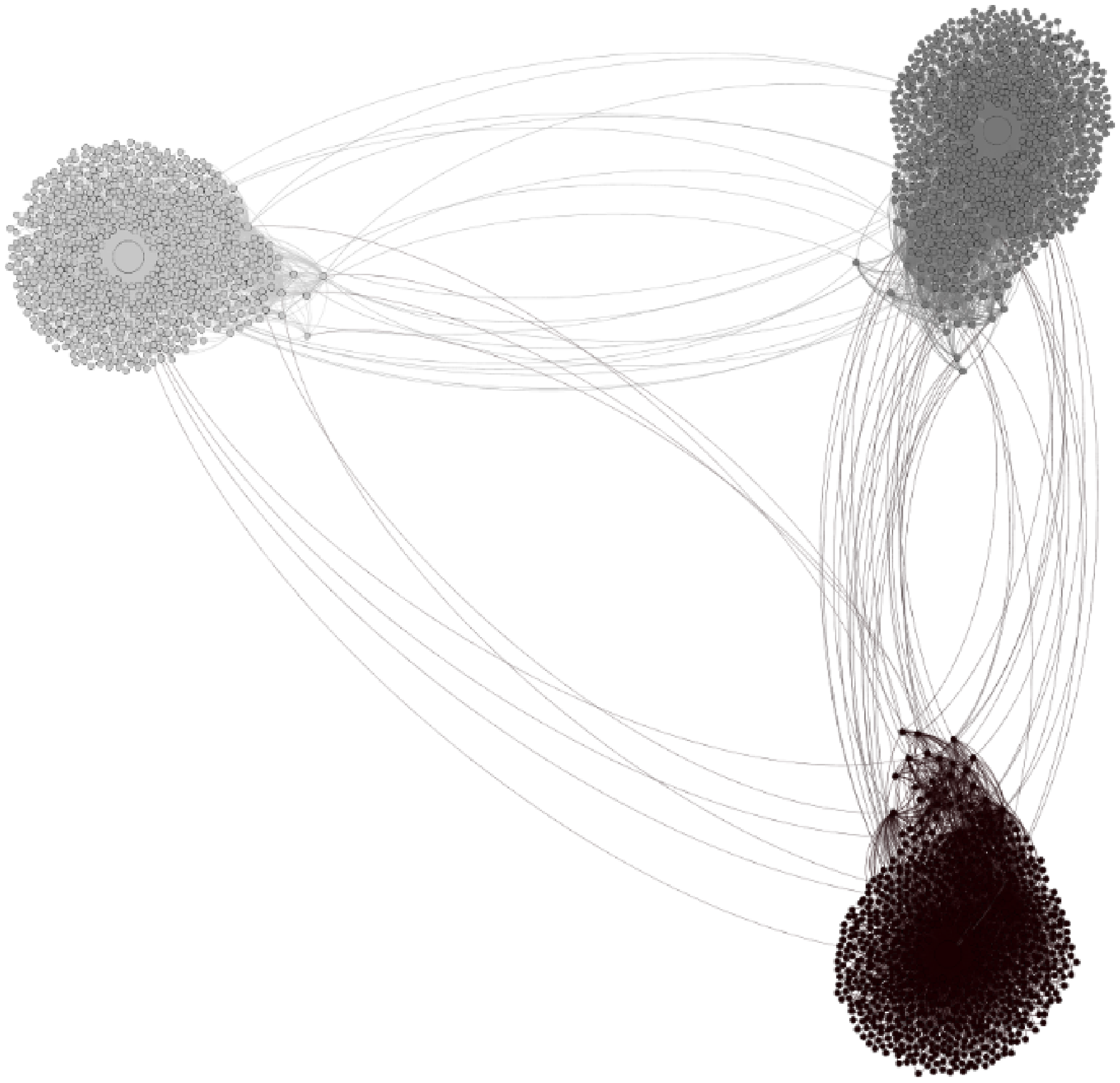} 
\label{fig:4a}}
\subfigure[]{\includegraphics[width=0.33\textwidth]{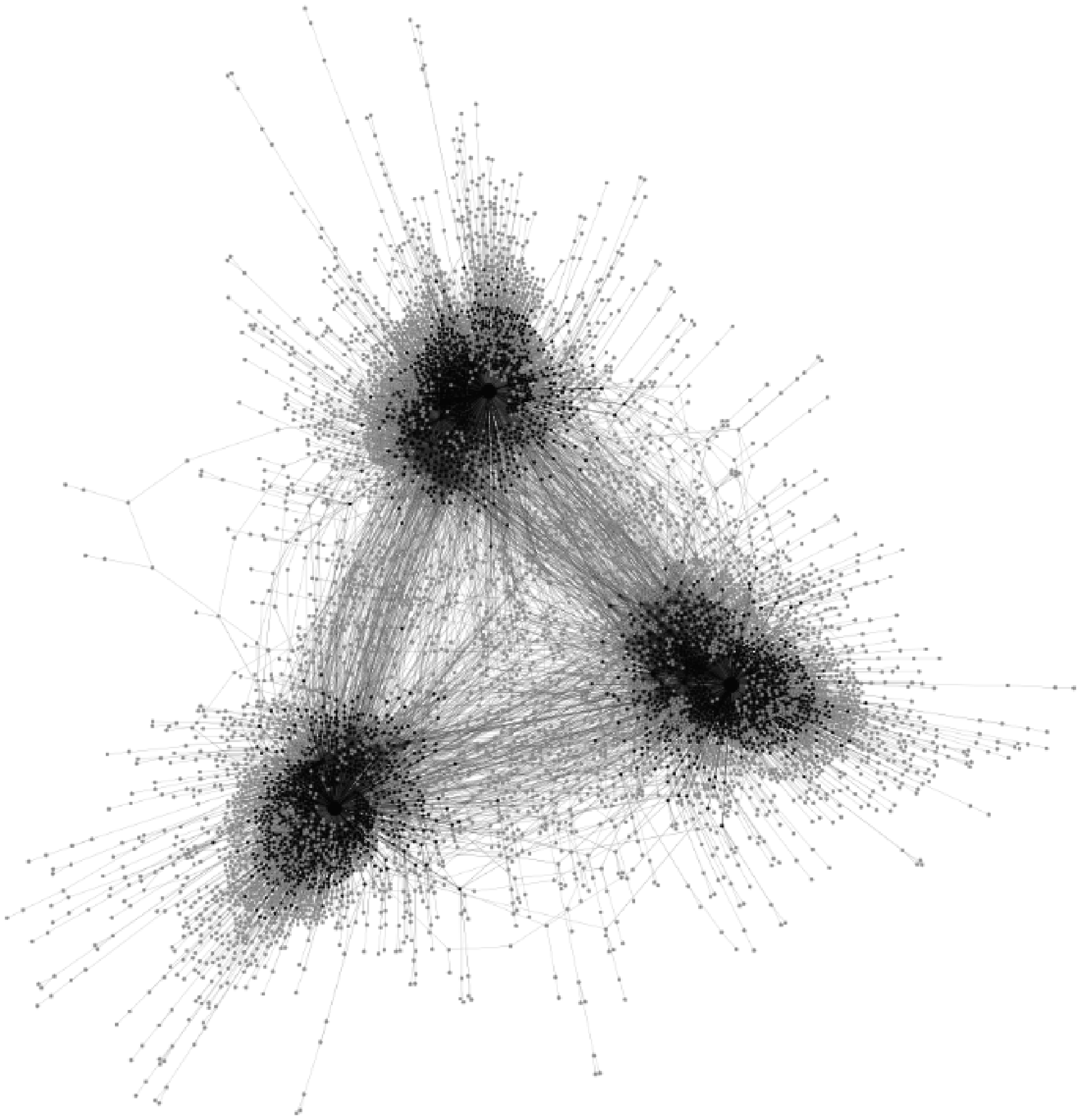}
\label{fig:4b}}
\subfigure[]{\includegraphics[width=0.33\textwidth]{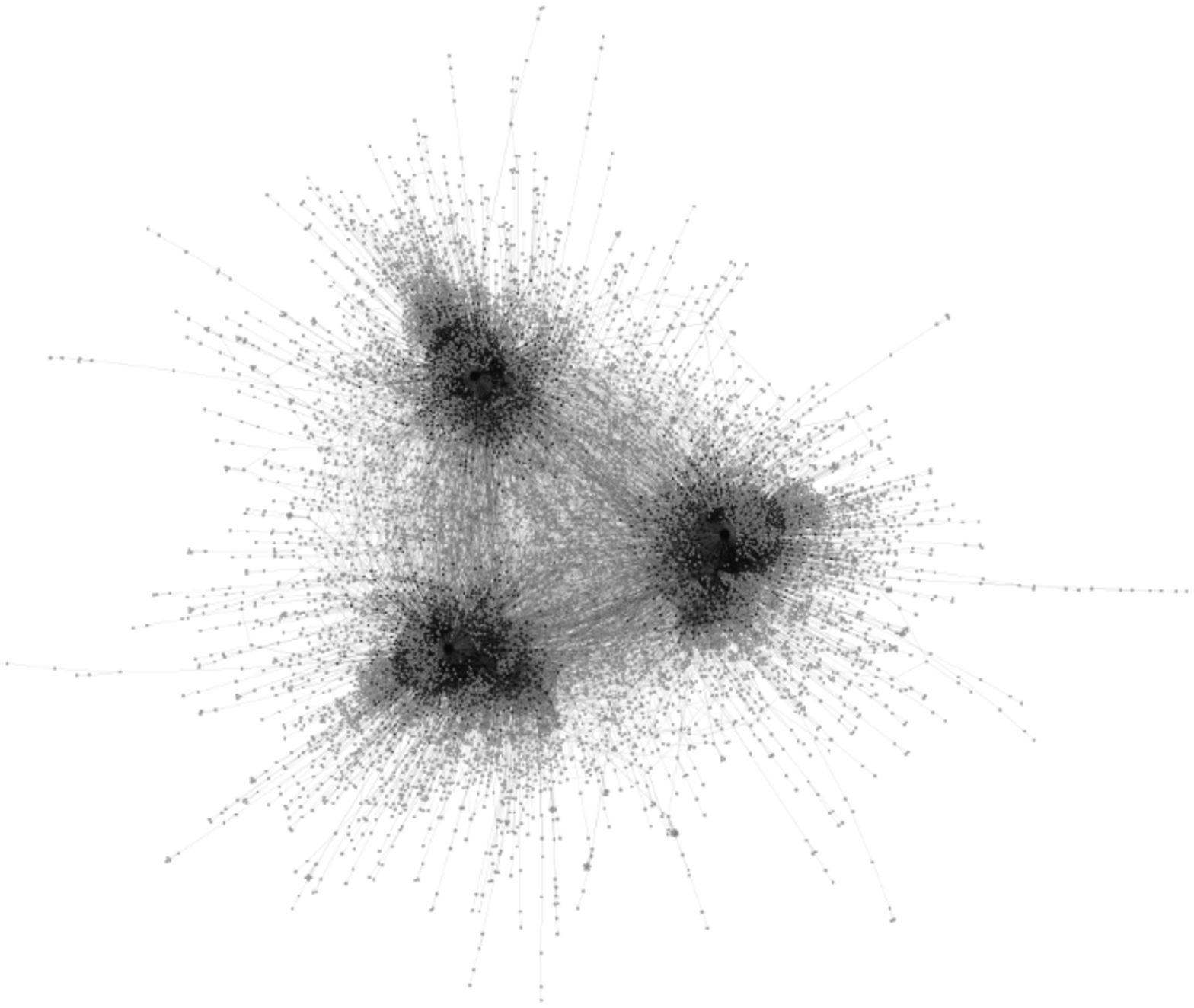}
\label{fig:4c}}
  \caption{The $TBT_1$ (black circles), $TBT_2$ (grey circles) and $TBT_3$ (white circles)
  before (Fig. \ref{fig:4a}) and after the  preferential attachment $PA(0.4;0.2;0.4)$ with parameters $\delta_{in}=\delta_{out}=1$ of $5493$ 
  (Fig. \ref{fig:4b}) and $10^4$ (Fig. \ref{fig:4c})
  new nodes, where 
  new nodes  are shown by grey 
  circles and "old" nodes by black ones; the circle sizes  are proportional to their PR values. 
  }
\label{fig:4}
\end{minipage}
\end{figure}

The TBTs are further connected by $100$ additional edges to simulate the connection 
between TBTs, see Fig. \ref{fig:4a}. 
Each TBT contains  
$800$ nodes.  
 To append new nodes and edges, $\alpha-$, $\beta-$ and $\gamma-$ schemes of the linear $PA(\alpha,\beta,\gamma)$ (see (\ref{al_sch})-(\ref{ga_sch}) in Appendix \ref{Attach}) 
 are used.
  The number of attached nodes $N_0$ is taken equal to $5493$ 
  and $10^4$, Fig. \ref{fig:4b} and \ref{fig:4c}.
  \\
We evaluate the tail index of the PRs by Hill's estimator (\ref{7}) where $k$ is selected by a bootstrap method.
The consistency of Hill's estimator for the PR
tail index has not been justified rigorously. The proof of the consistency is out of scope of this paper. The numerical consistency indirectly follows from Tab. \ref{Table1} since Hill's estimates decrease as the number of appending nodes $N_0$ increases. This reflects the appearance of giant nodes with a large node degree as the number of newly created edges to the TBTs increases.
\\
We consider the TBTs ("old" nodes) before and after  the PA of new nodes and edges,
and classes of newly appended nodes to the existing ones by the PA schemes ("new" nodes). 
Bootstrap confidence intervals for the Hill's estimates are calculated by $500$ bootstrap resamples (Markovich, 2007).
\par The classes of "new" nodes are built by an encoding of nodes as in Algorithm \ref{Algor2}. The "out-degree-classes" correspond to out-going edges from "new" nodes to "old" ones
and vice versa for the "in-degree-classes". 
For our $TBT_1$-$TBT_3$, a new node $v$ may be encoded by one of the codes $123$, $103$, $120$, $100$  and related to "$in-degree-class_1$" (or to "$out-degree-class_1$"), if the edge $(v,w)$ leads from an  "old" node $w$ belonging to $TBT_1$ to $v$ (or vice versa).
$Class_4$ contains new nodes that  have no links to "old" nodes from the TBTs but only to previously appended new nodes.
The number of newly appended nodes $N_0$ (and thus, the cumulative size of $Class_1$-$Class_4$) is random. 
It may be less or equal to the number of evolution steps due to the $\beta-$scheme. The latter creates a new edge between two existing nodes and no new node is added. The PA-scheme is selected by means of a trinomial r.v. (see Appendix \ref{PA}).
$Class_4$ can be re-encoded and  divided  into classes with regard to their edges to $Class_1$-$Class_3$. The latter are considered further as seed communities.
\par
Analyzing the tail indices of PRs of  "old" nodes in Tab. \ref{Table1}, one can conclude that the $TBT_1$  
is the 
"dominating"
community since it has the minimum tail index.\footnote{Since PRs are regularly varying distributed, the smaller positive tail index implies the heavier  tail by Breiman's theorem (Jessen \& Mikosch, 2006; Markovich, 2007).} 
\begin{table}[t]
\caption{The Hill's estimates 
$\widehat{\alpha}_{1}(n,k)$ and $\widehat{\alpha}_{2}(n,k)$ of the PR tail index 
 of "old" nodes in the $TBT_1-TBT_3$ 
before and after the attachment of $N_0$ new nodes by the  $PA(0.4;0.2;0.4)$ with  $\delta_{in}=\delta_{out}=1$, $\widehat{\alpha}_{3}(n,k)$ and $\widehat{\alpha}_{4}(n,k)$
correspond to the "in-degree-classes" and  "out-degree-classes"  of sizes $n$
 with  $97.5\%$ bootstrap
 confidence intervals $(u_1,u_2)$.
 }
 \centering
\tabcolsep=0.001cm
\begin{tabular}{lccc||ccccc}
  \hline
     & $(\iota_{in},\iota_{out})$ & $\widehat{\alpha}_{1}(n,k)$ & $\widehat{\alpha}_{2}(n,k)$ &   & $n$ & $\widehat{\alpha}_{3}(n,k)$ &$n$ & $\widehat{\alpha}_{4}(n,k)$ \\
     &                            & $(u_1,u_2)$             & $(u_1,u_2)$                     &   &     & $(u_1,u_2)$                 & & $(u_1,u_2)$\\
     \cline{2-9}
     &                            &                             &                             \multicolumn{6}{|c|}{\scriptsize {$N_0=5493$}}\\ 
     \cline{4-9}
  $TBT_1$    & (3.8,2.0)          & 3.3142                  &      2.2937        &   $Class_1$ & 967 & 3.1469              & 788 & 2.5651 \\
             &                    &(2.553,10.534)         & (1.824,3.093)    &             &     & (2.286, 3.501)    &     &(2.015, 3.772)
  \\
  $TBT_2$    & (2.5,2.5)          & 3.3513                  &    2.4600          &  $Class_2$  & 824 & 3.9242              & 728 &  2.7956\\
             &                    &(2.532,10.889)         & (1.831,3.182)    &             &     & (2.861, 4.234)    &     &(2.264, 4.003)
   \\
   $TBT_3$   & (3.0,4.5)          & 3.4697                  &    2.5518          &  $Class_3$  & 608 & 3.8494              & 557 &  2.7853 \\
             &                    &(1.839,18.238)         &(1.998,3.518)     &             &     & (3.237, 4.621)    &     &(2.128, 4.003)\\
             &                    &                         &                    & $Class_4$   & 3094& 3.2298              & 3420& 3.7226  \\
             &                    &                         &                    &             &     & (2.684, 3.854)    &     &(2.804, 3.868)
  \\
  \hline
  &                            &                             &                             \multicolumn{6}{|c|}{\scriptsize {$N_0=10000$}}\\
     \cline{4-9}
     $TBT_1$    &              &                   &      2.1315        &   $Class_1$ & 1575 & 3.1863              & 1267 & 3.0418 \\
             &                    &         & (1.445,2.645)    &             &     & (2.073, 3.274)    &     &(2.425, 3.998)
  \\
  $TBT_2$    &           &                  &    2.1596          &  $Class_2$  & 1301 & 3.2211              & 1101 &  3.0914\\
             &                    &         & (1.368,2.551)    &             &     & (2.375, 3.457)    &     &(2.283, 3.885)
   \\
   $TBT_3$   &           &                 &    2.5445          &  $Class_3$  & 973 & 3.7493              & 860 &  3.3251 \\
             &                    &         &(1.782,3.057)     &             &     & (3.083, 4.298)    &     &(2.192, 3.799)\\
             &                    &                         &                    & $Class_4$   & 6151 & 3.3036              & 6772 & 3.5173  \\
             &                    &                         &                    &             &     & (2.872, 3.731)    &     &(2.522, 3.367)
   \\
  \hline
\end{tabular}
 \label{Table1}
\end{table}
\begin{table}[t]
\caption{Intervals and $K-$gaps estimates  of  PR extremal index of
 "old" nodes in the $TBT_1-TBT_3$ before and after the $PA(0.4,0.2,0.4)$ with  $\delta_{in}=\delta_{out}=1$ of $5493$ 
 new nodes and of the
new nodes in the "in-degree-classes" and  "out-degree-classes".
  }
 \centering
\tabcolsep=0.01cm
\begin{tabular}{lcccc||
ccccc}
  \hline
     & $\widehat{\theta}^{IA1}$  & $\widehat{\theta}^{Idis}$ & $\widehat{\theta}^{K0dis}$& $\widehat{\theta}^{KIMT}$ & &$\widehat{\theta}^{IA1}$ & $\widehat{\theta}^{Idis}$ & $\widehat{\theta}^{K0dis}$ & $\widehat{\theta}^{KIMT}$\\
     \hline
   \multicolumn{5}{c}{\scriptsize {Before PA}}&  \multicolumn{5}{c}{\scriptsize {"in-degree-classes"}}\\
  $TBT_1$    & 0.9167   & 0.8143 & 0.8282 & 0.7793 &  $Class_1$ & 0.9952 
  & 0.9479 & 0.9551 & 0.8218\\
             &          & (0.9167) & (0.9311) & &        &   & (1)    & (1)
  \\
  $TBT_2$    & 0.8632   & 0.8098 &0.8157   &0.8573  & $Class_2$ & 0.9977 & 0.9861 & 0.9879 & 0.8723\\
             &          & (0.9387)& (0.9445) &  &         &        & (1)    & (1) &
   \\
   $TBT_3$   & 0.9804    & 0.9783 & 0.9808 
   & 0.8683& $Class_3$ & 0.9173 & 0.9356 & 0.9466 & 0.8907\\
              &          & (0.8506) & (0.8579) 
              & &       &   & (0.9453)    & (0.9494)&
  \\
  \hline
  \multicolumn{5}{c}{\scriptsize {After PA}} & \multicolumn{5}{c}{\scriptsize {"out-degree-classes"}}\\
  $TBT_1$    & 1   & 0.9345 & 0.9377 & 0.9525& $Class_1$ & 0.2772 
  & 0.3130 & 0.3121 &0.6554\\
             &          & (1) & (1) &  &      &   & (0.1961)    & (0.1883) &
  \\
  $TBT_2$    & 0.3098   & 0.3175 &0.3389 & 0.5646 & $Class_2$ & 0.1715 
  & 0.2974 & 0.2818 & 0.7461\\
             &          & (0.3257) & (0.3237) &  &      &   & (0.1140)    & (0.1091) &
   \\
   $TBT_3$   & 0.9618    &0.9070 &0.9211 & 0.8529 & $Class_3$ & 0.2421 & 0.3970 & 0.4276 & 0.7106\\
              &          & (0.8732) & (0.8817) &  &      &   & (0.1987)    & (0.1870) &
  \\
  \hline
\end{tabular}
 \label{Table2}
\end{table}
\begin{figure}[tbp]
 \begin{minipage}[t]{\textwidth}
\centering
 \subfigure[]{\includegraphics[width=0.32\textwidth]{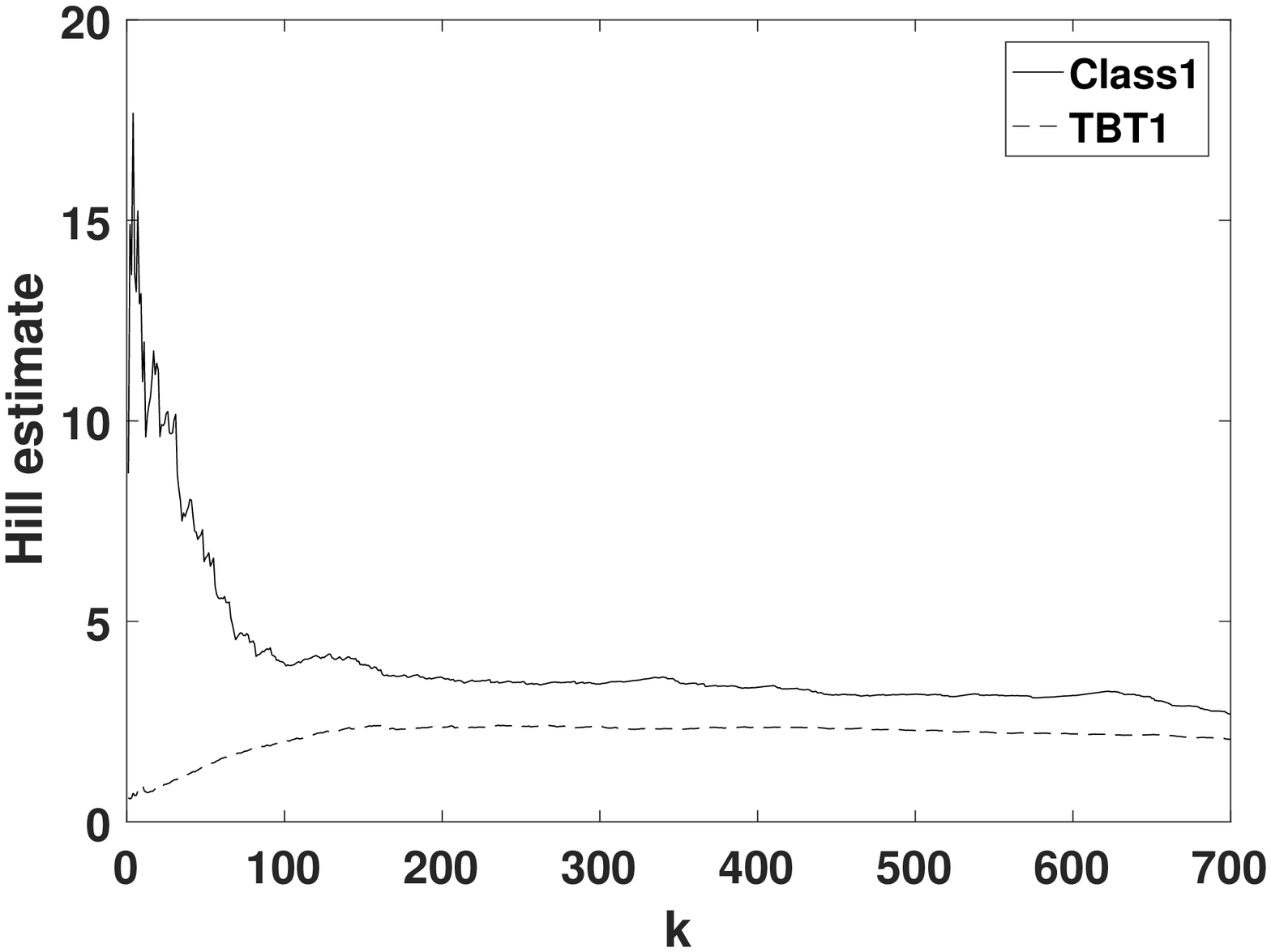}
\label{fig:5a}}
\subfigure[]{\includegraphics[width=0.32\textwidth]{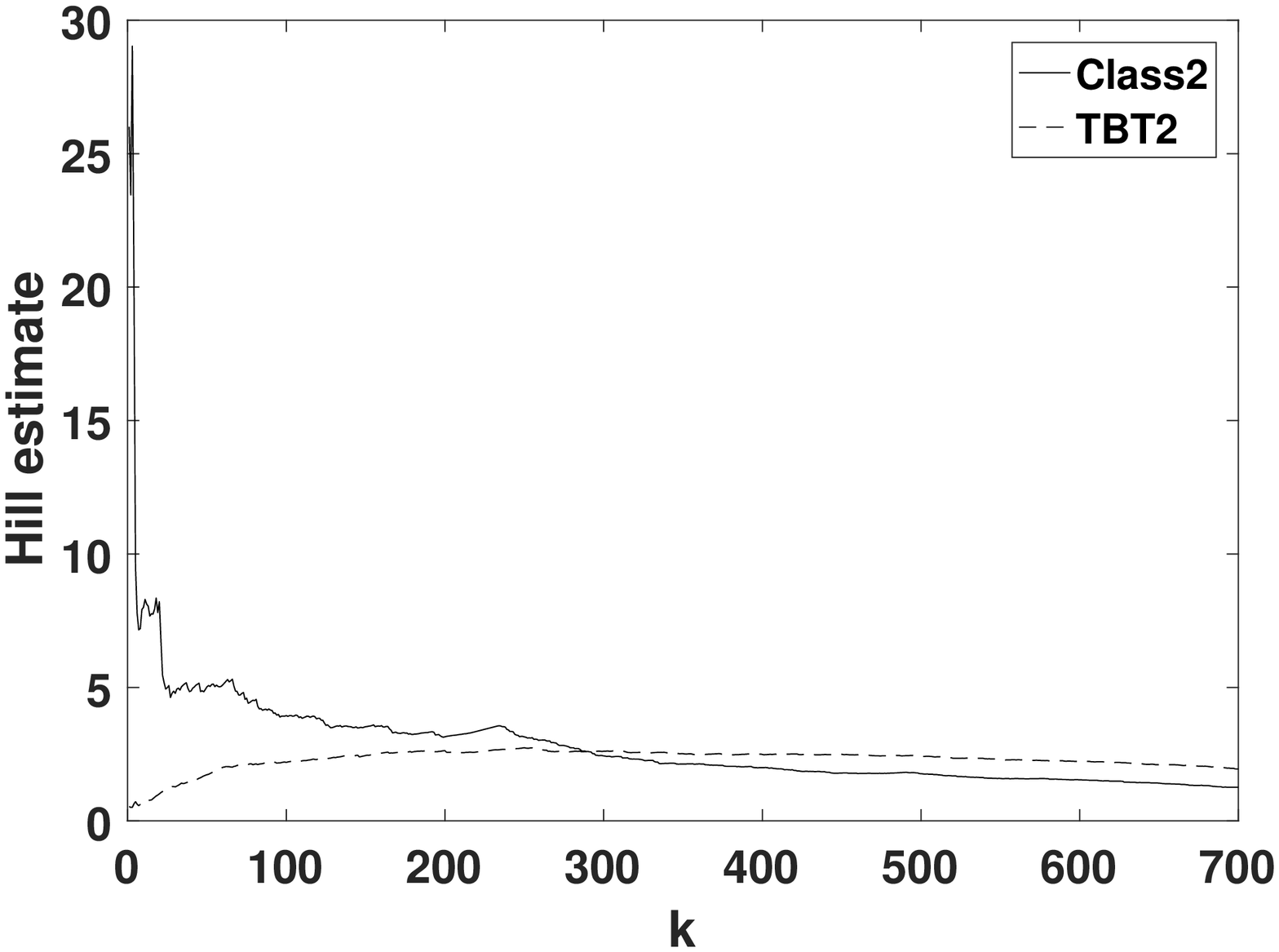}
\label{fig:5b}}
\subfigure[]{\includegraphics[width=0.32\textwidth]{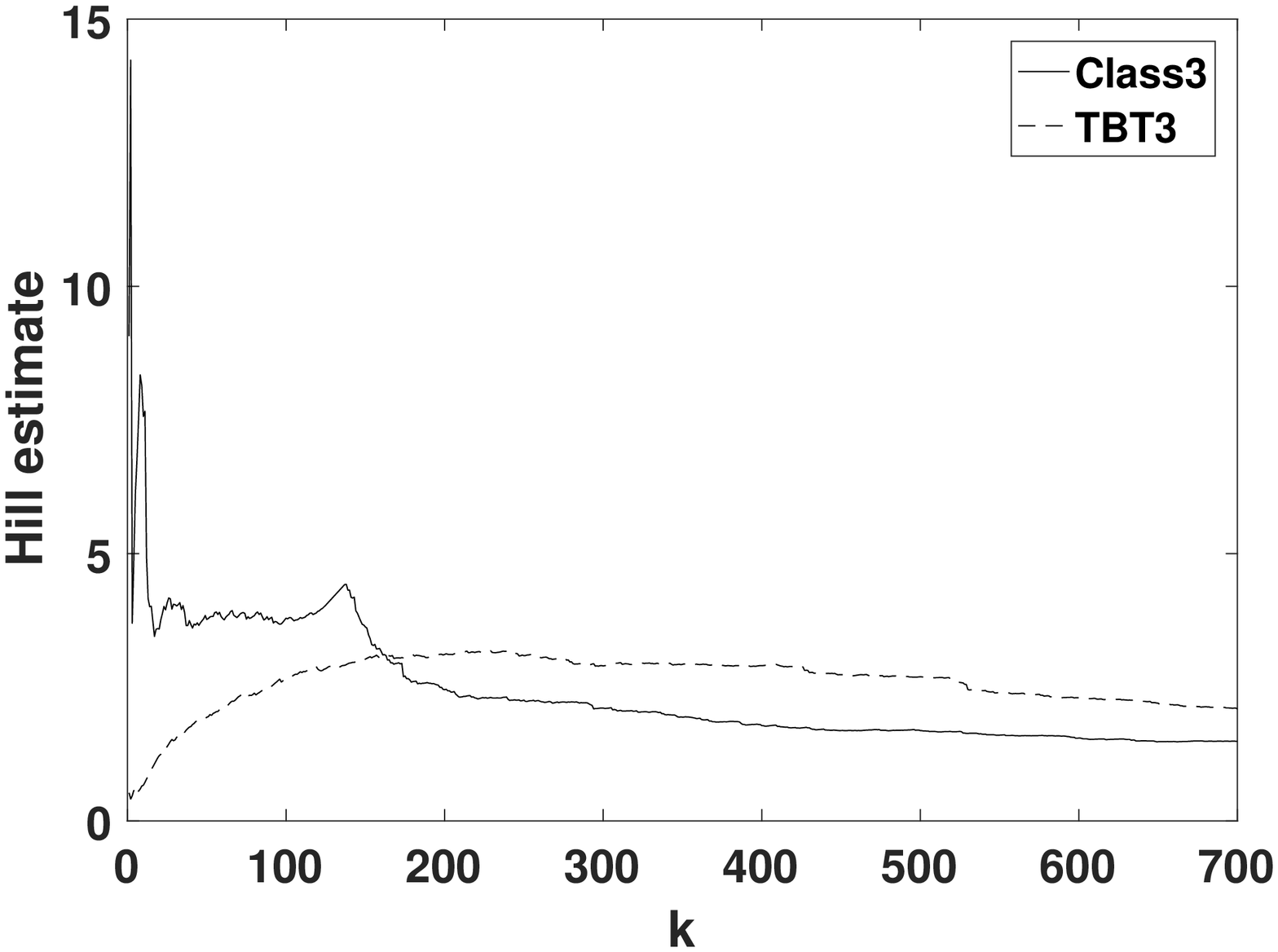}
\label{fig:5c}}
  \caption{The Hill's plots of node PRs of the $TBT_1-TBT_3$ after the appending of $5493$ 
  new nodes by $PA(0.4,0.2,0.4)$ with  $\delta_{in}=\delta_{out}=1$ and in-degree 
  $Class_1-Class_3$ versus the $k$ largest  order statistics  used in (\ref{7}).}
\label{fig:5}
\end{minipage}
\end{figure}
\begin{figure}[tbp]
 \begin{minipage}[t]{\textwidth}
\centering
\subfigure[]{\includegraphics[width=0.32\textwidth]{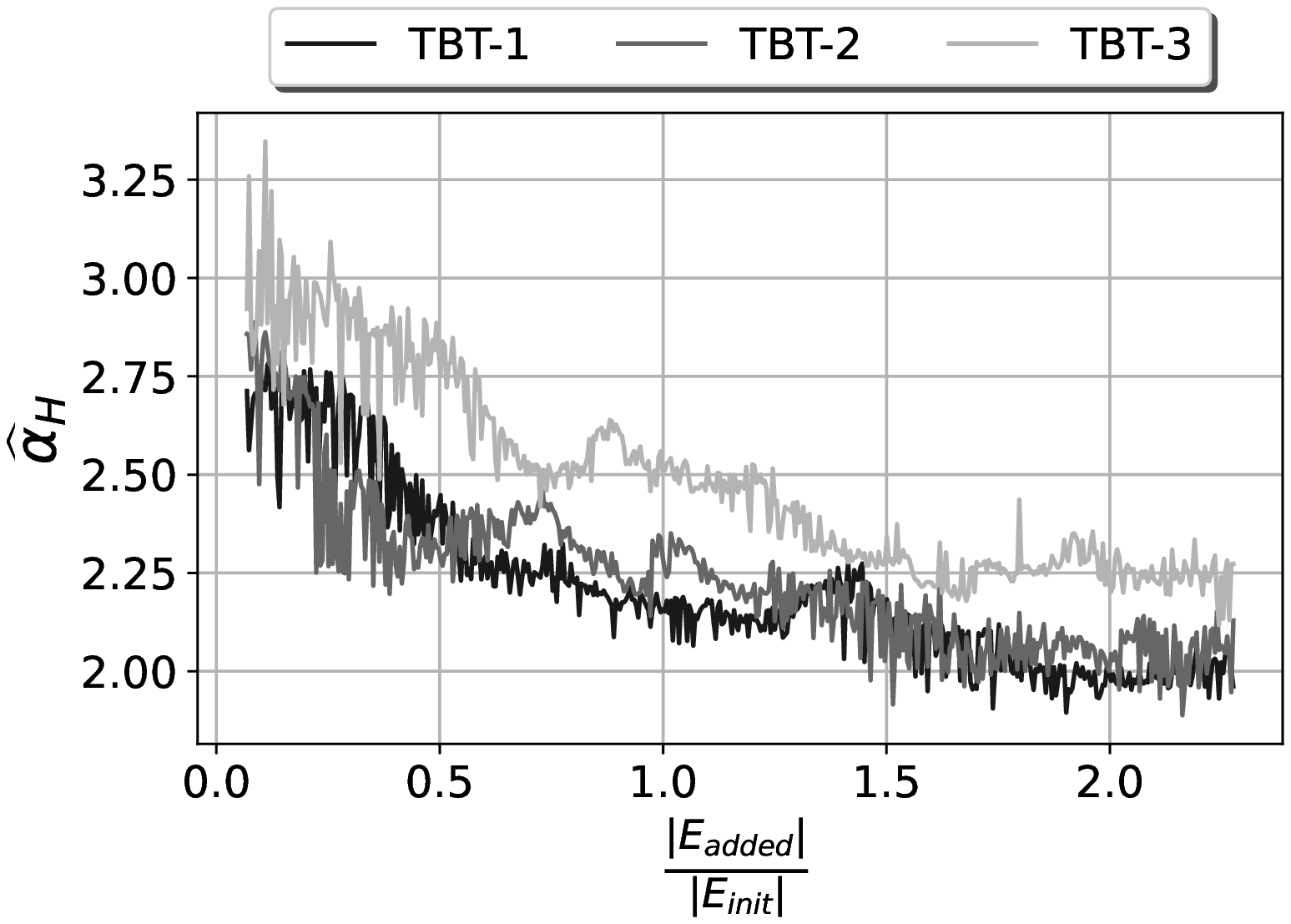} 
\label{fig:6a}}
\subfigure[]{\includegraphics[width=0.32\textwidth]{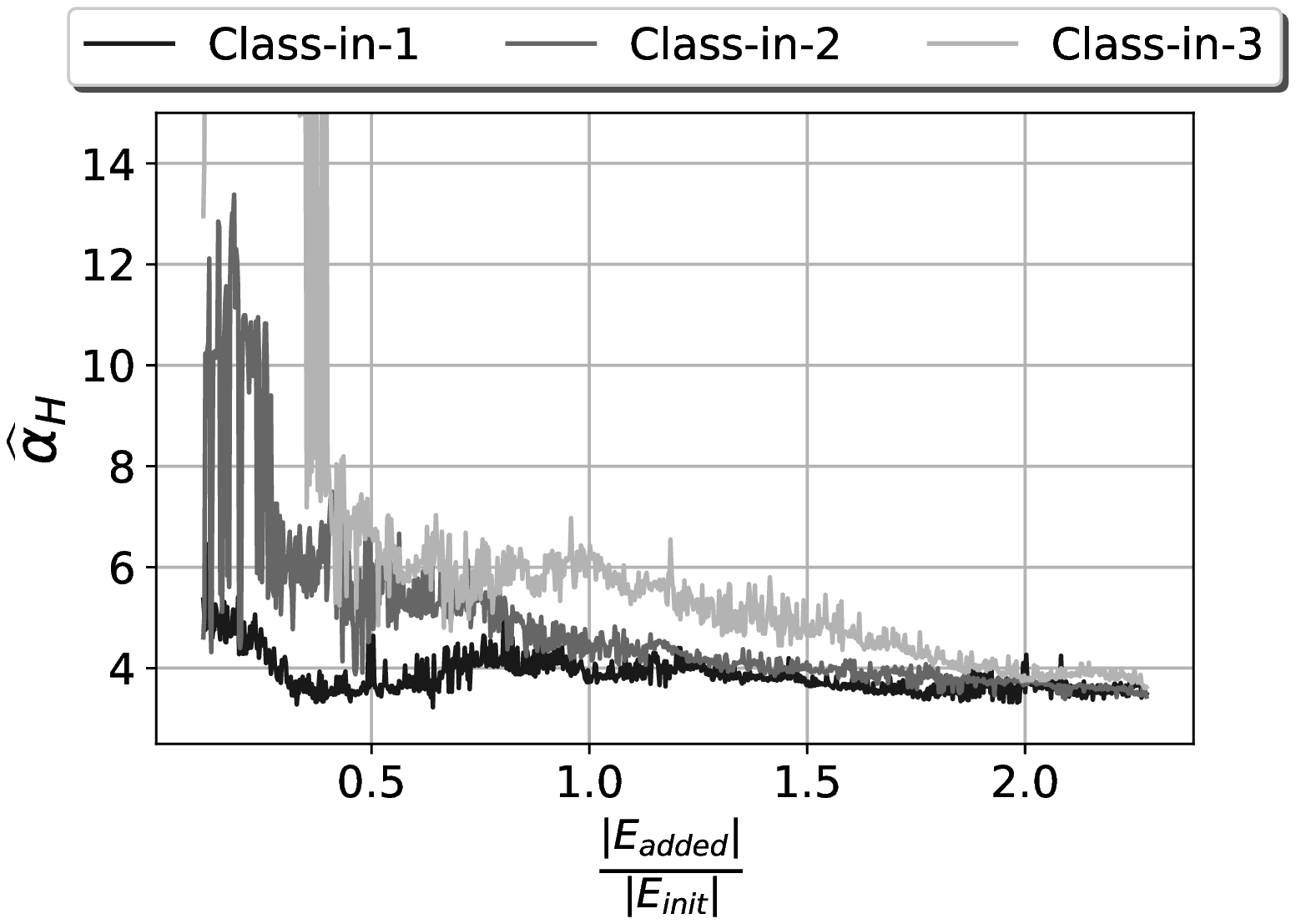}
\label{fig:6b}}
\subfigure[]{\includegraphics[width=0.32\textwidth]{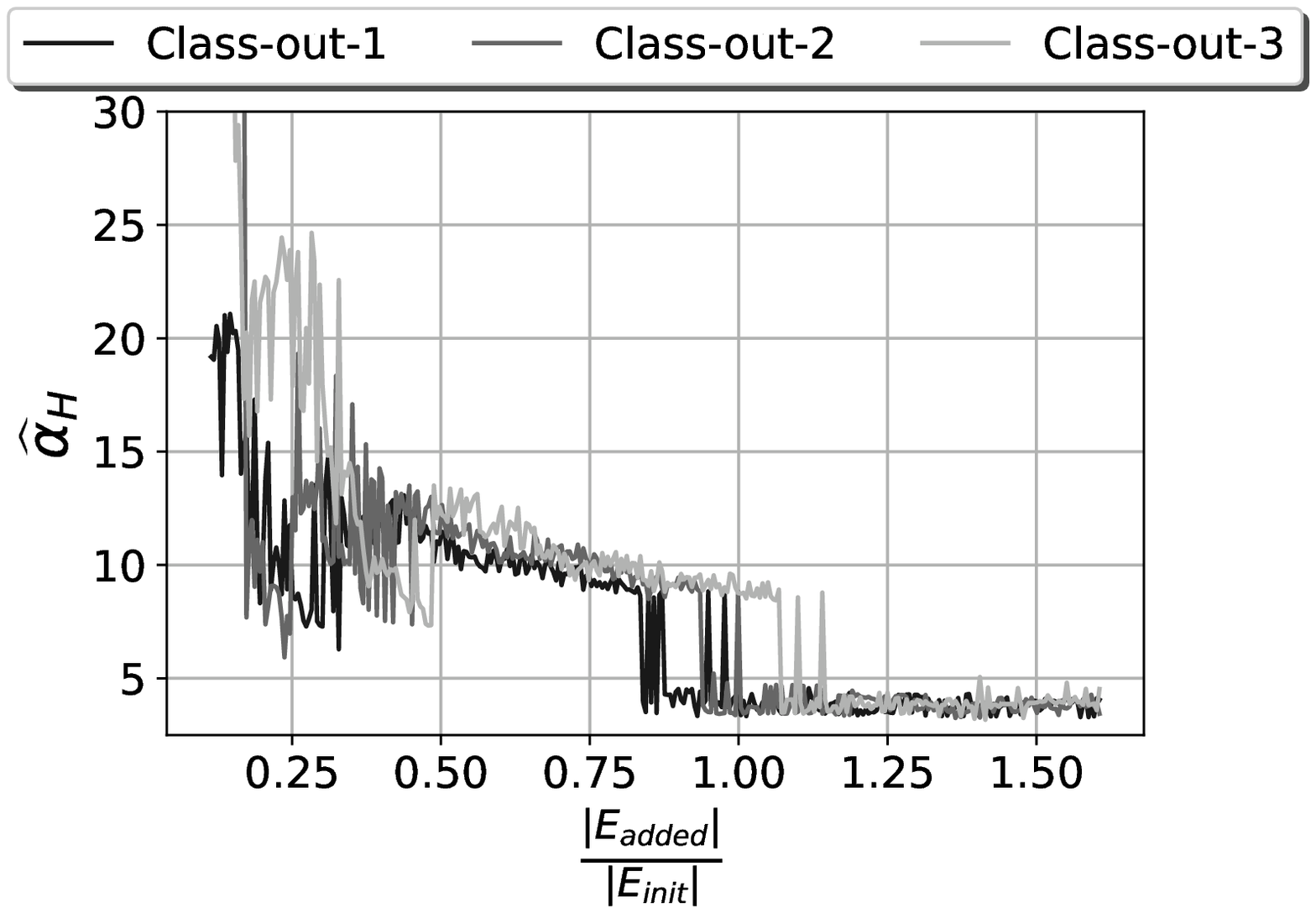}
\label{fig:6c}}
  \caption{The Hill's estimates $\widehat{\alpha}(n,k)$ of PR of the TBTs (Fig. \ref{fig:6a}) and   in- and out- degree classes (Fig. \ref{fig:6b}) and (Fig. \ref{fig:6c}) versus $|E|_{added}/|E|_{init}$, where $|E|_{added}$ is the number of edges added to the  TBTs during the PA
and the initial number of edges $|E|_{init}$ in  disjoint $TBT_1-TBT_3$.
  }
\label{fig:6}
\end{minipage}
\end{figure}
The Hill's estimates $\widehat{\alpha}_{1}(n,k)$ and $\widehat{\alpha}_{3}(n,k)$ as well as $\widehat{\alpha}_{2}(n,k)$ and $\widehat{\alpha}_{4}(n,k)$ of the corresponding  $TBT_i$ and $in(out)-degree-class_i$, $i\in\{1,2,3\}$ are close, see Tab. \ref{Table1} and Fig. \ref{fig:5}. The closeness grows up as the number of appended nodes increases.
The PR tail indices of the "old" nodes of TBTs decrease due to the appending of "new" nodes. Hence, their distribution tails become heavier.
 This property implies the appearance of giant components in the TBTs. It is in agreement with the  conclusion that the PA leads to giant nodes, since each new vertex tends to get connected to big others, 
 rather than to small ones (Krapivsky \& Redner, 2001;  Norros \& Reittu, 2006).
  The tail index drop is stabilized as the number of appended nodes increases, since the newly appearing nodes  gradually lose connection with the "old" ones. The increasing size of $Class_4$ reflects this effect also.
\par
The same tail index drop is shown in Fig. \ref{fig:6}, where the impact of the ratio $|E|_{added}/|E|_{init}$ is represented.
In total, $20000$ edges and $16083$ nodes were added.  
It is shown that the tail indices of both in- and out-degree classes tend to the tail indices of the TBTs before the evolution
as the number of attached edges grows up. 
This observation is in agreement with the item (i) of Theorem \ref{Prop1}  where the TBTs play the role of the columns of the matrix $A^{(0)}$. The values are also in agreement with formulae (2.9) in Wan et al. (2009) for tail indices of the in- and out-degrees: $\iota_{in}=(1+\delta_{in}(\alpha+\gamma))/(\alpha+\beta)$ and $\iota_{out}=(1+\delta_{out}(\alpha+\gamma))/(\beta+\gamma)$. Regarding our model $PA(0.4,0.2,0.4)$ with $\delta_{in}=\delta_{out}=1$ one obtains $\iota_{in}=\iota_{out}=3$. The tail indices of the PR and the in-degree are comparable (Litvak et al., 2007).
\par
\begin{figure}[tbp]
\begin{minipage}[t]{\textwidth}
\centering
\subfigure[]{\includegraphics[width=0.32\textwidth]{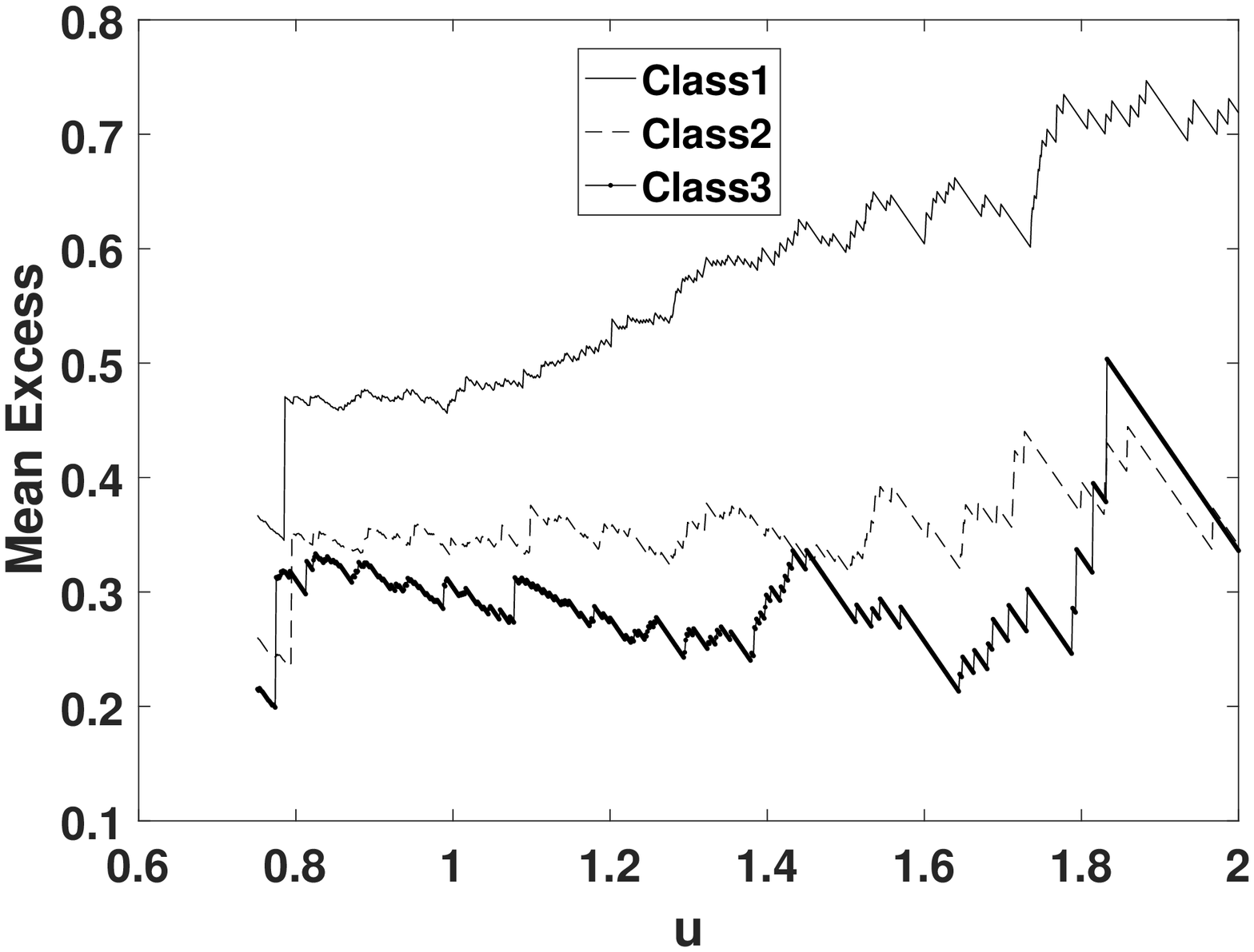}
\label{fig:7a}}
\subfigure[]{\includegraphics[width=0.32\textwidth]{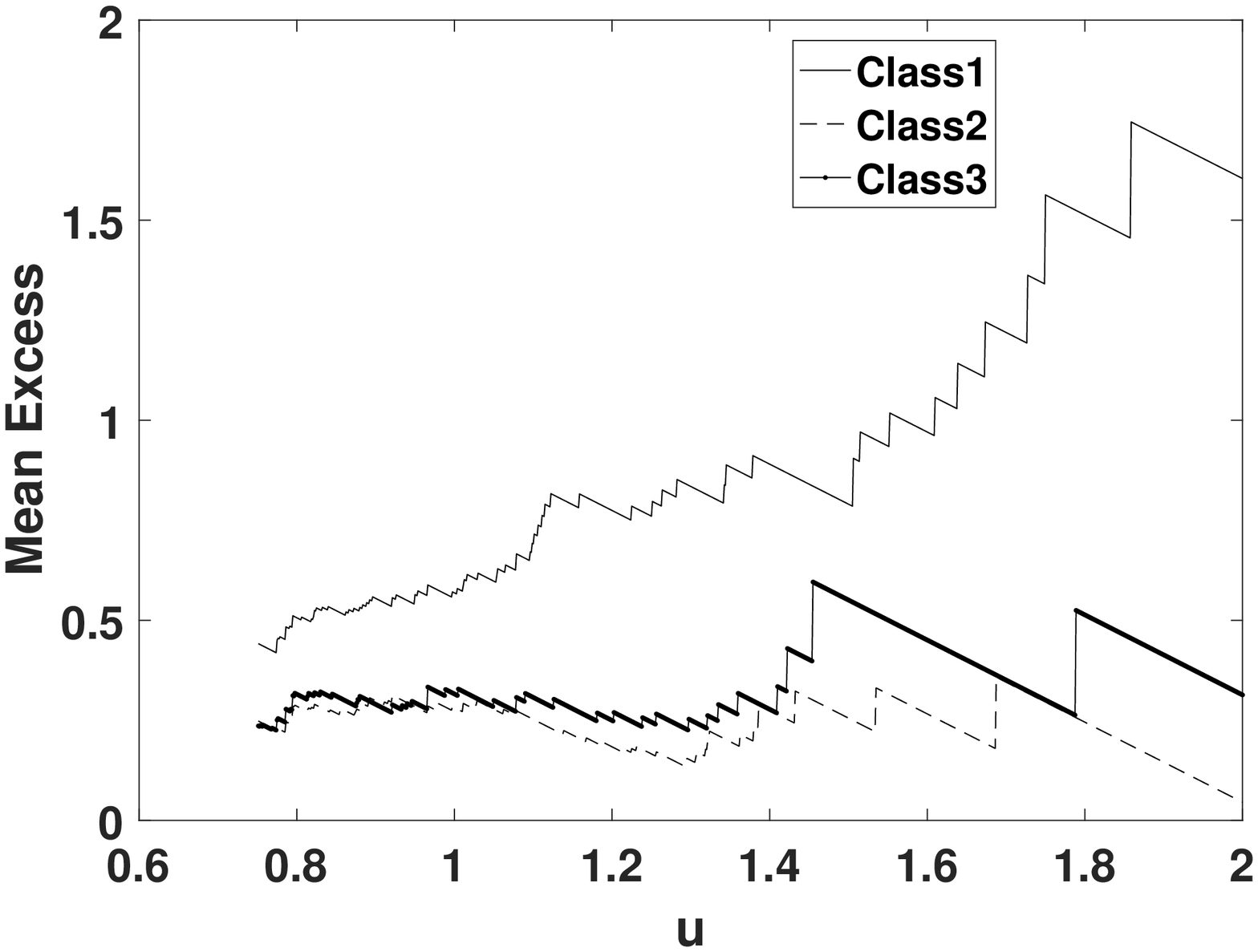}
\label{fig:7b}}
\subfigure[]{\includegraphics[width=0.32\textwidth]{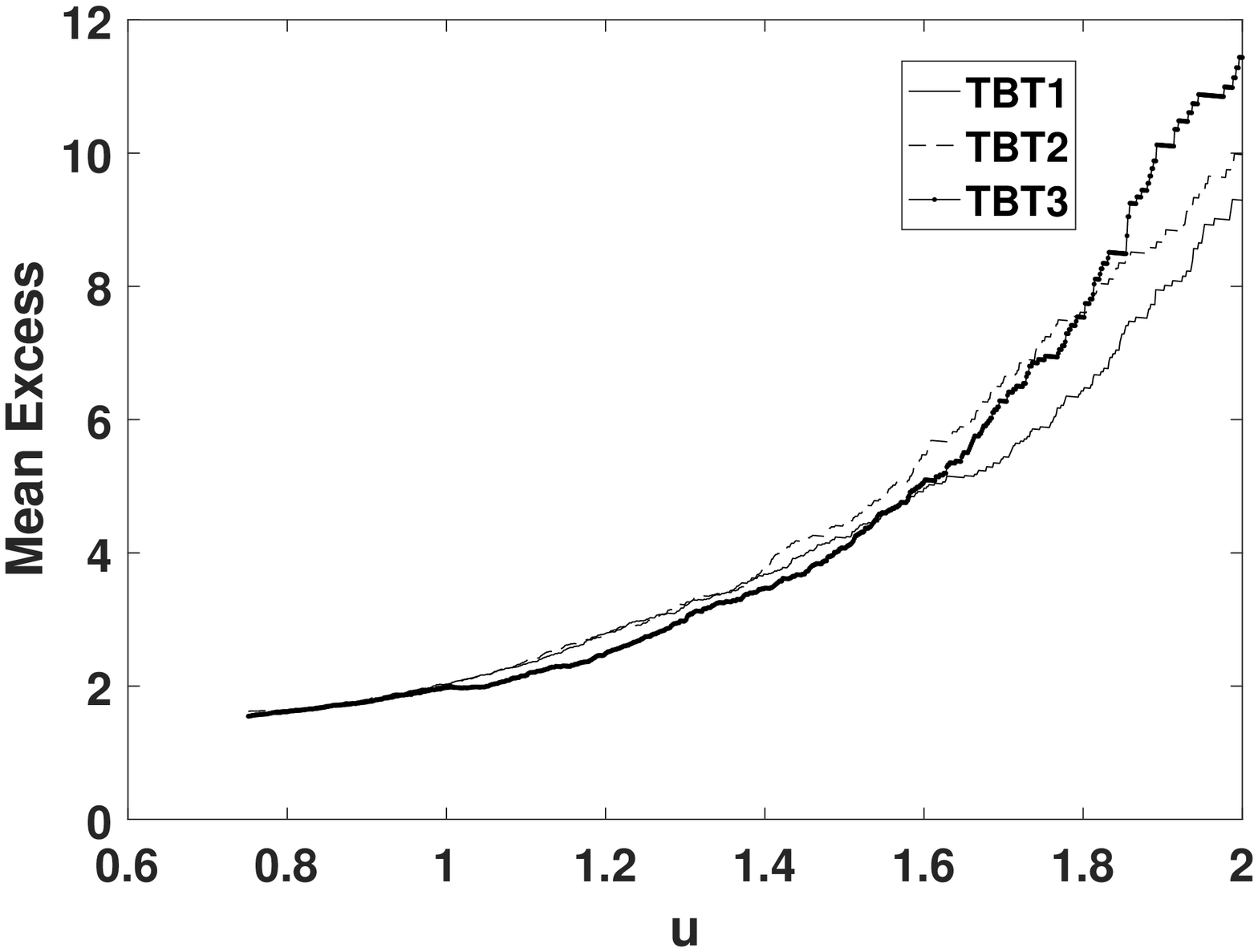}
\label{fig:7c}}
\caption{Mean excess function of node PRs in the "in-degree-classes" (Fig. \ref{fig:7a}), the "out-degree-classes" (Fig. \ref{fig:7b}) and the TBTs before 
the PA (Fig. \ref{fig:7c}).}
\label{fig:7}
\end{minipage}
\end{figure}

The TBTs before the PA are stationary distributed by the simulation. We check the stationarity of  PRs in the TBTs before 
the PA and in the classes of "new" nodes. 
To this end, we check the mean excess function $e(u)=E\{X-u|X>u\}$ and calculate its sample analogue $e_n(u)=\sum_{i=1}^n(X_i-u)\Ii(X_i>u)/\sum_{i=1}^n\Ii(X_i>u)$, see Fig.\ref{fig:7}.
The increasing, the decreasing and a constant value of $e(u)$ imply heavy-, light-tailed and exponential distributions, respectively (Markovich, 2007). The linear increase of $e(u)$ implies  Pareto-like distributions of the underlying random sequences. Since for Pareto distribution
\begin{eqnarray}\label{2}e(u)&=&(1+\gamma u)(1-\gamma), \qquad\gamma<1,\end{eqnarray} holds, where $\gamma=1/\alpha$ is the reciprocal of the tail index, one can see that the lower curves at the plot of the mean excess function correspond to the larger $\alpha$. Fig. \ref{fig:7} is in the agreement with Tab. \ref{Table1}.
\par
We estimate the extremal index of the TBTs before and after the PA as well as classes of the newly appearing nodes. 
The intervals and $K-$gaps estimators are used.
 The $K$-gaps  estimator proposed in S\"{u}veges and Davison (2010) 
 may show better accuracy than the intervals estimator according to a simulation study in Ferreira (2018), Markovich and Rodionov (2020b). It has the following form
 \begin{eqnarray*}\label{17a}\widehat{\theta}^K&=&0.5\left((a+b)/c+1-\sqrt{((a+b)/c+1)^2-4b/c}\right),
\end{eqnarray*}
with $a=L-N_C$, $b=2N_C$, $c=\sum_{i=1}^{L}\overline{F}(u_n)S(u_n)^{(K)}_i$. $N_C$ is the number of non-zero $K$-gaps. The $K$-gaps are determined by
 \[S(u)^{(K)}=(\max\left(T(u)-K,0\right)),~~K=0,1,2,...,\] where $T(u)$ is the same r.v. as for the intervals estimator. Both the intervals and $K$-gaps estimators require sufficiently large samples since they are based on inter-exceedance times $\{T(u)_i\}$ which can constitute a very moderate sample.
  In order to estimate the threshold $u$ that is a single parameter of these 
  estimators we use the $\omega^2-$discrepancy method proposed in Markovich and Rodionov (2020b). The latter estimates are denoted as $\widehat{\theta}^{Idis}(n,k)$ and $\widehat{\theta}^{K0dis}(n,k)$. 
  We take $k=\lfloor\hat{\theta}_0L\rfloor$, where $\hat{\theta}_0$ is a pilot intervals estimator. 
  \\
  Using the algorithm in Markovich and Rodionov (2020b, Sect. 4), we find $u$ as solution of the discrepancy inequality in formula (12) by Markovich and Rodionov (2020b) and calculate  \begin{eqnarray}\label{19}
   \widehat{\theta}_1 &=&\frac{1}{l}\sum_{i=1}^l\widehat{\theta}(u_i),\quad \widehat{\theta}_2 =\widehat{\theta}(u_{min})
    \end{eqnarray}
    as resulting estimates of the extremal index, where $u_1,...,u_l$ are  possible solutions of the discrepancy inequality, $l$ is their random number and
        $u_{min}=\min\{u_1,...,u_l\}$. $\widehat{\theta}_1$ and $\widehat{\theta}_2$ are shown in Tab. \ref{Table2}. $\widehat{\theta}_2$ is shown in brackets.
        \\
        To find $u$ for the intervals estimator, we apply also the plateau-finding Algorithm 1 proposed in  Ferreira (2018) denoted as $\widehat{\theta}^{IA1}(n,k)$. To find an optimal pair $(u,K)$ in the $K-$gaps estimator we use 
        the IMT method proposed in Fukutome et al. (2015). This estimate is denoted as $\widehat{\theta}^{KIMT}(n,k)$. All estimates provide 
        similar results, see 
        Tab.\ref{Table2}.  The "out-degree-classes" and the $TBT_2$ after the PA demonstrate  small estimates of the extremal indices except $\widehat{\theta}^{KIMT}(n,k)$. 
        This property may imply  strong local dependence in these data sets. The rest of the TBTs and the "in-degree-classes" have extremal indices close to $1$ that may mean nearly independence. The closeness of the extremal indices of the "in-degree-classes" and the TBTs before the PA is in the agreement with the item (ii) of Theorem \ref{Prop1}. Really, the PRs of "new" nodes belonging to an "in-degree-class" are obtained as sums of the PRs of "old" nodes (of the corresponding TBTs) 
        with in-coming links to "new" nodes.

%% file: realdata1.tex
\subsection{Real data analysis}\label{RealData}
We investigate the Berkeley-Stanford web graph from 2002 
with 685230 nodes and 7600595 edges (see snap.stanford.edu), which represents pages from the \\berkely.edu and stanford.edu domains that are connected in a union network by directed edges as hyperlinks between them (Leskovec et al., 2009).
\begin{figure}[tbp]
    \begin{minipage}[t]{\textwidth}
\centering
  \subfigure[]{\includegraphics[width=0.45\linewidth]{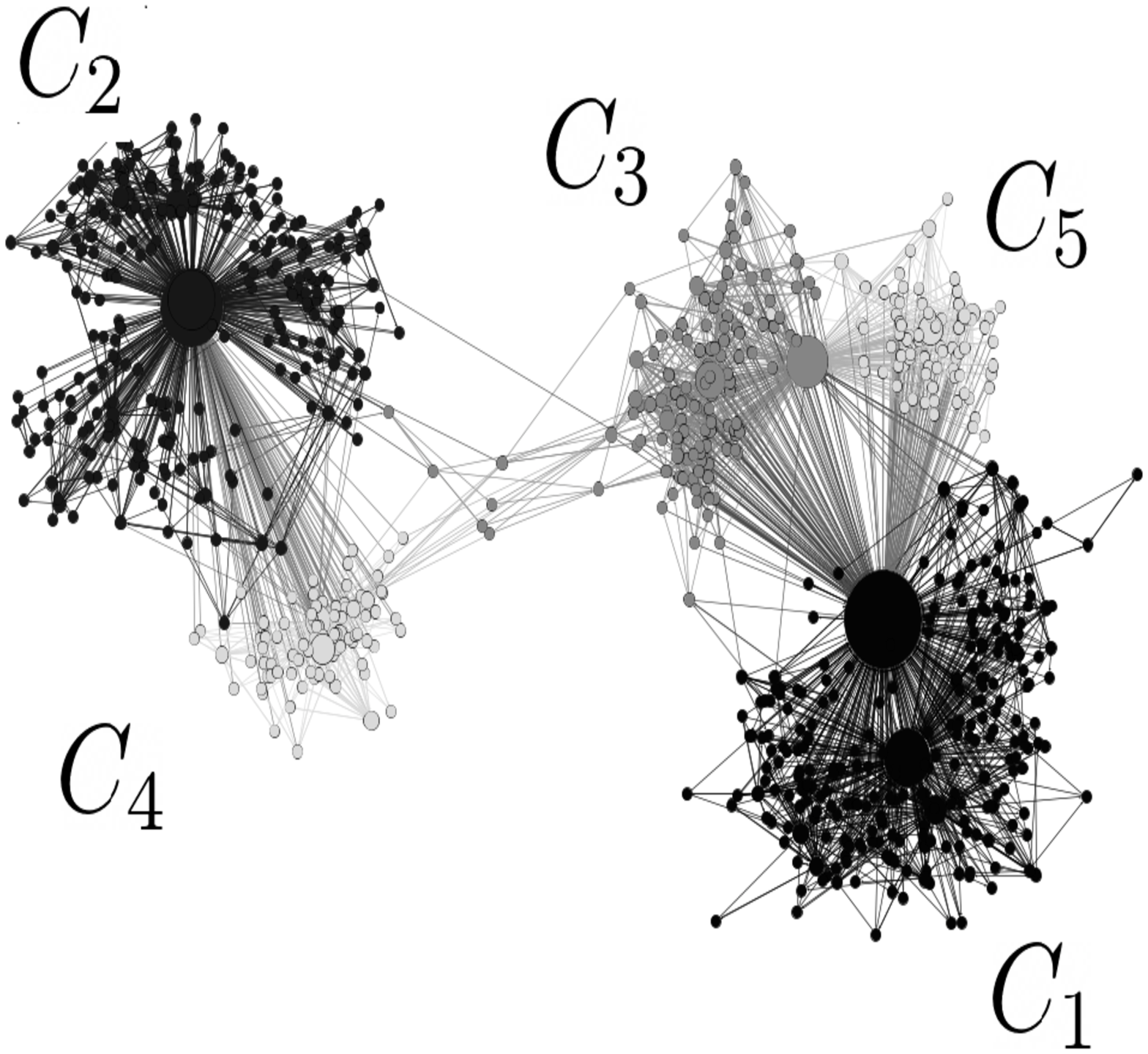}
  \label{fig:8a}}
\subfigure[]{\includegraphics[width=0.45\linewidth]{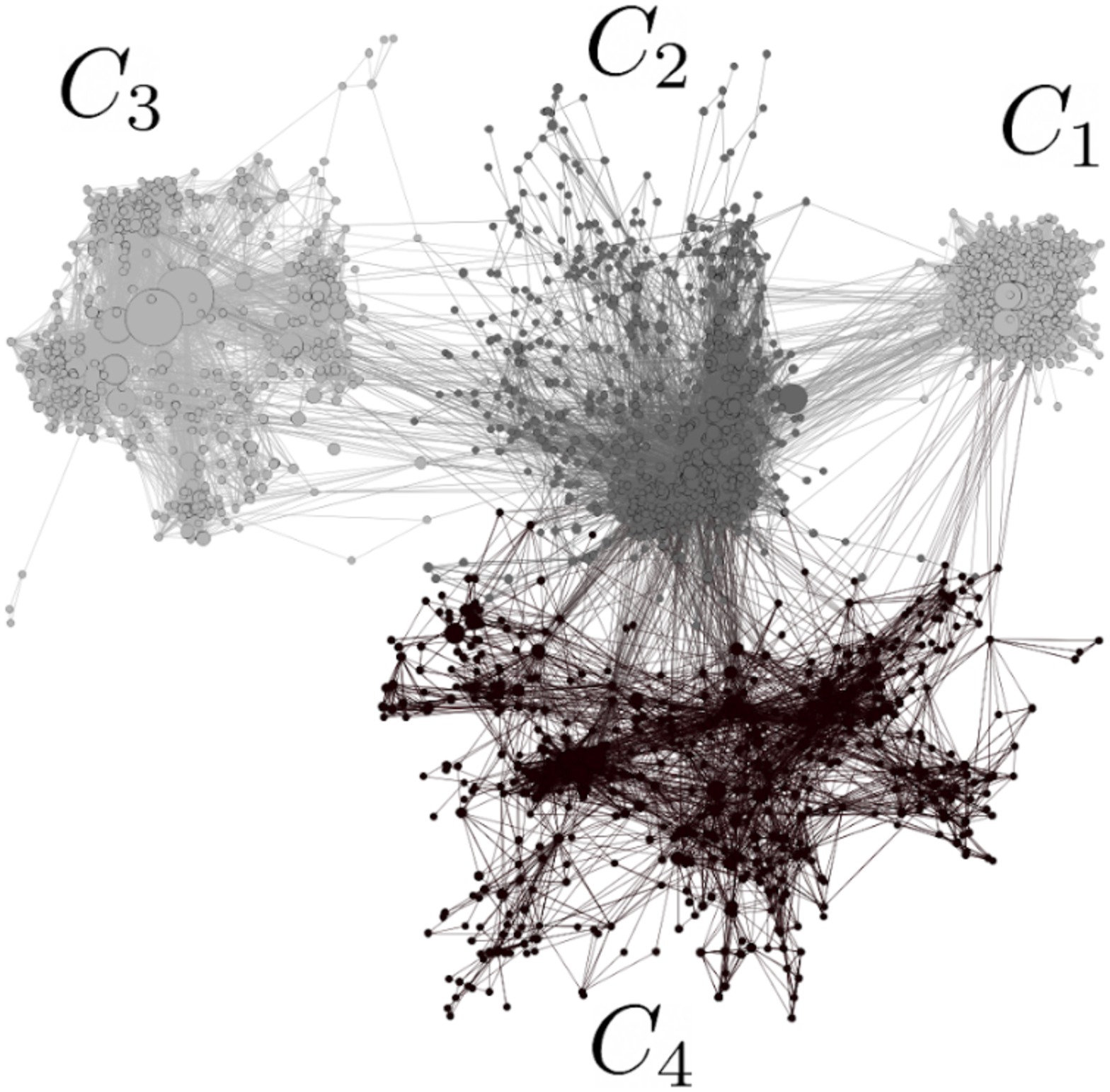}
\label{fig:8b}}
  \caption{Examples of subnetworks of smaller sizes (Fig. \ref{fig:8a}) and larger sizes (Fig. \ref{fig:8b}) divided into communities $\{C_i\}$ 
  received with Directed Louvain's algorithm
  by the Berkeley-Stanford data, where point sizes are proportional to the nodes' PRs.
  }\label{fig:8} 
  \end{minipage}
\end{figure}
\begin{figure}[tbp]
\begin{minipage}[t]{\textwidth}
\centering
\subfigure[]{\includegraphics[width=0.45\textwidth]{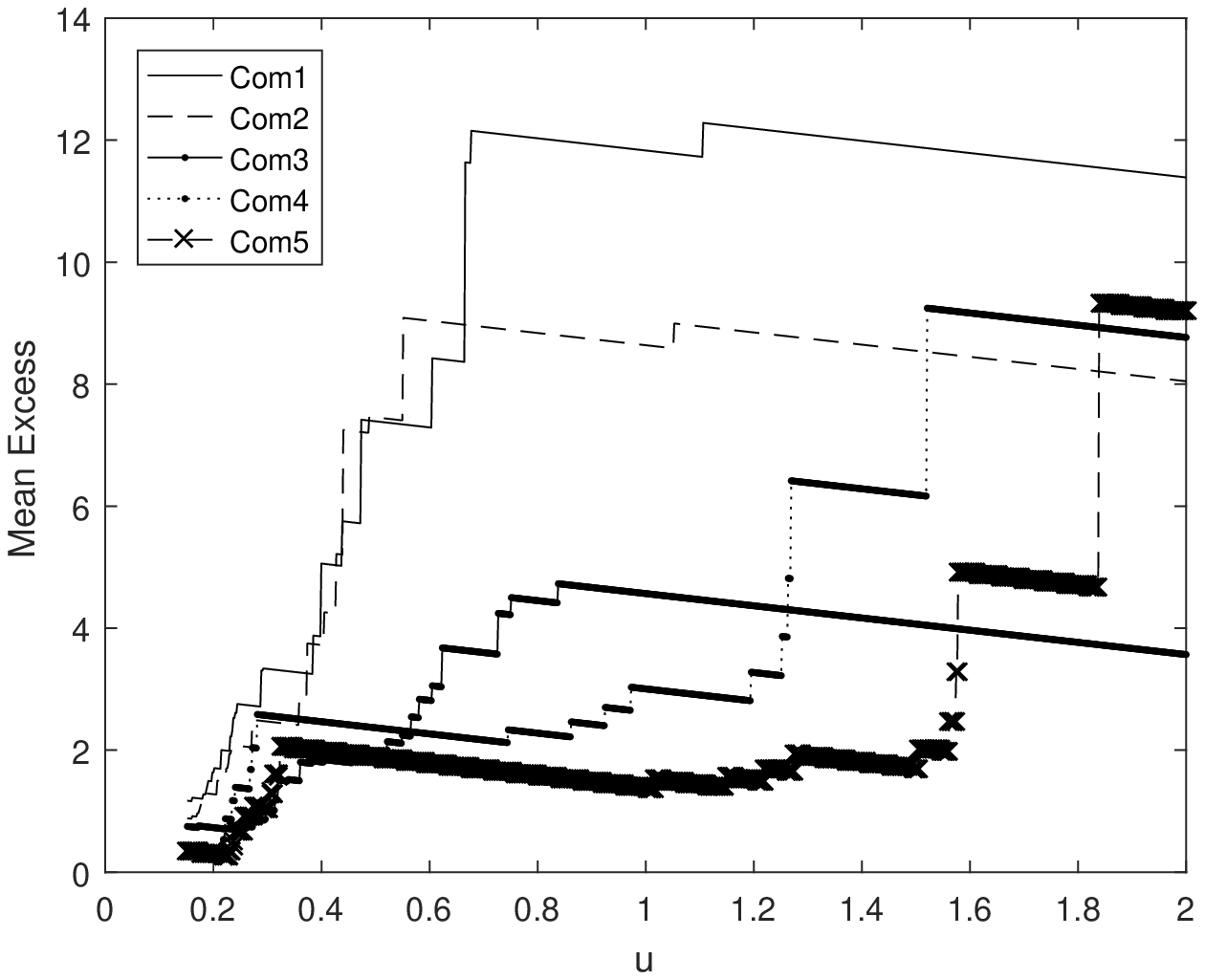}
\label{fig:9a}}
\subfigure[]{\includegraphics[width=0.45\textwidth]{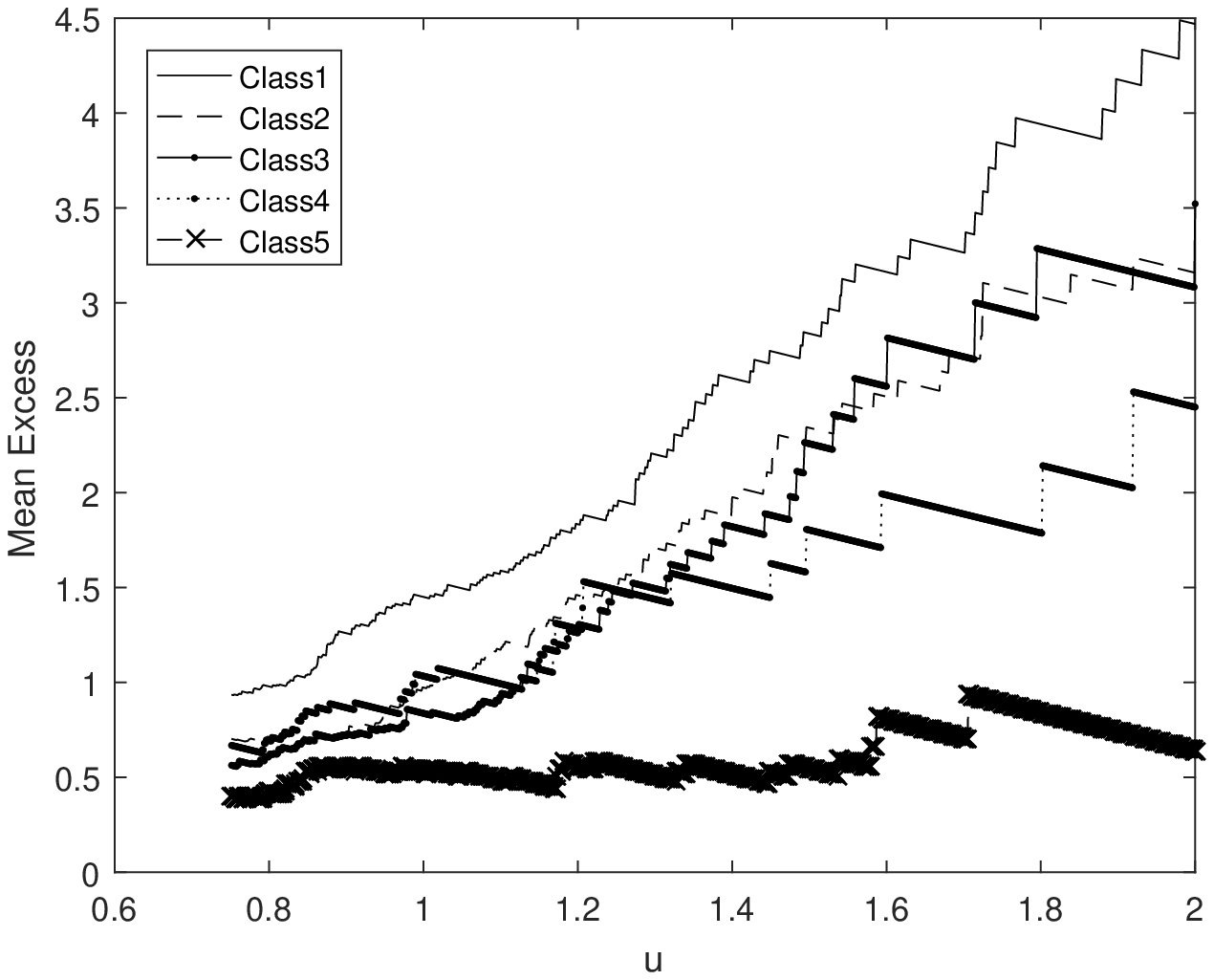}
\label{fig:9b}}
\\
\subfigure[]{\includegraphics[width=0.45\textwidth]{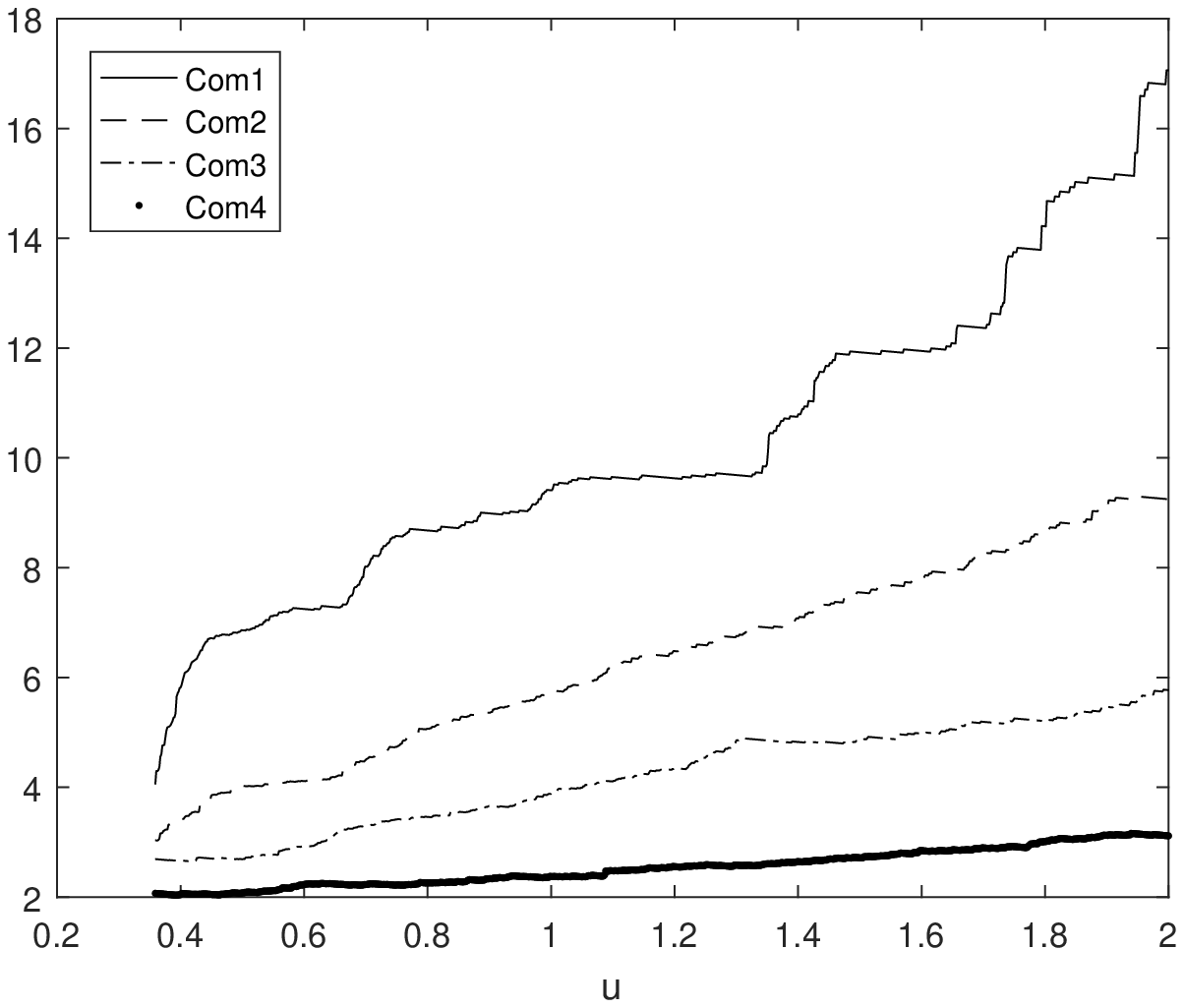} 
\label{fig:9c}}
\subfigure[]{\includegraphics[width=0.45\textwidth]{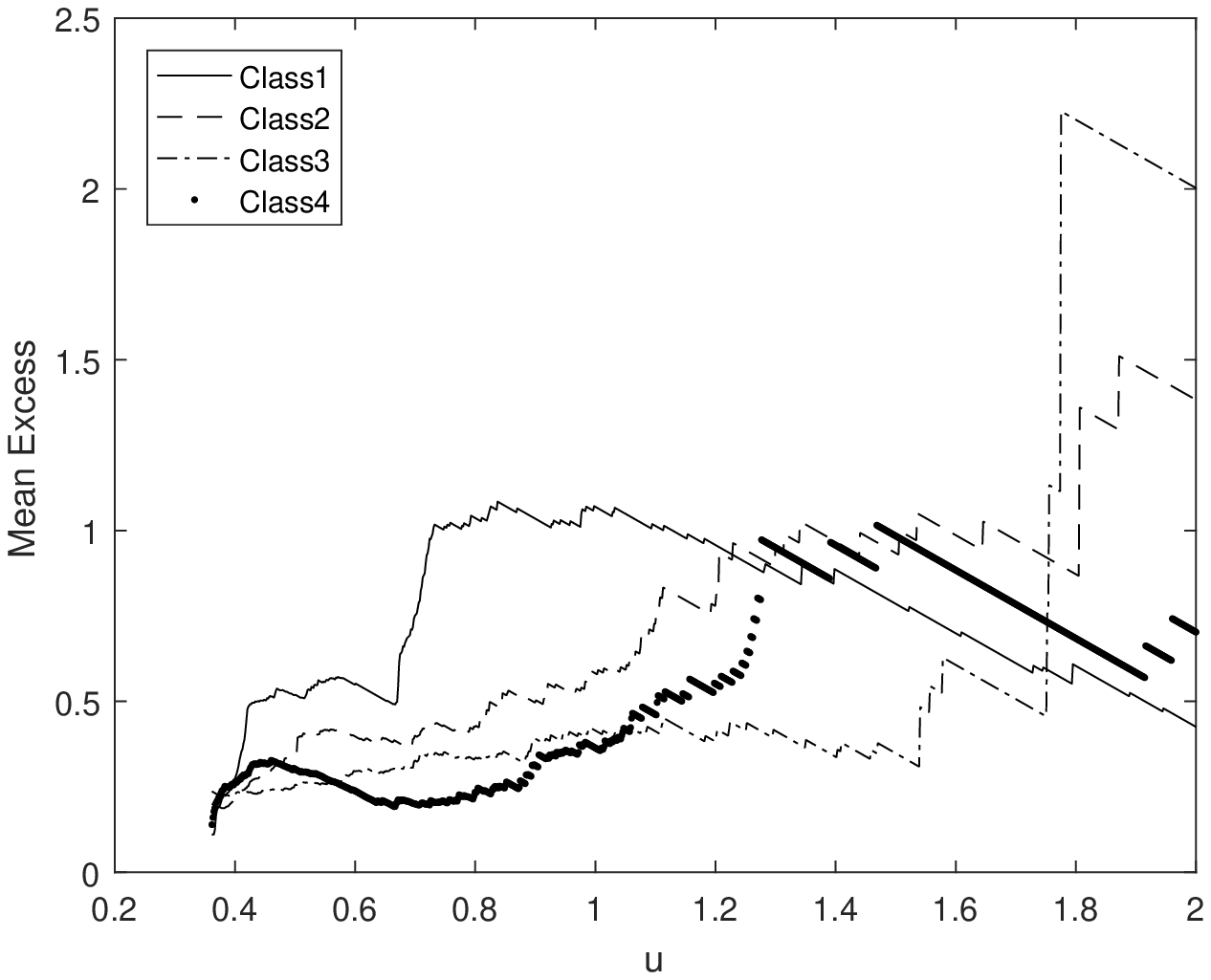} 
\label{fig:9d}}
\caption{Mean excess functions of the communities and the "in-degree-classes" of smaller sizes (upper line) 
and larger sizes (lower line) 
with appended edges obtained by $5000$ and $10^4$ newly appearing nodes, 
respectively.}
\label{fig:9} 
\end{minipage}
\end{figure}
\\
Within this network we select two examples of subnetworks of smaller and larger sizes (Fig. \ref{fig:8}) 
and
partition each of them into $m$ 
communities by
means of the Directed Louvain's algorithm, based on modularity maximization (Dugu\'{e} \& Perez, 2015). 
PRs  of  nodes are calculated by 
(\ref{8}). 
\\
These communities are used as  seed networks. Starting from the seeds we apply the $\beta-$ and $\gamma-$PA schemes to evolve graphs and to obtain "in-degree-classes" $Class_1-Class_m$
by Algorithm \ref{Algor2} (the $\alpha-$ and $\beta-$schemes lead to "out-degree-classes").  The 
$Class_{m+1}$
contains  new nodes which are 
appended to previously appearing new nodes 
but  not to the seed communities. The latter class can be further partition into subclasses by Algorithm \ref{Algor2} with regard to their links to $Class_1-Class_m$ taken now as the seed network.
\\
We check the stationarity of node PRs of the communities and the "in-degree-classes" by the estimation of the mean excess function, Fig. \ref{fig:9}. One may conclude that distributions of all communities and classes belong to the Pareto-type  due to the linear increasing of the mean excess functions. The latter show  the priority of the communities and classes by an ascending of their tail indices due to (\ref{2}). Particularly, the community $C_1$ has the largest $\gamma$ and thus, the smallest tail index $\alpha$ as one can see in Tab. \ref{Table3}, too. 
\\
\begin{figure}[tbp]
\begin{minipage}[t]{\textwidth}
\centering
\subfigure[]{\includegraphics[width=0.23\textwidth]{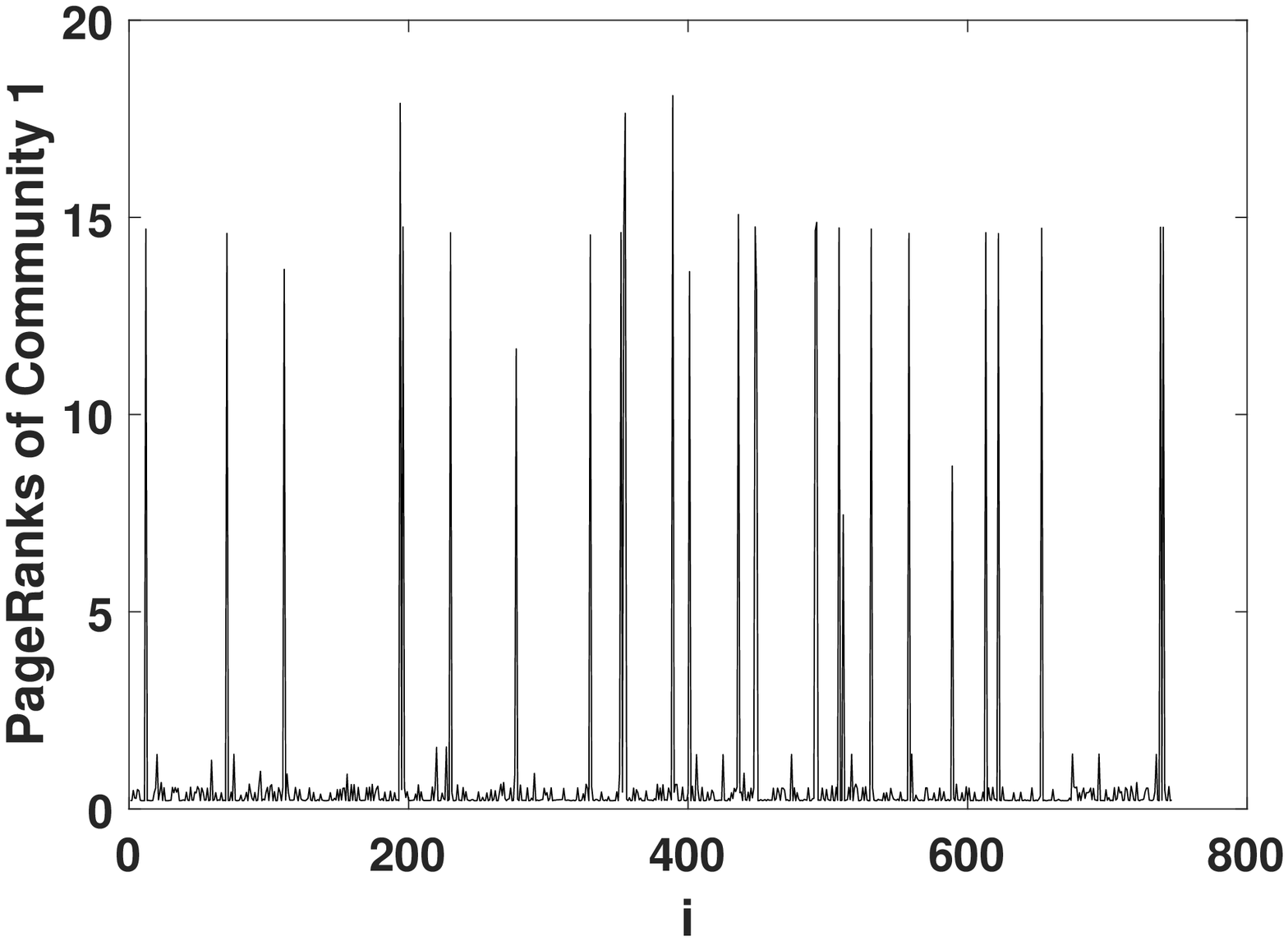} 
\label{fig:10a}}
\subfigure[]{\includegraphics[width=0.23\textwidth]{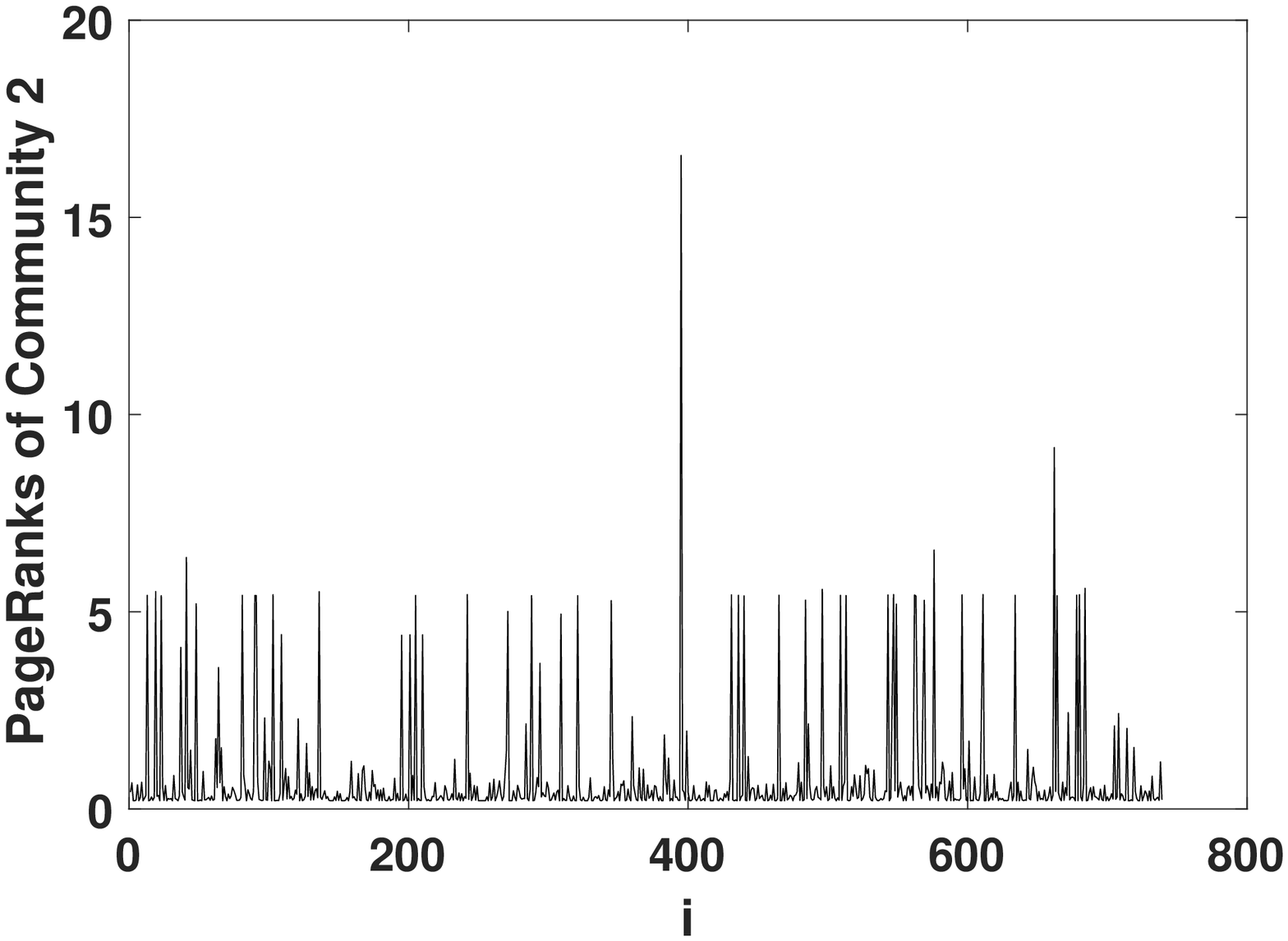} 
\label{fig:10b}}
\subfigure[]{\includegraphics[width=0.23\textwidth]{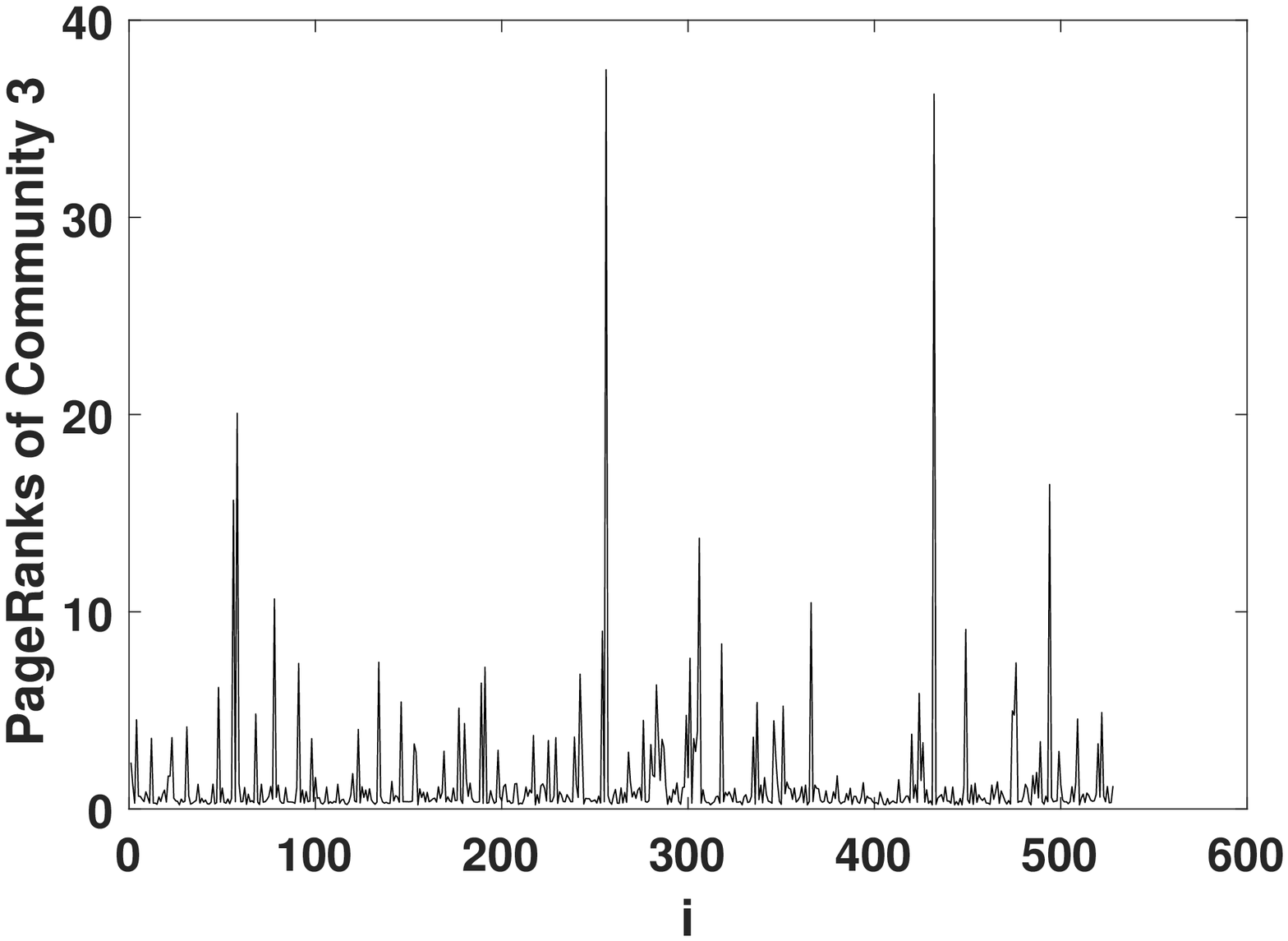} 
\label{fig:10c}}
\subfigure[]{\includegraphics[width=0.23\textwidth]{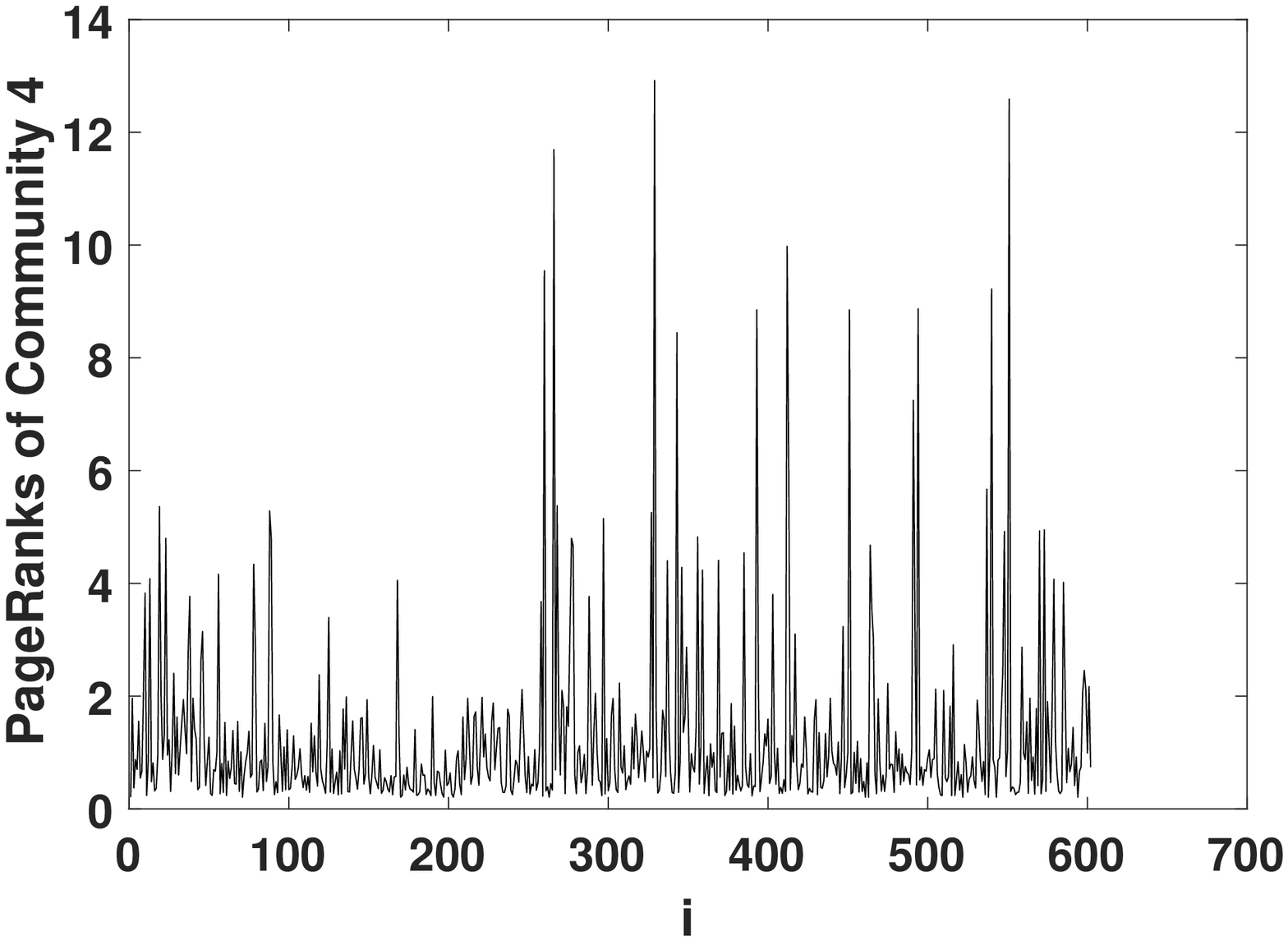} 
\label{fig:10d}}
\\
\subfigure[]{\includegraphics[width=0.23\textwidth]{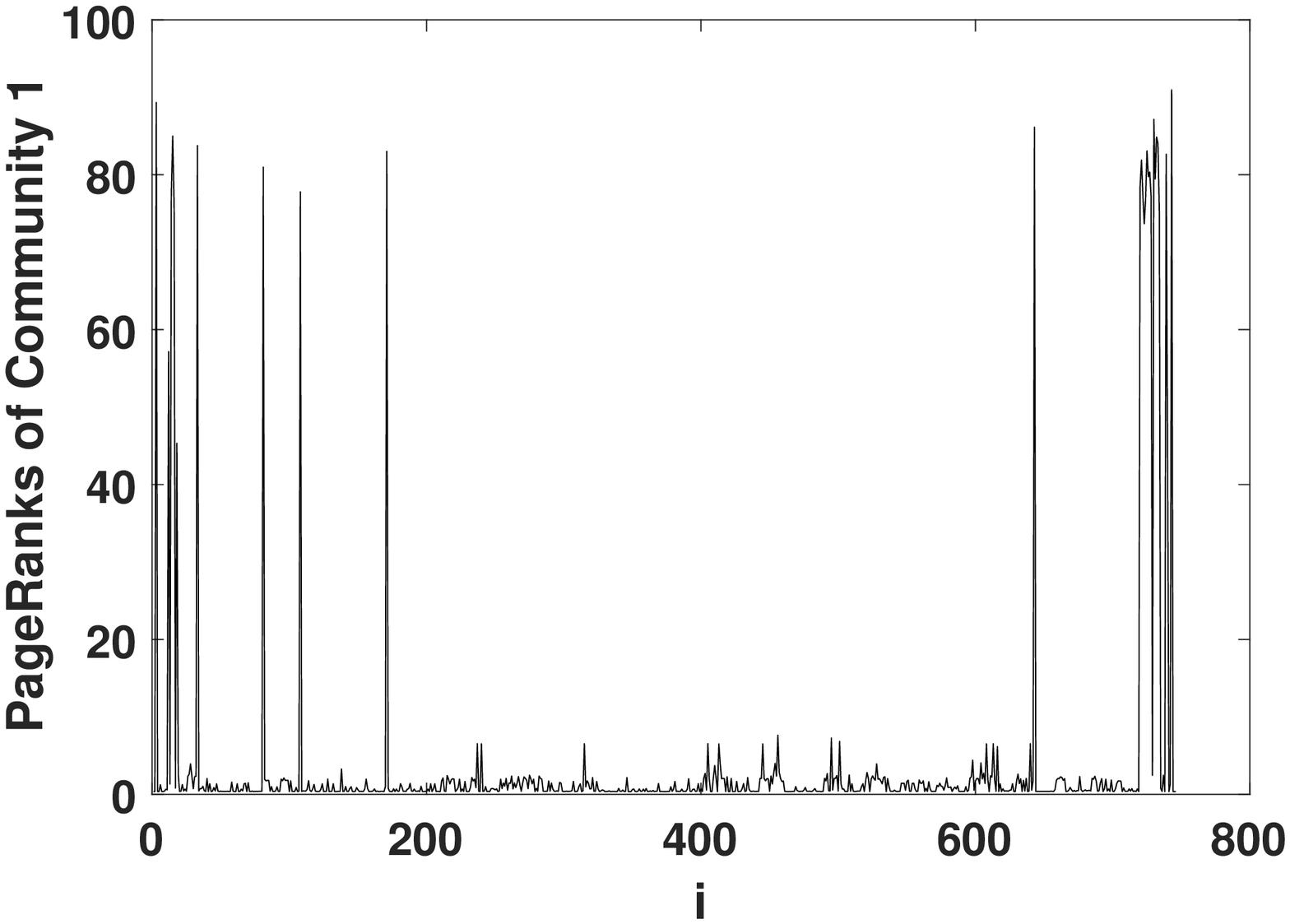} 
\label{fig:10e}}
\subfigure[]{\includegraphics[width=0.23\textwidth]{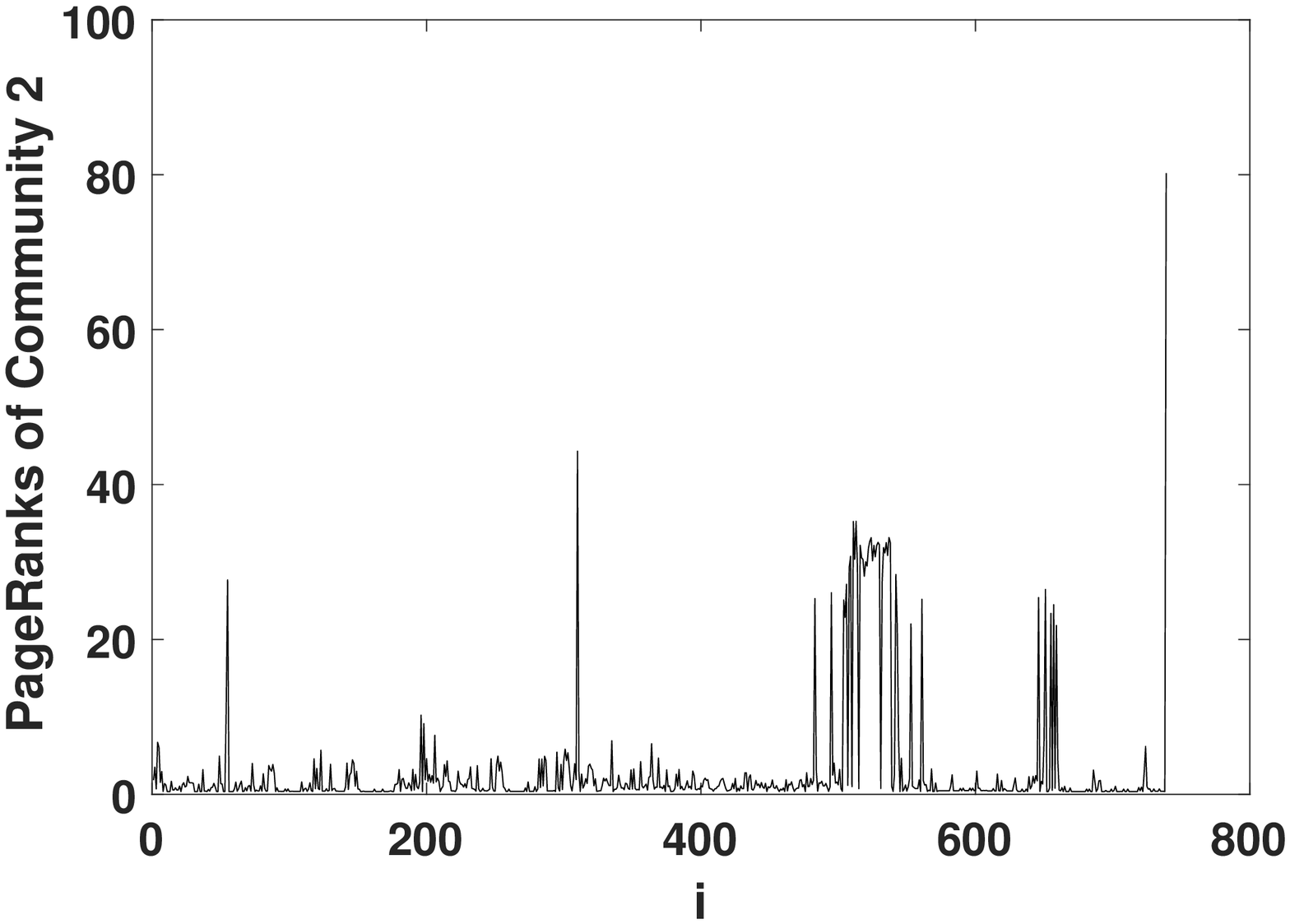} 
\label{fig:10f}}
\subfigure[]{\includegraphics[width=0.23\textwidth]{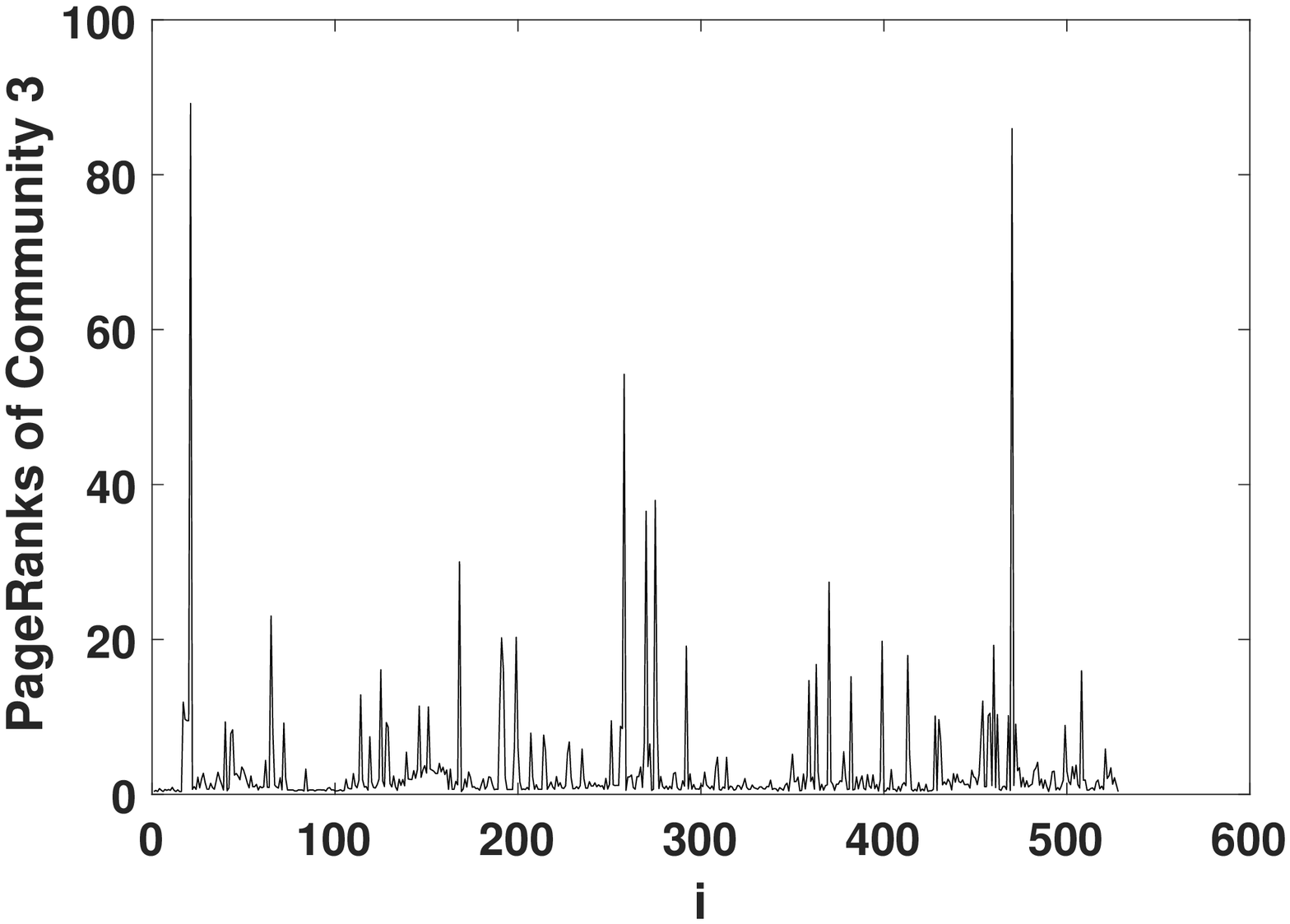} 
\label{fig:10g}}
\subfigure[]{\includegraphics[width=0.23\textwidth]{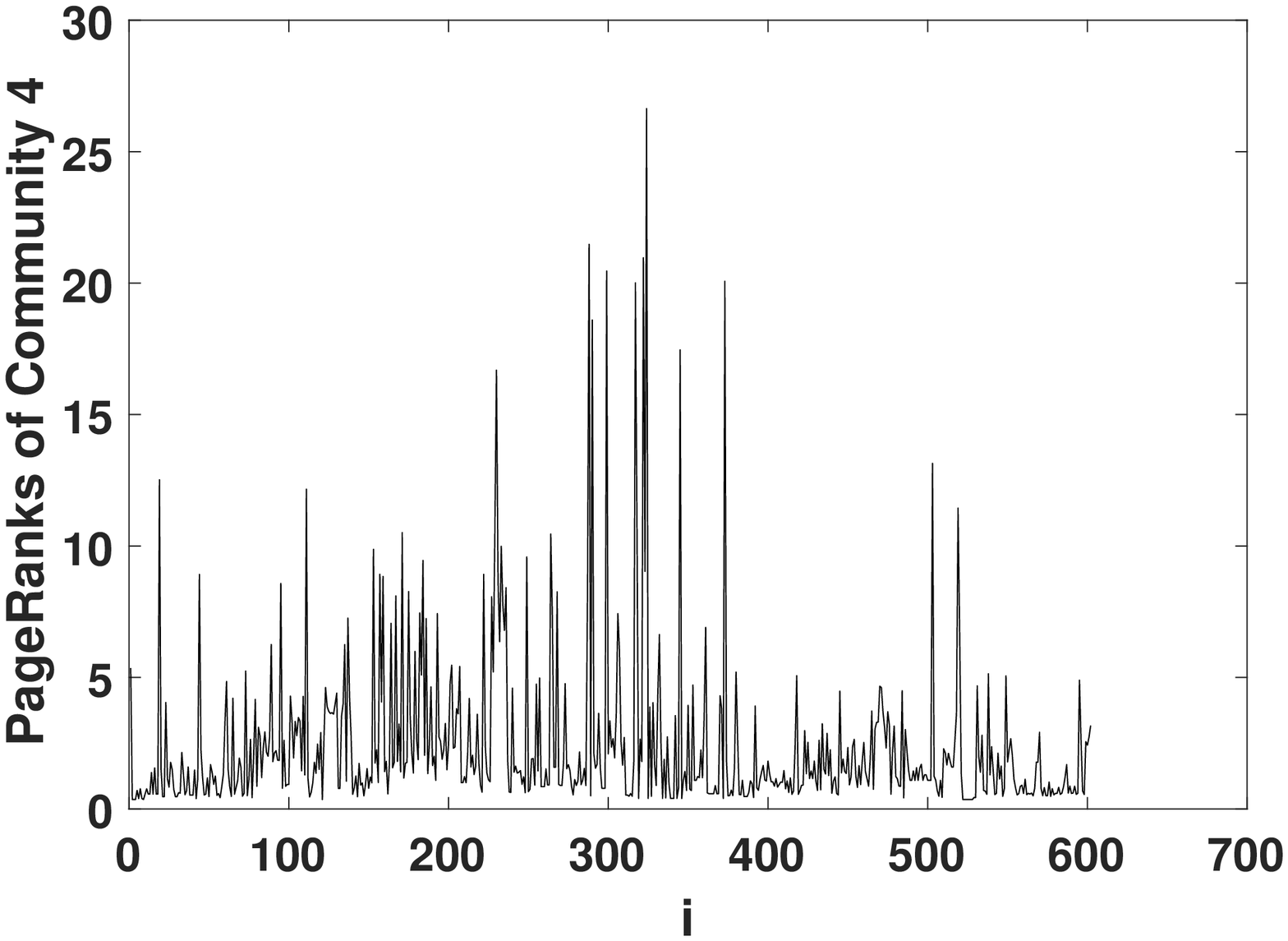} 
\label{fig:10h}}
\caption{The PRs of the "large" communities $1-4$ before (upper line from left to right) and after (lower line from left to right) the $PA(0.4,0.2,0.4)$ with  $\delta_{in}=\delta_{out}=1$ of $10^4$ new nodes.}
\label{fig:10} 
\end{minipage}
\end{figure}
\begin{table}[t]
\caption{The Hill's estimates $\widehat{\alpha}(n,k_b)$ and $\widehat{\alpha}^D(n,k_b)$ of the PR tail indices
 of "old" nodes in two sets of communities with appended new edges to the $N_0$ PA-appended "new" nodes 
 and of new nodes 
 in "in-degree-classes"
 with the $95\%$ bootstrap confidence intervals
 $(u_1,u_2)$ obtained by $5000$ bootstrap resamples, where the  $PA(0.4, 0.2, 0.4)$ with $\delta_{in} = \delta_{out} = 1$ is used. 
 }
 \centering
\tabcolsep=0.001cm
\begin{tabular}{lccc||ccccc}
  \hline

       $Commu$   & $n$ &  $~~\widehat{\alpha}(n,k_b)$ & $~~\widehat{\alpha}^D(n,k_b)$ && $Class$& $n$ &  $~~\widehat{\alpha}(n,k_b)$ & $~~\widehat{\alpha}^D(n,k_b)$ \\
       $nity$           &                          & $(u_1,u_2)$             & $(u_1,u_2)$                     &   &     & & $(u_1,u_2)$ & $(u_1,u_2)$\\\hline
              &                            &                                                        \multicolumn{7}{|c|}{\scriptsize {$N_0=5000$}}\\
     \cline{3-9}
       $1$ & 266 &  1.1782           & 1.2677 && $1$ & 512 &  1.1146 & 1.5520 \\
               &     & $(0.3647,1.9481)$ & $(0.3719,1.9921)$                     &   &     & &$(0.7713,1.7581)$ & $(0.7708,1.7741)$
       \\
      $2$ & 266 &  1.2582           & $1.3020$ && $2$ & 466 &  1.2280 & 1.4946\\
              &     & $(0.3886,2.1817)$ & $(0.3776,2.1500)$                     &   &     & & $(0.7804,1.8732)$& $(0.7613,1.8124)$
       \\
      $3$ & 135 &  1.8228           & 2.6509               && $3$ & 327 &  1.6247 &1.4708 \\
              &     & $(0.6234,9.2239)$ & $(0.6965,8.7533)$                     &   &     & &$(1.3935,2.7271)$ & $(0.7959,3.0236)$
       \\
      $4$ & 86 &  1.2493           & 1.0752 && $4$ & 182 &  1.8196 & 1.7822\\
              &    & $(0.5661,5.0906)$ & $(0.5725,5.0560)$                     &   &     & &$(0.7677,2.9515)$ & $(0.7671,3.0409)$
       \\
      $5$ & 85 &  1.9378             & 1.4255              && $5$ & 166 &  2.8095 & 3.6141\\
              &     & $(0.5846,14.1712)$ & $(0.5855,12.4387)$                     &   &     & & $(2.2747,3.3840)$ & $(1.6424,4.2651)$
       \\
              &    &         &                 && $6$ & 3347 &   &  \\
              \hline
                  &                            &                                                          \multicolumn{7}{|c|}{\scriptsize {$N_0=10000$}}\\
     \cline{3-9}
       $1$ & 746 &  1.1098  &  0.7504  && $1$ & 2035 &  1.7222 & 1.2845 
       \\
               &     & $(0.5333,1.6993)$ & $(0.6681,0.7808)$ &   &     & &$(1.5360,1.9250)$ & $(1.1958,1.3732)$
       \\
      $2$ & 739 &  0.8924 & 0.9440 && $2$ & 1394 &  3.2711 & 3.4538\\
              &     & $(0.4506,1.3905)$ & $(0.7743,0.9987)$                     &   &     & & $(2.8033,3.8688)$ & $(3.4520,3.4557)$ 
       \\
      $3$ & 528 &  1.1284 & 1.2029&& $3$ & 635 &  2.1855 & 2.0486
      \\
              &     & $(0.9853,1.3426)$ & $(0.9723,1.3260)$                     &   &     & &$(1.4684,3.6417)$ & $(2.0231,2.0740)$
      \\
      $4$ & 602 &  2.3255 & 3.0869 && $4$ & 480 &  2.9701 & 3.1583\\
             &     & $(1.6265,3.2978)$ & $(1.6220,3.3123)$                     &   &     & & $(2.2507,3.2001)$ & $(3.1430,3.1736)$
       \\
       &    &         &                        && $5$ & 5456 &   &  \\
                   \hline
    \end{tabular}
    \label{Table3}
\end{table}
\begin{figure}[tbp]
 \begin{minipage}[t]{\textwidth}
\centering
\subfigure[]{\includegraphics[width=0.45\textwidth]{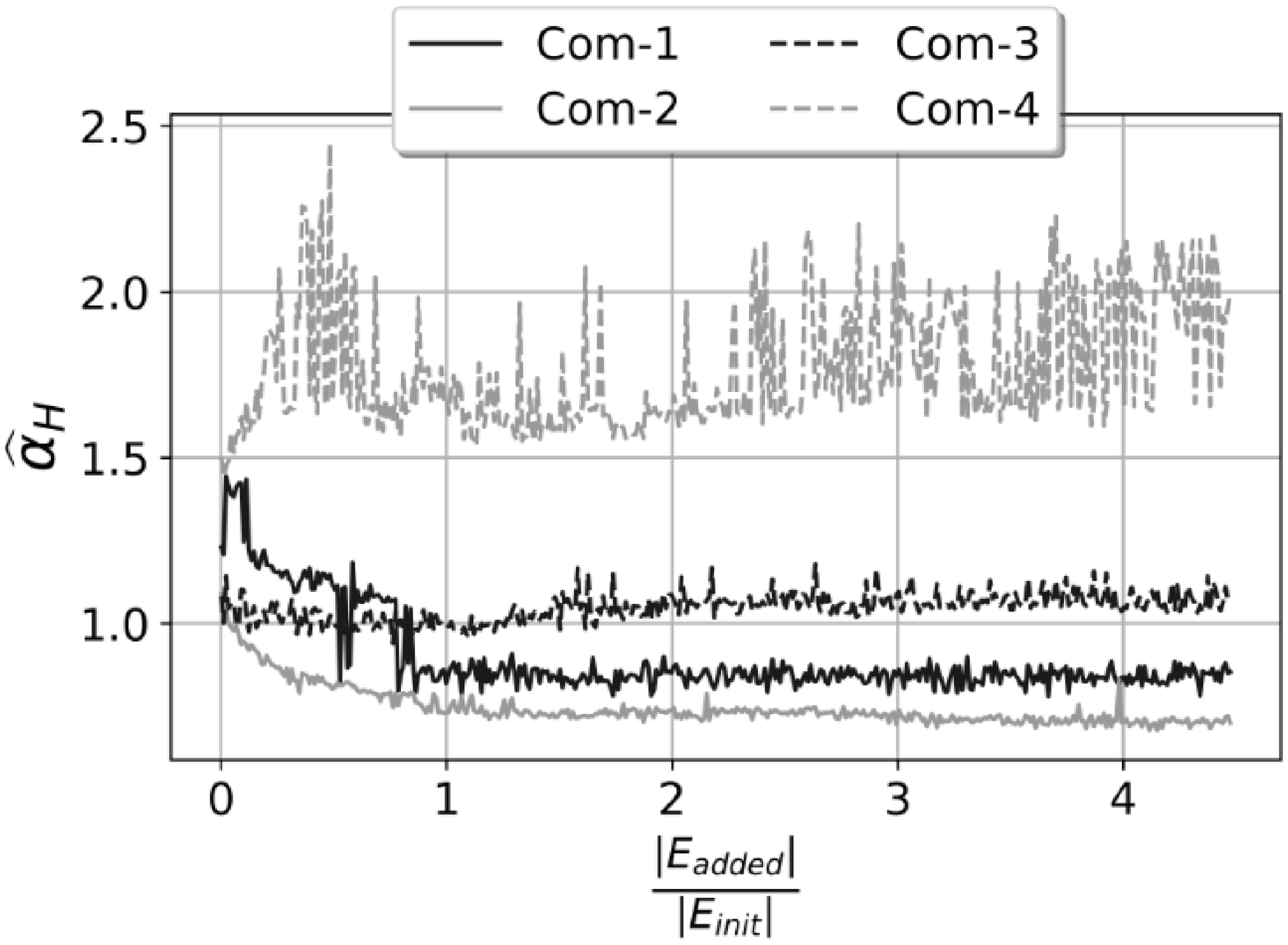} 
\label{fig:11a}}
\subfigure[]{\includegraphics[width=0.45\textwidth]{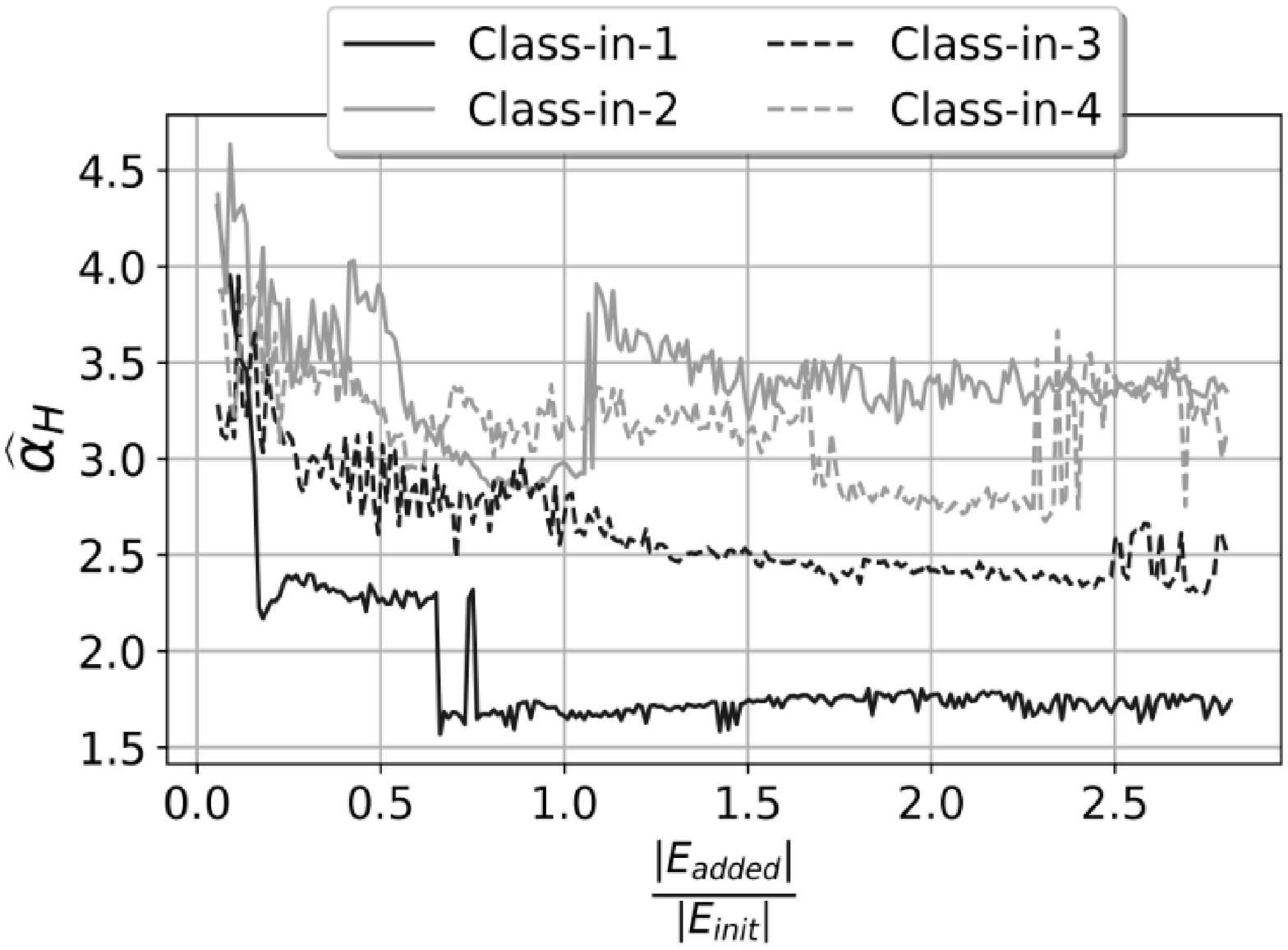}
\label{fig:11b}}
  \caption{The Hill's estimates $\widehat{\alpha}(n,k)$ of PR of the large size communities (Fig. \ref{fig:11a}) and   in-degree classes (Fig. \ref{fig:11b})  with $k$ selected by bootstrap method versus $|E|_{added}/|E|_{init}$, where $|E|_{added}$ is the number of edges added to the  communities  during the PA
and the initial number of edges $|E|_{init}$ in   the communities.
  }
\label{fig:11}
\end{minipage}
\end{figure}
Fig. \ref{fig:10} shows that the cluster structure of PRs before and after the attachment of new edges is different. The PRs are enumerated 
with regard to the node appearance in the dataset and the sequence of further attachment. The appearance of higher peaks after the PA corresponds to giant nodes with large PRs. The tail indices of the PRs after the PA may only decrease.
\\
In Tab. \ref{Table3} the Hill's estimates $\widehat{\alpha}(n,k_b)$ and $\widehat{\alpha}^D(n,k_b)$
of the PRs of the communities of both smaller and larger sizes after the appending of $N_0$ new nodes and the "in-degree-classes" of appended new nodes are shown. Due to a possible non-stationarity the tail index of the rest 
$Class_{m+1}$ is not estimated. The number of largest order statistics $k_b$ of $\widehat{\alpha}(n,k_b)$ is obtained by a bootstrap method and of $\widehat{\alpha}^D(n,k_b)$ by a double bootstrap method (see Markovich (2007) for details). The tail indices of the in-degree-classes are larger than ones of the evolved communities and correspond to the tail indices of the communities before the evolution.
The drop and stabilization of the tail indices of the communities and "in-degree-classes" is shown in Fig. \ref{fig:11} when the number of new edges increases.
\begin{table}[t]
\caption{Estimation  of the PR extremal index of the
 "old" nodes in the  "large" communities   and  the "in-degree-classes" of $N_0$ appended new nodes.
  }
 \centering
\tabcolsep=0.01cm
\begin{tabular}{lcccc||ccccc}
  \hline
   $Commu$  & $~~\widehat{\theta}^{IA1}$  & $~~\widehat{\theta}^{Idis}$ & $~~\widehat{\theta}^{K0dis}$& $~~\widehat{\theta}^{KIMT}$ & $Class$ &$~~\widehat{\theta}^{IA1}$ & $~~\widehat{\theta}^{Idis}$ & $~~\widehat{\theta}^{K0dis}$ & $~~\widehat{\theta}^{KIMT}$\\
     $nity$   &    &  &  &  &   &  & &  & \\
   \hline
   \multicolumn{5}{c}{\scriptsize {Before PA}}&  \multicolumn{5}{c}{\scriptsize {"in-degree-classes"}}\\
     \hline
     $1$    & 0.9983   &~ 0.9781 
     & 0.9842 & 0.8991 &  $1$ & 0.9356 & 0.9419 & 0.9475 & 0.8851\\
             &          & (0.9449) & (0.9455) & &        &   & (0.8033)    & (0.8065)
  \\
  $2$    & 0.9919   & 0.9874 
  &  0.9912 & 0.9124 & $2$ & 0.8924 & 0.9245 & 0.9284 & 0.8951\\
             &          &  (0.9249)  & (0.9324)&  &           &        & (0.8595)    & (0.8643) &
   \\
   $3$   & 1    & 0.9988 
   &     1   & 0.9278  & $3$    & 0.9583 & 0.9737 &  0.9745 &  0.9379  \\
             &      & (1) & (1)     &   &              &        & (1)    & (1)&
              \\
  $4$    & 0.6138        &  ~0.7001 
  & 0.7081  & 0.8399 & $4$ & 0.8982 & 0.9357 & 0.9424 & 1\\
             &          &  (0.6179)    & (0.6218)&  &           &        & (1)    & (1) &
   \\
  \hline
  \multicolumn{5}{c}{\scriptsize {After PA of $N_0=10^4$ nodes}} & \\
   \hline
     $1$    & 0.4877   & 0.4756 
                                        & 0.4684 & 0.6823 &   &  &  &  & \\
             &          & (0.4817) & (0.4785) & &        &   &     &
  \\
  $2$    & 0.2538   &   0.3409    
                                         &  0.3528 & 0.5523 &  &  &  &  & \\
             &          &  (0.3362)   & (0.3322)&  &           &        &     &  &
   \\
   $3$   & 0.8638    &  0.8957 
                               &     0.9039   & 0.7484  &     &  &  &   &    \\
             &      &       (0.6771) & (0.6947)     &   &              &        &     & &
              \\
  $4$    & 0.5185        &   0.6716 
                                              & 0.6769  & 0.7364 &  &  &  &  & \\
             &          &   (0.5985)   & (0.6193)&  &           &        &     &  &
     \\
  \hline
  \end{tabular}
 \label{Table5}
\end{table}
In Tab. 
\ref{Table5}  estimates (\ref{19}) of the extremal index are shown when the evolution starts from the 
 "large" communities by the $PA(0.4, 0.2, 0.4)$ with $\delta_{in} = \delta_{out} = 1$. The same estimators of the extremal index as in Tab. \ref{Table2} are used.
 The evolution from the "small" communities is not considered due to a lack of data.
The extremal indices of the "in-degree-classes" are close to those of the communities before the PA. After the PA of new edges the extremal indices of the communities are decreased which implies the increasing of their clustering. 

%% file: conclusions1.tex
\section{Conclusion}\label{Sec4}
\par
The prediction of the tail and extremal indices of node influence characteristics of an evolving network is studied. Assuming that the network can be partitioned into stationary distributed communities of nodes, we classify the newly appended nodes according to their edges to the communities. 
Ranking the communities by their tail indices in ascending order, one can select the most heavy-tailed ("dominating") community with a minimum tail index.
We assign a set of new nodes to the first class with  the tail and extremal indices of the latter community if each node of this set has at least one edge to nodes of the "dominating" community and the latter  is unique. In the next step, we repeat the procedure with the rest of the communities finding the "dominating" one among them and classifying the rest of the newly appearing nodes. Clearly, the procedure can be done in discrete time moments.
The same procedure can be applied if there are a random number of "dominating" communities since the communities of the seed network are independent or weak dependent due to a few links between them.
\par
The assumption regarding the uniqueness of the community with a minimum tail index is plausible and not restrictive for graphs since one has anyway to estimate the tail indices of PRs. The tail index estimates are likely different. This uniqueness property simplifies the analysis. Since the most heavy-tailed "column" series is likely unique, the checking of the homogeneous pair-wise dependence between the components of the most heavy-tailed communities required for an application of Theorems \ref{Prop1} and \ref{T3} and thus, for a prediction of the tail and extremal indices of new classes can be omitted. In case, the community with the minimum tail index is not unique, the pair-wise dependence can be investigated as in Appendix \ref{Sec5.3}. In fact, the "dominating" community with the largest maximum PR determines the extremal index of the PR and MLM of the set of newly appended nodes.
\par
There are some problems of the graph analysis that make it complicated and rough: (a) nodes are not enumerated; (b) the dependence between nodes may be complex and non-homogeneous; (c) the node characteristics in  communities may be non-stationary distributed.
The estimation of the tail index is based on the largest order statistics of the sample and hence, it does not require any enumeration of the nodes. In contrast, the estimation of the extremal index depends on the node enumeration.
\par
To estimate the extremal index we apply the intervals and $K-$gaps estimators. 
To modify the intervals estimator for graphs, we propose to take the number of nodes at the paths between the nodes with exceedances of some feature (e.g., PR) over a sufficiently high threshold as the inter-exceedance times. The node attachment like the preferential attachment provides a natural enumeration of nodes. Thus, one can 
use the 
intervals estimator for random sequences of the PRs or the MLMs for newly appended classes of nodes of evolving networks. 
\\
The stationarity of the communities is proposed to be checked by  the mean excess function that does not require the enumeration of nodes.

%% file: appendix.tex
\appendix
\section{The PageRank and Max-Linear Model}\label{Sec5.1}
\par PR (Langville \& Meyer, 2006) and the MLM  may be considered as node influence 
characteristics (Gissibl \& Kl\"{u}ppelberg, 2018; Markovich et al., 2017). The PR $R$ of a randomly chosen Web page (a node in the Web graph) is viewed as a r.v.. It was considered as the solution to the fixed-point problem
\begin{equation}\label{6} R=^D\sum_{j=1}^{N}A_{j}R_{j}+Q
 \end{equation}
in Jelenkovic and Olvera-Cravioto (2010), Volkovich and Litvak (2010). 
$=^D$ denotes equality in distribution. 
The r.v.s $\{R_j\}$ are assumed to be 
iid copies of $R$ and $E(Q)<1$ holds. $(Q,N,\{A_j\})$ is a real-valued vector, $\{A_j\}$ are independent non-negative r.v.s. distributed as some r.v. $A$ with $E(A)<1$.  $N$ denotes the in-degree of a node. $Q$ is a personalization value of the vertex. Under the assumptions (we shall call them  Assumptions A)
 that  $\{R_j\}$ are regularly varying iid and independent of $(Q, N, \{A_j\})$ with $\{A_j\}$ independent of $(N,Q)$, $N$ is regularly varying r.v., and that $N$ and $Q$ are allowed to be dependent, it is stated in Jelenkovic and Olvera-Cravioto (2010), Volkovich and Litvak (2010) that the stationary distribution of $R$ in (\ref{6}) is regularly varying and its tail index is determined by the most heavy-tailed distributed term in the pair $(N, Q)$. The approach implicitly assumes that the underlying graph is an infinite tree, an assumption that is not plausible in real-world networks. In Chen et al. (2014), the behavior of the PR is considered on a directed configuration model, which is a tree-like graph in a sense that the first loop is observed at a distance of order $\log n$, where $n$ is the size of the
graph. It is derived that the PR in the latter model is well approximated by the PR of the root node of a suitably
constructed tree as $n\to\infty$.
\\
In the same way, a MLM is considered as the 'minimal/endogeneous' solution   of the  equation
  \begin{equation}\label{6a}R=^D\left(\bigvee_{j=1}^{N}A_{j}R_j\right)\vee Q,\end{equation}
(Jelenkovic \& Olvera-Cravioto, 2015). Assuming that all r.v.s in the triple $(R_j,Q,N)$ are regularly varying and mutually independent and $\{R_j\}$ are iid, PR $R$ was proved to have a regularly varying tail in Jelenkovic and Olvera-Cravioto (2010), Volkovich and Litvak (2010) at the directed configuration model. A similar statement was proved in Jelenkovic and Olvera-Cravioto (2015) with regard to the MLM.

\section{Important results from extreme value analysis}\label{Sec5.2}
\par We recall the theorems derived in Markovich (2022) that are important for the prediction of the tail and extremal indices of evolving random graphs. These theorems generalize Theorems 3 and 4 in Markovich and Rodionov (2020a). The latter state the conditions when the sequences of sums $Y_{n}(z, N_n)$ and maxima $Y_{n}^*(z, N_n)$ (see (\ref{3})) have the same tail  and  extremal indices. There are the following constrains for these statements. 
The slowly varying functions $\{\ell_i(x)\}$ in (\ref{11a}) are uniformly upper bounded in $i$ by a polynomial function, i.e. for all constants
$A>1$, $\delta>0$
there exists $x_0(A, \delta)$ such that for all $i\geq 1$
\begin{eqnarray} \ell_i(x)\leq A x^\delta, \qquad  x>x_0(A, \delta) \label{uniform}
\end{eqnarray}
holds. Despite $N_n$ is integer-valued, one can accept a distribution with regularly varying tail with tail index $\alpha>0$ as a relevant model for $N_n$, i.e. it holds
 \begin{eqnarray}\label{15a} &&P(N_n>x) = x^{-\alpha} \tilde{\ell}_n(x),\end{eqnarray}
 $\tilde{\ell}_n(x)$ is a slowly varying function.
This model is motivated in several papers, see Jessen and Mikosch (2006), Robert and Segers (2008) and Volkovich and Litvak (2010) among them.
\par In Markovich and Rodionov (2020a) it is assumed that there is a unique "column" sequence with a minimum tail index $k_1<k$, 
$k:= \lim_{n\to\infty} \inf_{2\leq i\leq l_n} k_i$,  and $N_n$ has a lighter tail than $Y_{n,i}$, i.e. it holds
\begin{eqnarray}\label{4a}P\{N_n>l_n\}&=&o\left(P\{Y_{n,1}>u_n\}\right), ~~n\to\infty,\end{eqnarray}
 where
 the sequence of thresholds $u_n$ is taken as 
 $u_n=yn^{1/k_1}\ell_1^{\sharp}(n)$, $y>0$, $\ell^{\sharp}(x)$ is the de Brujin conjugate of $\ell(x)$,
 the sequence $l_n$ satisfies
 \begin{eqnarray}\label{27}&& l_n=[n^\chi],
  \qquad
   \end{eqnarray}
  and  $\chi$ satisfies
    \begin{equation}
0<\chi<\chi_0, \qquad\chi_0 = \frac{k-k_1}{k_1(k+1)}.
\label{chi}
\end{equation}
An arbitrary dependence between "column" sequences and between $\{Y_{n,i}\}$ and $\{N_n\}$ is allowed. In Theorem 4 in Markovich (2022) recalled here in Theorem \ref{T3}  the number $d$ of "column" sequences with a minimum tail index is allowed to be random. This is realistic for random graphs since a random number of communities considering as the "column" sequences may have a minimum tail index (or a tail index close to that).
\par
Let us recall  the following conditions for a fixed $d>1$ proposed in Markovich (2022).
\begin{enumerate}
\item[(A1)] The stationary sequences $\{Y_{n,i}\}_{n\ge 1}$, $i\in\{1,...,d\}$ are mutually independent, and  independent of the sequences
$\{Y_{n,i}\}_{n\ge 1}$, $i\in\{d+1,...,l_n\}$.
\\
\item[(A2)] Assume $\{Y_{n,i}\}_{n\ge 1}$, $i\in\{1,...,d\}$ satisfy the following conditions as $x\to\infty$
\begin{eqnarray*}\label{6b}\frac{P\{Y_{n,i}>x\}}{x^{-k_1}\ell_1(x)}&\rightarrow& c_i, ~~i\in\{1,...,d\},
\end{eqnarray*}
for some non-negative numbers $c_i$,
\begin{eqnarray*}\label{6c}\frac{P\{Y_{n,i}>x, Y_{n,j}>x\}}{x^{-k_1}\ell_1(x)}&\rightarrow& 0, ~~ i\neq j, ~~i,j\in\{1,...,d\}.
\end{eqnarray*}
\item[(A3)] Assume that for each $n\ge 1$ there exists $i\in\{1,...,d\}$ such that 
\begin{eqnarray*}
\!\!\!\!\!\!\!\!&& 
P\{\max_{1\le j\le d, j\neq i}(z_{j}Y_{n,j})>x, z_{i}Y_{n,i}\le x\}
=o(P\{z_iY_{n,i}>x\}),~~x\to\infty
\end{eqnarray*}
holds.
\\
\item[(A4)] Assume that there exists $i\in\{1,...,d\}$ such that it holds
\begin{eqnarray}\label{11b}\!\!\!\!\!&& 
P\{\max_{1\le j\le d, j\neq i}(z_{j}M_n^{(j)})>u_n, z_{i}M_n^{(i)}\le u_n\}
=o(1),~~n\to\infty.
\end{eqnarray}
\end{enumerate}
Let us denote $M_n^{(i)}= \max\{Y_{1,i}, Y_{2,i},...,Y_{n,i}\}, ~i\in\{1,..,l_n\}$. 
\begin{theorem}\label{T3} (Markovich, 2022)
Let the sets of slowly varying functions $\{\tilde{\ell}_n(x)\}_{n\geq 1}$ in (\ref{15a}) and $\{\ell_i(x)\}_{i\ge 1}$ in (\ref{11a}) satisfy the condition (\ref{uniform}), and (\ref{4a}), (\ref{27}), (\ref{chi}) hold.
Assume that $d$ and $\{Y_{n,i}\}$ are  independent.
\begin{enumerate}
\item [(i)]
Let $d$ be a bounded discrete r.v.  such that
$1<d<d_n = \min(C, l_n)$, 
$C>1$ holds.
\begin{enumerate}
\item
    If (A1) or (A2) for any $d\in\{2,3,...,\lfloor d_n-1\rfloor\}$
holds and
$N_n$ and $\{Y_{n,i}\}$ are  independent, then $Y_{n}(z,N_n)$ and $Y^*_{n}(z,N_n)$ have the  tail index $k_1$. If, instead of (A1) and (A2), (A3) holds, then $Y_n^*(z,N_n)$ has the same tail index.
\item If (A4) where in (\ref{11b}) $d$ is replaced by $\lfloor d_n-1\rfloor$
holds, then $Y_n^*(z,N_n)$ has the  extremal index $\theta_i$. If, in addition, (A1) (or (A2)) for any 
$d\in\{2,3,...,\lfloor d_n-1\rfloor\}$ holds, then $Y_n(z,N_n)$ has the same extremal index.
\end{enumerate}
\item [(ii)] Suppose that $d>1$ is a bounded discrete r.v. equal to a positive integer  a.s.. Then all statements of Item (i) are fulfilled.
\end{enumerate}
\end{theorem}
 It follows by Example 2 in Markovich (2022) that (\ref{11b}) is valid for all $d$ "column" sequences such that
 \begin{eqnarray}\label{5}
 &&M_n^{(1)}\le M_n^{(2)}\le M_n^{(3)}\le...\le M_n^{(d)}\end{eqnarray}
 holds. Theorem \ref{T3} means that if there are a random number $d$ of "column" sequences with a minimum tail index, then $Y^*_{n}(z,N_n)$ has the extremal index $\theta_i$ of the $i$th "column" satisfying (\ref{11b}),
 $1\le i\le d$. If the latter "column" sequences are independent (see, (A1)) or weakly dependent (see, (A2)), then $Y_{n}(z,N_n)$ has the same extremal index. For random networks this implies that the extremal index of the MLMs of the newly appended nodes is equal to the extremal index of the community with the minimum tail index that has a largest maximum PR among $d$ dominating communities.  The extremal index of the PRs of newly appended nodes is the same if the dominating communities satisfy (A1) or (A2). Since the communities are not enumerated, their maxima can be reordered as (\ref{5}). The statements of Theorem \ref{T3} are asymptotic. Thus, the approximation can be applied for sufficiently large size communities.
\\
In case of different pair-wise dependency among elements of the $d$ "column" series with the minimum tail index, the extremal index of the maxima $Y^*_{n}(z,N_n)$ and sums $Y_{n}(z,N_n)$ may not exist due to a non-stationarity of these sequences.
\begin{remark}\label{Rem1} The theorems in Markovich and Rodionov (2020a), Markovich (2022) are valid if there are non-zero 
elements in each row corresponding to the "column" sequences $\{Y_{n,i}: n\ge 1\}$ with minimum tail index. 
If the most heavy-tailed column is unique and at least one element in the latter sequence is equal to zero, the sequences of the sums and maxima of the "row" elements are not stationary. This feature plays a role for graphs.
\end{remark}
\begin{corollary}\label{Cor2} The statements of Theorems \ref{T3} 
remain valid if the tail indices $\{k_{n,i}\}$ of the elements
in the "columns" $\{Y_{n,i}: n\ge 1\}$  are different,
apart of those
columns with the minimum tail
index.
\end{corollary}
Corollary \ref{Cor2} is very important for practice. It implies that  the columns with non-minimum tail indices may be non-stationary distributed and hence, their extremal indices may not exist. Its proof follows from the proofs of Theorem 3 in Markovich and Rodionov (2020a) and Theorems 3 and 4 in Markovich (2022). 
The columns with the minimum tail index impact on the distribution and dependence structure of the sequence of sums and maxima over rows.

\section{Dependence structures in graphs}\label{Sec5.3}
\par We have to investigate the dependence of PRs of two  communities.
One of the approaches 
is to consider Pearson's correlation of two r.v.s belonging to two graphs. Each r.v. shows whether there is an edge between two nodes in a graph or not. Each  edge may be sampled iid  from a Bernoulli distribution with some
parameter $p$ (Xiong et al., 2020). A distance correlation is an extension of  Pearson's correlation both  to linear and nonlinear associations between two r.v.s or random vectors (Shen et al., 2020). It takes values in $[0,1]$. The distance correlation equal to zero does imply independence. Since nodes can be enumerated arbitrarily, the distance correlation has to be combined with a permutation test to check the dependence hypothesis. 
The distance correlation is calculated first for an original pair of vectors. It is compared with those ones calculated by shuffles of these vectors.
\par
In contrast to Shen et al. (2020), in our setting pairs of observations relating to two stationary distributed communities can be dependent and not necessarily identically distributed. One can use the distance correlation and the permutation test with regard to the row-column pairs. The $p-$value of the permutation test is the proportion of the number of the correlation measures from the samples with permutated pairs that are larger than the distance 
correlation that was calculated from the original data. 
\par To measure dependencies in heavy tailed graph data using statistical inference
for multivariate regular variation one can apply the polar coordinate transform to the examined random vectors $\{X_i\}$ and $\{Y_i\}$, $i=1,...,n$ (Resnick \& St\v{a}ric\v{a}, 1999; Volkovich et al., 2008; Wan et al., 2020). One can estimate the empirical distribution function (edf)  of the angular coordinates for the $k$ largest values of the  radial coordinate. The total dependence (or total independence) corresponds to the concentration of the edf to $\pi/4$ (or, to $0$ or $\pi/2$). A Starica plot  can be used   to find a suitable value of $k$.
\section{Preferential attachment}\label{Attach}\label{PA}
Let us consider a network growth where each node is attached to a small seed network at a unit time. The well-known  tool is a linear PA. 
A node $i$ can be attached randomly to existing nodes according to a probability $P_{PA}(i)=d_i/\sum_{s=1}^N d_s$ proportional to its degree $d_i$, or the number of its neighbors, where $N$ is the number of nodes. Nodes $i$ and $j$ may be connected with probability $d_id_j/\sum_{s=1}^N d_s$ (Norros \& Reittu, 2006). A kind of PA with a Poisson random number of new edges to the new vertex is proposed in Norros and Reittu (2006). The PA provides the "rich-get-richer" mechanism since earlier appearing nodes may increase their numbers of edges longer. This property leads to a power-law degree distribution $P(i)\sim i^{-(1+\alpha)}$ of node degrees 
(Newman, 2018; Wan et al., 2020). In Wan et al. (2020) it is derived that the  linear and superstar linear PA models on directed graphs lead to networks with power-law distributed in- and out-degrees.
\par
The $\alpha-$, $\beta-$ and $\gamma-$schemes of the linear PA  provide proportions of new nodes with incoming (outgoing) links to (from) existing nodes (scheme $\alpha-$ ($\gamma-$)) or  the directed edges between pairs of existing nodes ($\beta-$scheme) (Samorodnitsky et al., 2016; Wan et al., 2020). Let $I_n(v)$ and $O_n(v)$ be in- and out-degree of vertex $v\in V_n$ in a graph $G_n$, $n$  and $N(n-1)$ denote the numbers of edges and nodes  in $G_n$, respectively. Appending a new node $v$ to an existing graph $G_{n-1}$, one can select one of three scenarios by generating an iid sequence of trinomial r.v.s with cells marked $1,2,3$ with probabilities $\alpha, \beta, \gamma$. The probability to generate the edge $v\to w$ from $v$ to an existing node $w$ is given by
\begin{eqnarray}\label{al_sch}P\{\mbox{choose}~w\in V(n-1)\}&=&\frac{I_{n-1}(w)+\delta_{in}}{n-1+\delta_{in}N(n-1)}\end{eqnarray} by the $\alpha-$scheme;
between the existing nodes $v$ and $w$
\begin{eqnarray}\label{be_sch}P\{\mbox{choose}~(v,w)\}&=&\left(\frac{I_{n-1}(w)+\delta_{in}}{n-1+\delta_{in}N(n-1)}\right)\left(\frac{O_{n-1}(w)+\delta_{out}}{n-1+\delta_{out}N(n-1)}\right)\end{eqnarray} by the $\beta-$scheme; from the existing node $w$ to $v$
\begin{eqnarray}\label{ga_sch}P\{\mbox{choose}~w\in V(n-1)\}&=&\frac{O_{n-1}(w)+\delta_{out}}{n-1+\delta_{out}N(n-1)}\end{eqnarray} by the $\gamma-$scheme,
where $\delta_{in}$ and $\delta_{out}$ are parameters of the PA method. The latter may be estimated by the semi-parametric extreme value method  based on the maximum-likelihood method (Wan et al., 2020).

\section{Proof of Theorem \ref{Prop1}}\label{Sec5}
\begin{proof}
Part $(i)$ follows by Theorem 4 in Markovich \& Rodionov (2020). 
 We start the induction with $m=1$.
 The columns of matrix $A^{(1)}$ of the first iteration are obtained using submatrices of  $A^{(0)}$ (\ref{4c}) for different $j$
 by 
 recursions (\ref{1b}) and (\ref{1e}). 
 By the latter theorem 
 the $j$th columns  $\{Y^{(1)}_{i,j}\}$ and  $\{X^{(1)}_{i,j}\}$ have the same tail indices  $\{k^{(0)}_j\}$
 and the same  extremal indices $\{\theta^{(0)}_j\}$, $j\ge 1$.
 Getting matrices $A^{(2)}, A^{(3)},...$ for the next iterations both for sums and maxima similarly we obtain the same pairs
 of indices $(k^{(0)}_j,\theta^{(0)}_j)$ for their $j$th columns. 
 \\
 Part $(ii), (a)$.  The independence condition (A1) for $A^{(0)}$ is valid since communities (the "column" series of $A^{(0)}$) are nearly disconnected. The condition (A2) follows by (A1). The random number of communities $N_n$ is evidently independent on the PRs of nodes within the communities. Then the first $d_1^{(0)}$ "column" series of matrix  $A^{(1)}$ have the tail index $k_1^{(0)}$ by Theorem \ref{T3}. The next $d_2^{(0)}$ columns have the tail index $k_2^{(0)}>k_1^{(0)}$, etc. The columns of $A^{(1)}$ have the same tail indices as $A^{(0)}$.
 \\
 Note that the "column" series $\{Y^{(m)}_{i,j}\}$ and  $\{X^{(m)}_{i,j}\}$ of $A^{(m)}$, $m\ge 1$ are dependent due to their definition as partial sums and maxima of row elements of $A^{(m-1)}$ by the "domino" principle. Hence, we have $Y^{(m)}_{i,1}\ge Y^{(m)}_{i,2}\ge ...\ge Y^{(m)}_{i,N_n}$ and $X^{(m)}_{i,1}\ge X^{(m)}_{i,2}\ge ...\ge X^{(m)}_{i,N_n}$.
Elements of matrices $A^{(m)}$, $m\ge 1$ may be represented by elements of $A^{(0)}$. Really, for $m=2$ we have
\begin{eqnarray}\label{3a}Y^{(2)}_{n,i}&=& \sum_{j=i}^{N_n}Y^{(1)}_{n,j}=\sum_{j=i}^{N_n}jY^{(0)}_{n,j},\qquad i\ge 1,
\end{eqnarray}
and similarly for  $X^{(2)}_{n,i}$. Considering weights $z_j=j$ in (\ref{3a}) as in (\ref{3}), we obtain by Theorem \ref{T3} that the first $d_1^{(0)}$ sequences $\{Y^{(2)}_{n,i}\}$ and $\{X^{(2)}_{n,i}\}$ have  the tail index $k_1^{(0)}$, the next $d_2^{(0)}$ ones - $k_2^{(0)}$, etc. in the same way as "column" series of matrix $A^{(1)}$. The same is valid for $A^{(m)}$ with $m>2$ by induction.
\\
Part $(ii), (b)$. 
In order to find the extremal indices of the  "column" series of  $A^{(1)}$, let us enumerate the $d_j^{(0)}$ columns (the communities) in an descending order of the PR maxima  of $A^{(0)}$ over columns, i.e. $M_n^{(1)}\ge M_n^{(2)}\ge ...\ge M_n^{(d_j^{(0)})}$ for each fixed value of r.v. $d_j^{(0)}\in\{2,3,...,\lfloor d_n-1\rfloor$ and $j\in\{1,2,...\}$. Then (A4) is fulfilled consistently for $i\in\{1,2,...,d_j^{(0)}-1\}$.  By Theorem  \ref{T3} the first $d_j^{(0)}$, $j\ge 1$  "column" series $\{X_{i,j}^{(1)}\}_{i\ge 1}$ of $A^{(1)}$
have the extremal indices $\theta_{d_{j-1}^{(0)}+1}^{(0)},..., \theta^{(0)}_{d_{j-1}^{(0)}+d_{j}^{(0)}}$ for any values of $d_j^{(0)}\in\{2,3,...,\lfloor d_n-1\rfloor\}$, $d_0^{(0)}=0$. Since, in addition, (A1) (or (A2)) for $A^{(0)}$ and any
$d_j^{(0)}\in\{2,3,...,\lfloor d_n-1\rfloor\}$ holds, then $\{Y_{i,j}^{(1)}\}_{i\ge 1}$  have the same extremal indices.
Since elements of $A^{(m)}$ may be represented as weighted sums or maxima of elements of $A^{(0)}$, the extremal indices of the  "column" series of  $A^{(m)}$, $m\ge 1$ are the same as ones of $A^{(0)}$.
  \end{proof}